\documentclass[12pt, reqno]{amsart}
\usepackage[english]{babel}
\usepackage{amsmath,amscd,amssymb,latexsym,graphicx}
\usepackage[all]{xy}

\usepackage{color}

\newtheorem{theorem}{Theorem}[section]

\newtheorem{proposition}[theorem]{Proposition}

\newtheorem{corollary}[theorem]{Corollary}

\newtheorem{lemma}[theorem]{Lemma}

\theoremstyle{definition}
\newtheorem{definition}[theorem]{Definition}

\theoremstyle{definition}
\newtheorem{examples}[theorem]{Examples}

\theoremstyle{definition}
\newtheorem{remark}[theorem]{Remark}

\theoremstyle{definition}
\newtheorem{exercise}[theorem]{Exercise}

\numberwithin{equation}{section}

\newcommand{\cA}{\mathcal{A}}

\newcommand{\cC}{\mathcal{C}}

\newcommand{\cE}{\mathcal{E}}
\newcommand{\cF}{\mathcal{F}}
\newcommand{\cG}{\mathcal{G}}
\newcommand{\cH}{\mathcal{H}}

\newcommand{\cK}{\mathcal{K}}
\newcommand{\cL}{\mathcal{L}}
\newcommand{\cM}{\mathcal{M}}

\newcommand{\cP}{\mathcal{P}}

\newcommand{\cR}{\mathcal{R}}

\newcommand{\cT}{\mathcal{T}}
\newcommand{\cU}{\mathcal{U}}

\newcommand{\cX}{\mathcal{X}}

\newcommand{\fS}{ \mathsf{S}}

\newcommand{\NN}{\mathbb N}
\newcommand{\ZZ}{\mathbb Z}

\newcommand{\RR}{\mathbb R}
\newcommand{\CC}{\mathbb C}

\newcommand{\im}{\mathop{\mathrm{Im}}\nolimits}

\newcommand{\id}{\mathop{\mathrm{Id}}\nolimits}

\newcommand{\ot}{\otimes}
\newcommand{\mor}{\mathop{\mathrm{Mor}}\nolimits}
\newcommand{\pb}[1]{{}^*#1^*}

\newcommand{\smX}{X^{\circ}}

\newcommand{\ds}{\displaystyle}
\newcommand\smx{X^{\circ}}

\newcommand\smE{E^{\circ}}

\newcommand\smpi{\pi^{\circ}}

\title{Index theory and groupoids} 

\author{Claire Debord and Jean-Marie Lescure } 
 
\date{\today} 
 
\begin{document} 

\maketitle 

\begin{abstract} These lecture notes are mainly devoted to a proof
  using groupoids and $KK$-theory of Atiyah and Singer's index theorem on
  compact smooth manifolds. We first
  present an elementary introduction to groupoids, $C^*$-algebras,
  $KK$-theory and pseudodifferential calculus on groupoids. We then
  show how the point of view adopted here generalizes to the case of
  conical pseudo-manifolds.
\end{abstract}

\setcounter{tocdepth}{3}
\tableofcontents

\centerline{\bf \Large INTRODUCTION}

\bigskip During this course we intend to give the tools
involved in our approach of index theory for singular spaces. The
global framework adopted here is Noncommutative Geometry, with a
particular focus on groupoids,  $C^{*}$-algebras and bivariant $K$-theory.

The idea to use $C^{*}$-algebras to study {\sl spaces} may be
understood with the Gelfand theorem which asserts that Hausdorff
locally compact spaces are in one to one correspondance with
commutative $C^{*}$-algebras. 

A starting point in Noncommutative Geometry is then to think of
noncommutative $C^{*}$-algebras as corresponding to a wider class of
spaces, more singular than Hausdorff locally compact spaces.  

As a first consequence,  given a geometrical or topological object
which is badly behaved with respect to classical tools,
Noncommutative Geometry suggests to define a $C^{*}$-algebra
encoding relevant information carried by the original object. 

Refining this construction, one may try to define this
$C^{*}$-algebra as the $C^{*}$-algebra of a groupoid \cite{Re,Re1}. That is, one
can try to build a groupoid directly, encoding the original object
and regular enough to allow the construction of its $C^{*}$-algebra. 
In the ideal case where the groupoid is smooth,  one gets much more
than a $C^{*}$-algebra,  which only reflects topological properties:
the groupoid has a geometrical and analytical flavor enabling many
applications.  

An illuminating example is the study of the space of leaves of a
foliated manifold $(M, \cF)$ \cite{Co2,Co3,CoS}.  While this space $M/\cF$ is usually
very singular,  the holonomy groupoid of the foliation leads to a
$C^{*}$-algebra $C^{*}(M, \cF)$ replacing with great success the algebra of
continuous functions on the space $M/\cF$. Moreover,  the holonomy
groupoid is smooth and characterizes the original foliation. 

Once a $C^{*}$-algebra is built for the study of a given problem, one
can look for {\sl invariants} attached to it. For ordinary spaces,
basic invariants live in the homology or cohomology
of the space.  When dealing with $C^{*}$-algebras, the suitable
homology theory is $K$-theory, or better the $KK$-theory developped by
G. Kasparov \cite{Ka1,Ka2,ska} (when a smooth subalgebra of the
$C^{*}$-algebra is specified, which for instance is the case if a
smooth groupoid is available, one may also consider cyclic (co-)homology, but
this theory is beyond the scope of these notes). 

There is a fundamental theory which links the previous
ideas,  namely index theory. In the 60's, M. Atiyah and I. Singer
\cite{AS1} showed their famous index theorem. Roughly speaking, they
showed that,  given a closed manifold,  one can associate to any elliptic
operator an integer called the {\it index} which can be described in two
different ways: one purely analytic and the other one purely
topological. This result is stated with the help of $K$-theory of
spaces. Hence using the Swan-Serre theorem,  it can be formulated
with $K$-theory of (commutative) $C^{*}$-algebras. This point,  and the
fact that the index theorem can be proved in many ways using
$K$-theoretic methods,  leads to the attempt to generalize it to more
{\sl singular} situations where appropriate $C^{*}$-algebras are available.
In this a way, Noncommutative Geometry is a very general framework in which
one can try to state and prove index theorems. The case of
foliations illustrates this perfectly again: elliptic operators along
the leaves and equivariant with respect to 
the holonomy groupoid,  admit an analytical index living in the
$K$-theory of the $C^{*}$-algebra $C^{*}(M, \cF)$. Moreover this
index can also be described in a topological way and this is the
contents of the index theorem for foliations of A. Connes and
G. Skandalis \cite{CoS}. 

A. Connes \cite{Co0} also observed the important role played by groupoids in the
definition of the index map: in both cases of closed manifolds and
foliations, the analytical index map  can be described with the use of a
groupoid, namely a {\it deformation groupoid}. This approach has been
extented by the authors and V. Nistor \cite{DLN} who showed that the topological
index of Atiyah-Singer can also be described using deformation
groupoids. This leads to a {\sl geometrical proof} of the index theorem of
Atiyah-Singer; moreover this proof easily apply to a class of singular spaces (namely,
pseudomanifolds with isolated singularities).

The contents of this serie of lectures are
divided into three parts.  Let us briefly describe them:

\smallskip\noindent{\bf Part 1: Groupoids and their $C^{*}$-algebras.}

\smallskip As mentioned earlier, the first problem in the study of a singular geometrical situation is to associate
to it a mathematical object which carries the information one wants to study and which is regular enough to be analyzed in
a reasonable way.   In noncommutative geometry, answering this
question amounts to looking for a good {\it groupoid} and constructing its {\it $C^*$-algebra}. These points will
be the subject of the first two sections. 

\smallskip\noindent {\bf  Part 2: $KK$-theory.}

\smallskip Once the situation is desingularized,  say trough the
construction of a groupoid and its $C^{*}$-algebra, one may look for 
invariants which capture the basic properties. Roughly speaking, the {\it $KK$-theory} groups are convenient groups of
invariants for $C^*$-algebras and $KK$-theory comes with
powerful tools to carry out computations. Kasparov's bivariant $K$-theory will be the main topic of
sections 3 to \ref{finKK}.

\smallskip\noindent  {\bf  Part 3: Index theorems.}

\smallskip We first briefly explain in section \ref{pseudo} the pseudo-differential calculus on
groupoids. Then in Section \ref{debutInd},  we give a geometrical
proof of the Atiyah Singer index theorem for closed manifolds
using the language of groupoids and $KK$-theory.
Finally we show in the last section
how these results can be extended to {\it conical pseudo-manifolds}.

\smallskip \noindent {\bf Prerequists. } The reader interested in this
course should have background on several domains. Familiarity with
basic differential geometry (manifolds, tangent spaces) is
needed. The notions of fibre bundle, of $K$-theory for locally
compact spaces should  be known. Basic functional analysis like  analysis
of linear operators on Hilbert spaces should be familiar. The
knowledge of pseudodifferential calculus (basic
definitions, ellipticity) is necessary. Altough it is not absolutely
necessary, some familiarity with $C^{*}$-algebras is preferable.   

\smallskip \noindent {\bf Acknowledgments } We would like to thank
Georges Skandalis who allowed us to use several of his works to write
up this
course, in particular the manuscript of one of his courses
\cite{ska-DEA,ska}. We would like to warmly thank Jorge Plazas
for having typewritten a part of this course during the summer
school and J\'er\^ome Chabert who carefully read these notes and
corrected several mistakes. We are grateful to all the organizers for their kind
invitation to the extremely stimulating summer school held
at Villa de Leyva in July 2007 and we particulary thank Sylvie Paycha
both as an organizer and for her valuable comments on this document.

\newpage \part*{GROUPOIDS AND THEIR $C^*$-ALGEBRAS}

\medskip This first part will be devoted to the notion of groupoid,
specifically that of 
differentiable groupoid. We provide definitions and consider standard
examples. The interested reader may look for example at
\cite{Mc2,CW}. We then recall the definition of
$C^*$-algebras and see how one can associate a $C^*$-algebra to a groupoid.
The theory of $C^*$-algebra of groupoid was initiated by Jean Renault
\cite{Re}. A really good reference for the construction of groupoid
$C^*$-algebras is \cite{KhSk} from which the end of this section is inspired.

\section{Groupoids} 
\subsection{Definitions and basic examples of groupoids}

\begin{definition} Let $G$ and $G^{(0)}$ be two sets.  A structure of
  {\it groupoid} on $G$ over $G^{(0)}$ is given by the following homomorphisms:\begin{itemize}
\item[$\circ$] An injective map $u:G^{(0)}\rightarrow G$. The map $u$
  is called the {\it unit map}. We often identify $G^{(0)}$ with its image
  in $G$. The set $G^{(0)}$ is called the {\it  set of units} of the groupoid.
\item[$\circ$] Two surjective maps: $r,s: G\rightarrow G^{(0)}$,
  which are respectively the {\it range} and {\it source} map. They are equal to
  identity on the space of units.
\item[$\circ$] An involution: $$\begin{array}{cccc} i: & G & \rightarrow & G
  \\ & \gamma & \mapsto & \gamma^{-1} \end{array}$$ called the {\it inverse}
map. It satisfies: $s\circ i=r$.
\item[$\circ$] A map $$\begin{array}{cccc} p: & G^{(2)} & \rightarrow & G
  \\ & (\gamma_1,\gamma_2) & \mapsto & \gamma_1\cdot \gamma_2 \end{array}$$
called the {\it product}, where the set 
$$G^{(2)}:=\{(\gamma_1,\gamma_2)\in G\times G \ \vert \
s(\gamma_1)=r(\gamma_2)\}$$ is the set of {\it composable pairs}. Moreover for $(\gamma_1,\gamma_2)\in
G^{(2)}$ we have $r(\gamma_1\cdot \gamma_2)=r(\gamma_1)$ and $s(\gamma_1\cdot \gamma_2)=s(\gamma_2)$.
\end{itemize}
The following properties must be fulfilled:
\begin{itemize}
\item[$\circ$] The product is associative: for any $\gamma_1,\
  \gamma_2,\ \gamma_3$ in $G$ such that $s(\gamma_1)=r(\gamma_2)$ and
  $s(\gamma_2)=r(\gamma_3)$ the following equality
  holds $$(\gamma_1\cdot \gamma_2)\cdot \gamma_3= \gamma_1\cdot
  (\gamma_2\cdot \gamma_3)\ .$$
\item[$\circ$] For any $\gamma$ in $G$: $r(\gamma)\cdot
  \gamma=\gamma\cdot s(\gamma)=\gamma$ and $\gamma\cdot
  \gamma^{-1}=r(\gamma)$.
\end{itemize}
A groupoid structure on $G$ over $G^{(0)}$ is usually denoted by
$G\rightrightarrows G^{(0)}$,  where the arrows stand for the source
and target maps. 
\end{definition}
 
\noindent  We will often use the following notations: $$G_A:=
s^{-1}(A)\ ,\ G^B=r^{-1}(B) \ \mbox{ and } G_A^B=G_A\cap G^B \ .$$
If $x$ belongs to $G^{(0)}$, the {\it $s$-fiber} (resp. {\it
  $r$-fiber}) of $G$ over $x$ is $G_x=s^{-1}(x)$ (resp. $G^x=r^{-1}(x)$).

\medskip \noindent The groupoid is {\it topological} when
$G$ and $G^{(0)}$ are topological spaces with $G^{(0)}$ Hausdorff, the
structural homomorphisms are continuous and $i$ is an
homeomorphism. We will often require that our topological groupoids be {\it
  locally compact}. This means that $G\rightrightarrows G^{(0)}$ is a
topological groupoid, such
that $G$ is second countable, each point $\gamma$ in $G$ has a compact (Hausdorff)
neighborhood, and the map $s$ is open. In this situation the map $r$
is open and the $s$-fibers of $G$ are Hausdorff. 

\smallskip \noindent The groupoid is {\it smooth} when
$G$ and $G^{(0)}$ are second countable smooth manifolds with $G^{(0)}$ Hausdorff, the
structural homomorphisms are smooth, $u$ is an embedding, $s$ is a
submersion and $i$ is a 
diffeomorphism.

\smallskip \noindent When $G\rightrightarrows G^{(0)}$ is  at least
topological, we say that $G$ is {\it $s$-connected} when for any $x\in
G^{(0)}$, the $s$-fiber of $G$ over $x$ is connected. The
$s$-connected component of a groupoid $G$ is $\cup_{x\in G^{(0)}} CG_x$
  where $CG_x$ is the connected component of the $s$-fiber $G_x$
  which contains the unit $u(x)$.

\medskip \noindent {\bf Examples}

\smallskip \noindent {\bf 1.}  A space $X$ is a groupoid over itself with
  $s=r=u=\mbox{Id}$.

\smallskip \noindent {\bf 2.} A group $G\rightrightarrows \{e\}$ is a groupoid over its unit
  $e$, with the usual product and inverse map.

\smallskip \noindent {\bf 3.} A group bundle : $\pi: E\rightarrow X$ is a groupoid
  $E\rightrightarrows X$ with $r=s=\pi$ and algebraic operations
  given by the  group structure of each fiber $E_x$, $x\in X$.

\smallskip \noindent {\bf 4.} If $\cR$ is an equivalence relation on a space $X$, then the
  graph of $\cR$: $$G_{\cR}:=\{(x,y)\in X\times X \ \vert \ x\cR y\}$$
  admits a structure of groupoid over $X$,  which is given by:
$$ u(x)=(x,x)\ , \ s(x,y)=y\ ,\ r(x,y)=x \ , \\
  (x,y)^{-1}=(y,x) \ ,\ (x,y)\cdot(y,z)=(x,z) \ $$ for
$x,\ y,\ z$ in $X$.

\smallskip \noindent When $x\cR y$ for any $x,\ y$ in $X$,
$G_{\cR}=X\times X \rightrightarrows X$ is called the {\it pair
  groupoid}.

\smallskip \noindent {\bf 5.} If $G$ is a group acting on a space $X$, the {\it groupoid of
  the action} is $G\times X \rightrightarrows X$ with the following
structural homomorphisms 
$$\begin{array}{cc} u(x)=(e,x)\ , \ s(g,x)=x\ ,\ r(g,x)=g\cdot x \ , \\
  (g,x)^{-1}=(g^{-1},g\cdot x) \ ,\ (h,g\cdot x)\cdot(g,x)=(hg,x) \
  ,\end{array}$$ for 
$x$ in $X$ and $g,\ h$ in $G$. 

\smallskip \noindent {\bf 6.} Let $X$ be a topological space the {\it
  homotopy groupoid} of 
$X$ is $$\Pi(X):=\{\bar{c} \ \vert \ c:[0,1] \rightarrow X \mbox{ a
  continuous path} \} \rightrightarrows X $$ where $\bar{c}$ denotes
the homotopy class (with fixed endpoints) of $c$. We let
$$\begin{array}{ccc} 
u(x)=\overline{c_x} \mbox{ where } c_x \mbox{ is the constant path
  equal to } 
x,\ s(\overline{c})=c(0),\ r(\overline{c})=c(1) \\
 \overline{c}^{-1}=\overline{c^{-1}} \mbox{ where } c^{-1}(t)=c(1-t), \\
  \overline{c_1}\cdot \overline{c_2}=\overline{c_1\cdot c_2} \mbox{ where } c_1\cdot
  c_2(t)=c_2(2t) \mbox{ for } t\in[0,\frac{1}{2}] \mbox{ and } c_1\cdot
  c_2(t)=c_1(2t-1) \mbox{ for } t\in[\frac{1}{2},1] \ . 
\end{array}$$
When $X$ is a smooth manifold of dimension $n$, $\Pi(X)$ is naturally endowed with a
smooth structure (of dimension $2n$). A neighborhood of $\bar{c}$
is of the form $\{\bar{c_1}\bar{c}\bar{c_0} \ \vert \ c_1(0)=c(1),\
c(0)=c_0(1),\ Im c_i \subset U_i \ i=0,1\}$ where $U_i$ is a given
neighborhood of $c(i)$ in $X$.

\subsection{Homomorphisms and Morita equivalences} \

\smallskip \noindent {\bf Homomorphisms}\\
Let $G\rightrightarrows G^{(0)}$ be a groupoid of source $s_G$ and range
$r_G$ and $H\rightrightarrows H^{(0)}$ be a groupoid of source $s_H$ and range
$r_H$. A groupoid {\it homomorphism} from $G$ to $H$ is given by two maps
: $$ f:G\rightarrow H \mbox{ and }  f^{(0)}: G^{(0)}\rightarrow
H^{(0)} $$ such that 
\begin{itemize}\item[$\circ$] $r_H\circ f = f^{(0)} \circ r_G$,
\item[$\circ$] $f(\gamma)^{-1}=f(\gamma^{-1})$ for any $\gamma\in G$,
\item[$\circ$] $f(\gamma_1 \cdot \gamma_2)=f(\gamma_1)\cdot
  f(\gamma_2)$ for $\gamma_1,\ \gamma_2$ in $G$ such that 
  $s_G(\gamma_1)=r_G(\gamma_2)$.
\end{itemize}

\smallskip \noindent We say that $f$ is a {\it homomorphism over}
$f^{(0)}$. When $G^{(0)}=H^{(0)}$ and $f^{(0)}=\mbox{Id}$ we say that $f$ is a
{\it homomorphism over the identity}. 

\smallskip \noindent The homomorphism $f$ is an {\it isomorphism} when the maps $f$,
$f^{(0)}$ are bijections and $f^{-1} 
: H\rightarrow G$ is a homomorphism over $(f^{(0)})^{-1}$.

\smallskip \noindent As usual, when dealing with topological
groupoids we require that $f$ to
be continuous and when dealing with smooth groupoids, that $f$ be smooth.

\medskip \noindent {\bf Morita equivalence}\\ In most situations, the right notion of ``isomorphism of locally compact groupoids'' is
the weaker notion of Morita equivalence.

\begin{definition} Two locally compact groupoids $G\rightrightarrows
  G^{(0)}$ and $H\rightrightarrows H^{(0)}$ are {\it Morita
    equivalent} if there exists a locally compact groupoid
  $P\rightrightarrows G^{(0)} \sqcup H^{(0)}$ such that 

\begin{itemize} \item[$\circ$] the restrictions of $P$ over
$G^{(0)}$ and $H^{(0)}$ are respectively $G$ and $H$:
 $$ P_{G^{(0)}}^{G^{(0)}}=G \hbox{ and } P_{H^{(0)}}^{H^{(0)}}=H$$ 
\item[$\circ$] for any $\gamma \in P$ there exists $\eta$ in 
  $P_{G^{(0)}}^{H^{(0)}}\cup P_{H^{(0)}}^{G^{(0)}}$ such that 
  $(\gamma,\eta)$ is a composable pair (ie $s(\gamma)=r(\eta)$). 
\end{itemize} 
\end{definition}

\noindent {\bf Examples} {\bf 1.} Let $f:G \rightarrow H$ be an isomorphism of
locally compact groupoid.  Then the following groupoid defines a Morita
equivalence between $H$ and $G$: $$P=G\sqcup \tilde{G} \sqcup
\tilde{G}^{-1} \sqcup H \rightrightarrows G^{(0)} \sqcup H^{(0)}$$ where
with the obvious notations we have  $$\begin{array}{ccc}
G=\tilde{G}=\tilde{G}^{-1}\\
s_P= \left\{ \begin{array}{llll} s_G \mbox{ on } G \\ s_H\circ f \mbox{
      on } \tilde{G} \\ r_G \mbox{ on }
\tilde{G}^{-1} \\ s_H \mbox{ on } H \end{array} \right. ,\ 
r_P=\left\{ \begin{array}{lll} r_G \mbox{ on } G \sqcup \tilde{G} \\ s_H\circ f \mbox{
      on } \tilde{G}^{-1} \\ r_H \mbox{ on } H \end{array} \right. , \
u_P=\left\{ \begin{array}{cc} u_G \mbox{ on } G^{(0)} \\ u_H \mbox{ on }
    H^{(0)}\end{array} \right. \\
i_P(\gamma)=\left\{ \begin{array}{llll} i_G(\gamma) \mbox{ on } G \\
    i_H(\gamma) \mbox{ on } H \\ \gamma \in \tilde{G}^{-1} \mbox{ on }
    \tilde{G} \\ \gamma \in \tilde{G} \mbox{ on }
    \tilde{G}^{-1} \end{array} \right. ,
p_P(\gamma_1,\gamma_2)=\left\{ \begin{array}{lllllll}
    p_G(\gamma_1,\gamma_2) \mbox{ on } G^{(2)} \\
  p_H(\gamma_1,\gamma_2) \mbox{ on } H^{(2)} \\
p_G(\gamma_1,\gamma_2)\in \tilde{G} \mbox{ for } \gamma_1 \in G,\
\gamma_2\in \tilde{G} \\ p_G(\gamma_1,f^{-1}(\gamma_2))\in \tilde{G} \mbox{ for } \gamma_1 \in \tilde{G},\
\gamma_2\in H \\ p_G(\gamma_1,\gamma_2)\in {G} \mbox{ for } \gamma_1 \in \tilde{G},\
\gamma_2\in \tilde{G}^{-1} 
\\ f\circ p_G(\gamma_1,\gamma_2)\in H \mbox{ for } \gamma_1 \in \tilde{G},\
\gamma_2\in \tilde{G}^{-1} \end{array} \right.
\end{array}$$
{\bf 2.} Suppose that $G\rightrightarrows G^{(0)}$ is a locally
compact groupoid and $\varphi: X \rightarrow G^{(0)}$ is an open
surjective map, where $X$ is a locally compact space. The {\it pull
  back groupoid} is the groupoid: $$ 
\pb{\varphi}(G) \rightrightarrows X$$ where
$$\pb{\varphi}(G)=\{ (x,\gamma,y) \in X\times G\times X
\ \vert \ \varphi(x)=r(\gamma) \mbox{ and } \varphi(y)=s(\gamma) \}$$
with $s(x,\gamma,y)=y$, $r(x,\gamma,y)=x$, $(x,\gamma_1,y)\cdot
(y,\gamma_2,z)=(x,\gamma_1\cdot \gamma_2,z)$ and
$(x,\gamma,y)^{-1}=(y,\gamma^{-1},x)$. \\
One can show that this endows $\pb{\varphi}(G)$ with a structure of locally compact groupoid.
Moreover the groupoids $G$ and $\pb{\varphi}(G)$ are Morita
equivalent,  but not isomorphic in general.

\smallskip \noindent To prove this last point, one can put a structure of locally
compact groupoid on $P=G\sqcup X\times_r G \sqcup G\times_s X \sqcup
\pb{\varphi}(G)$ over $X\sqcup G^{(0)}$ where $X\times_r
G=\{(x,\gamma)\in X\times G \ \vert \ \varphi(x)=r(\gamma)\}$ and 
$G\times_s X=\{(\gamma,x)\in G\times X \ \vert \ \varphi(x)=s(\gamma)\}$.

\subsection{The orbits of a groupoid} \ 

\noindent Suppose that $G\rightrightarrows G^{(0)}$ is a groupoid of source $s$ and
range $r$.
\begin{definition} The {\it orbit} of $G$
  passing trough $x$ is the following subset of $G^{(0)}$:
  $$Or_x=r(G_x)=s(G^x)\ .$$
We let $G^{(0)}/G$ or $Or(G)$ be the {\it space of orbits}.

\smallskip \noindent The {\it isotropy group} of $G$ at $x$ is $G_x^x$, 
which is naturally endowed with a group structure with $x$ as
unit. Notice that multiplication induces a free left (resp. right) action of $G_x^x$
on $G^x$ (resp. $G_x$). Moreover the orbits space of this action is
precisely $Or_x$ and the restriction $s:G^x\rightarrow Or_x$ is the quotient map.
\end{definition}

\noindent {\bf Examples and remarks} {\bf 1.} In Example $4.$
above,  the orbits of $G_{\cR}$ 
correspond exactly to the orbits of the equivalence relation
$\cR$. In Example $5.$ above the orbits of the groupoid of the action
are the orbits of the action.\\
{\bf 2.} The second assertion in the definition of Morita equivalence
 precisely means that both $G^{(0)}$ and $H^{(0)}$ meet all the orbits of
$P$. Moreover one can show that the map $$\begin{array}{ccc} Or(G) &
  \rightarrow & Or(H) \\ Or(G)_x & \mapsto & Or(P)_x\cap H^{(0)}
\end{array}$$ is a bijection. In other word, when two groupoids are
Morita equivalent, they have the same orbits space. 

\medskip \noindent Groupoids are often used in Noncommutative Geometry
for the study of geometrical singular situations. 
In many geometrical situations, the topological space which arises
is strongly non Hausdorff and the standard tools do not
apply. Nevertheless, it is sometimes possible to associate to such
a space  $X$  a relevant $C^*$-algebra as a substitute for
$C_0(X)$. Usually we first associate a groupoid 
$G\rightrightarrows G^{(0)}$ such that its space of orbits $G^{(0)}/G$
is (equivalent to) $X$. If the groupoid is  regular enough (smooth for
example) then we can associate natural $C^*$-algebras to $G$.  This
point will be discussed later. \\ In other
words we desingularize a singular space by viewing it as coming from the
action of a nice groupoid on its space of units. To illustrate this
point let us consider two examples.

\subsection{Groupoids associated to a foliation}

Let $M$ be a smooth manifold.

\begin{definition} A (regular) smooth {\it foliation} $\cF$ on $M$ of
  dimension $p$ is a partition $\{F_i\}_I$ of $M$ where each $F_i$ is
  an immersed sub-manifold of dimension $p$ called a {\it
    leaf}. Moreover the manifold $M$ admits charts of the following
  type: $$\varphi: U \rightarrow \RR^p \times \RR^q$$ where $U$ is
  open in $M$ and such
  that for any  connected component $P$ of $F_i\cap U$ where $i\in I$,
  there is a $t\in \RR^q$ such that $\varphi(P)=\RR^p\times
  \{t\}$.

\smallskip \noindent In this situation the {\it tangent space to the
  foliation},  $T\cF:=\cup_{I}TF_i$,  is a sub-bundle of $TM$ stable under Lie bracket.

\smallskip \noindent The {\it space of leaves} $M/\cF$ is the quotient of $M$ by the
equivalence relation: being on the same leaf.
 \end{definition}

\noindent {\bf A typical example.} Take $M=P\times T$ where $P$ and $T$ are
connected smooth manifolds with the partition into leaves given by
$\{P\times \{t\}\}_{t\in T}$. Every foliation is locally of this
type.

\medskip \noindent  The space of leaves of a foliation is often difficult to
study,  as it appears in the following two examples:

\smallskip \noindent {\bf Examples} {\bf 1.} Let $\tilde{\cF}_a$ be the foliation
on the plane $\RR^2$ by lines
$\{y=ax+t\}_{t\in \RR}$ where $a$ belongs to $\RR$. Take the torus
$T=\RR^2/\ZZ^2$ to be the quotient of $\RR^2$ by translations of
$\ZZ^2$. We denote by $\cF_a$ the foliation induced by $\tilde{\cF}_a$ on
$T$. When $a$ is rational the space of leaves is a circle but when
$a$ is irrational it is topologically equivalent to a point (ie: each
point is in any neighborhood of any other point). \\
{\bf 2.} Let $\CC \setminus \{(0)\}$ be foliated by: 
$$ \{S_t\}_{t\in
  ]0,1]} \cup \{D_t\}_{t\in ]0,2\pi] } $$ where $S_t=\{z\in \CC \ \vert \
  \vert z\vert =t\}$ is the circle of
  radius $t$ and $D_t=\{z=e^{i(x+t)+x} \ \vert \ x\in \RR^+_* \}$.
\begin{center} 
\includegraphics[width=3.5cm]{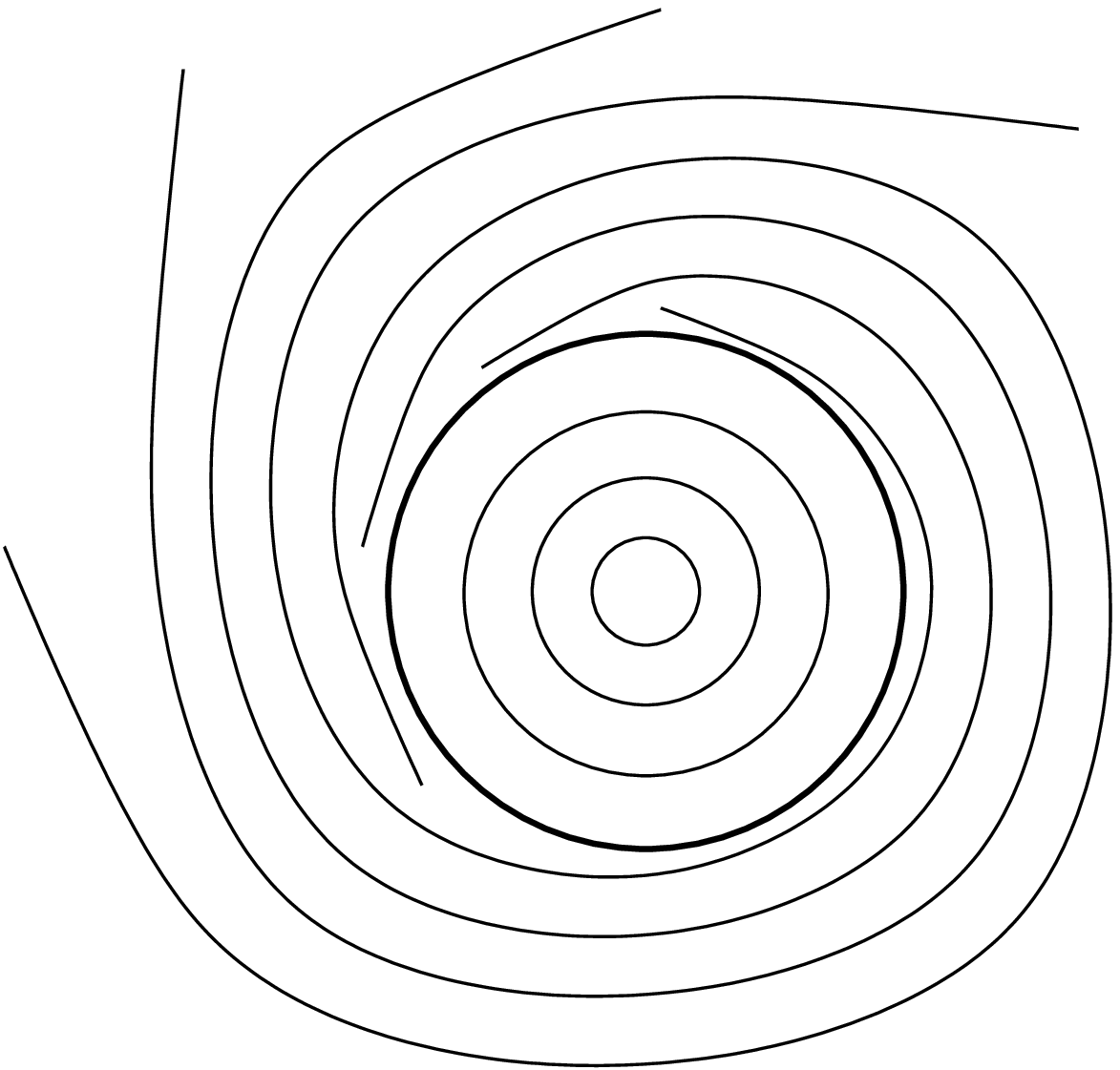} 
\end{center} 

\medskip \noindent  The {\it holonomy groupoid} is a smooth groupoid which
desingularizes the space of leaves of a foliation. Precisely, if $\cF$ is a smooth
foliation on a manifold $M$ its {\it holonomy groupoid} is the smallest
$s$-connected smooth groupoid 
$G\rightrightarrows M$ whose orbits are precisely the leaves of the
foliation. \\ Here,  {\it smallest} means that if $H\rightrightarrows M $
is another $s$-connected smooth groupoid whose orbits are the leaves of the
foliation then there is a surjective groupoid homomorphism
: $H\rightarrow G$ over identity.

\smallskip \noindent The first naive attempt to define such a groupoid is
to consider the graph of the equivalence relation being on the same
leaf. This does not work: you get a groupoid but it may be not
smooth. This fact can be observed on example 2. below.
Another idea consists in looking at the {\it homotopy groupoid}. Let
$\Pi(\cF)$ be the set of homotopy classes of smooth paths lying on
leaves of the foliation. It is naturally endowed with a groupoid
structure similarly to the homotopy groupoid of Section 1. Example
6. Such a groupoid can be naturally equipped with a smooth structure (of
dimension $2p+q$) and the holonomy groupoid is a quotient of this
homotopy groupoid. In particular, when the leaves have no homotopy,
the holonomy groupoid is the graph of the equivalence relation of being
in the same leaf.

\subsection{The noncommutative tangent space of a conical pseudomanifold}\label{ncgtangent}

It may happen that the underlying topological space which is under
study is a nice compact space which is ``almost'' smooth. This is the case of 
pseudo-manifolds \cite{GoMa,Mat,Ver}, for a review on the
subject see \cite{BHS,HW}. In such a situation we can desingularize
the tangent space \cite{DL2,Nous2}. Let us see how this works in the
case of a conical pseudomanifold with one singularity.

\medskip \noindent Let $M$ be an $m$-dimensional compact manifold with 
a compact boundary $L$. We attach to $L$ the cone  $cL= L\times [0,1] / L\times \{0\}$, using the 
obvious  map  $L\times\{1\} \to L\subset\partial M$. 
The new space $X=cL \cup M$ is a compact pseudomanifold 
with a singularity \cite{GoMa}. In general, there is no 
manifold structure around the vertex $c$ of the cone.

\smallskip \noindent 
We will use the following notations: $\smX=X\setminus\{c\}$ is the {\it regular part}, $ X^+$ denotes 
$M\setminus L = X\setminus cL$, $\overline{X_+}=M$ its closure in $X$ and 
$ X^-=L\times ]0,1 [$.  If $y$ is a point of the cylindrical part of
$X\setminus\{c\}$, we write $y=(y_L,k_y)$ where $y_L \in L$ and 
$k_y\in ]0,1] $ are the tangential and  radial 
coordinates. The map $y\to k_{y}$ is extended into a smooth defining
function for the boundary of $M$. In particular, $k^{-1}(1)=L=\partial M$ and $k(M) \subset[1,+\infty[$. 
\begin{center} 
\includegraphics[height=3cm]{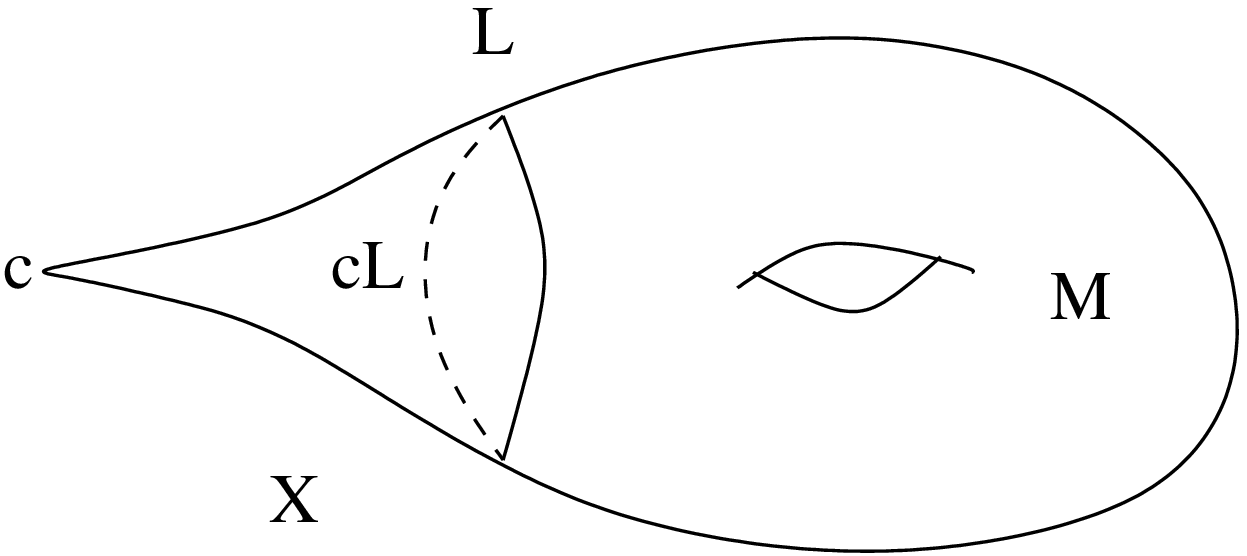} 
\end{center}

\smallskip \noindent Let us consider $T\overline{X^+}$, the restriction to $\overline{X^+}$ 
of the tangent bundle of $\smX$. As a $\cC^\infty$ vector bundle, it is a 
smooth groupoid with unit space $\overline{X^+}$. We define the 
groupoid $T^\fS X$ as the disjoint union: 
$$T^{\fS}X=X^-\times X^-\ \cup \ 
T\overline{X^+} \overset{s}{\underset{r}{\rightrightarrows}} \smX,$$ 
where  $X^-\times X^-\rightrightarrows X^-$ is the pair groupoid. 
 
\smallskip \noindent In order to endow $T^{\fS}X$ with a smooth structure, 
compatible with the usual smooth structure on $X^-\times X^-$ and 
on $T\overline{X^+}$, we have to take care of what happens around points of 
$T\overline{X^+}|_{\partial \overline{X^+}}$. \\ Let $\tau$ be a 
smooth positive function on $\RR$ such that
$\tau^{-1}(\{0\})=[1,+\infty[$. 
We let $\tilde{\tau}$ be the smooth 
map from $\smX$ to $\RR^{+}$ given by $\tilde{\tau}(y)=\tau\circ k(y)$. \\ 
Let $(U,\phi)$ be a local chart for 
$\smX$ around $z\in\partial \overline{X^+}$. Setting $U^-=U\cap X^-$ and 
$\overline{U^+}= U \cap \overline{X^+}$, we define a local chart of $T^{\fS}X$ by: 
$$ \begin{array}{cccc} \tilde{\phi}: & 
 U^-\times U^- \cup \ T \overline{U^+} & \longrightarrow &  \RR^{m}\times \RR^{m} \end{array}$$ 
\begin{equation}\label{local-charts} \tilde{\phi}(x,y)= 
(\phi(x),\frac{\phi(y)-\phi(x)}{\tilde{\tau}(x)}) \mbox{ if } (x,y)\in 
U^-\times U^- \mbox{ and} 
\end{equation} 
$$\tilde{\phi}(x,V)=(\phi(x),(\phi)_*(x,V)) \mbox{ elsewhere.}$$ 

\noindent We define in this way a structure of smooth  groupoid on $T^{\fS}X$. Note
that at the topological level, the space of orbits of $T^{\fS}X$ is 
equivalent to $X$: there is a canonical isomorphism 
between the algebras $C(X)$ and $C(\smX/T^{\fS}X)$. 

\medskip \noindent The smooth groupoid 
  $T^{\fS}X \rightrightarrows \smX$ is called the {\it noncommutative
    tangent space} of $X$.

\subsection{Lie Theory for smooth groupoids} 
Let us go into the more specific world of smooth groupoids. Similarly to
Lie groups which admit Lie algebras, any smooth groupoid has a {\it
  Lie algebroid} \cite{Pr4,Pr5}.

\begin{definition} A {\it Lie
  algebroid} $\cA =(p:\cA\rightarrow TM,[\ ,\ ]_{\cA})$ on a smooth
manifold $M$ is a vector bundle $\cA \rightarrow M$
equipped with a bracket $[\ ,\ ]_{\cA}:\Gamma(\cA)\times \Gamma(\cA)
\rightarrow \Gamma(\cA)$ on the module of sections of $\cA$ together
with a homomorphism of fiber bundle $p:\cA \rightarrow TM$ from $\cA$ to the
tangent bundle $TM$ of $M$ called the {\it anchor}, such that:
\begin{enumerate}
\item[i)] the bracket $[\ ,\ ]_{\cA}$ is $\RR$-bilinear, antisymmetric
  and satisfies the Jacobi identity,
\item[ii)] $[X,fY]_{\cA}=f[X,Y]_{\cA}+p(X)(f)Y$ for all $X,\ Y \in
\Gamma(\cA)$ and $f$ a smooth function of $M$, 
\item[iii)] $p([X,Y]_{\cA})=[p(X),p(Y)]$ for all
$X,\ Y \in \Gamma(\cA)$.
\end{enumerate}
\end{definition}

\smallskip \noindent
Each Lie groupoid admits a Lie algebroid. Let us recall this
construction.

\smallskip \noindent 
Let $G \overset{s}{\underset{r}\rightrightarrows}G^{(0)}$ be a Lie
groupoid. We denote by $T^sG$ the subbundle of $TG$ of $s$-vertical tangent
vectors. In other words, $T^sG$ is the kernel of the differential $Ts$ of
$s$.

\smallskip  \noindent
For any $\gamma$ in $G$ let $R_{\gamma}: G_{r(\gamma)} \rightarrow
G_{s(\gamma)}$ be the right multiplication by $\gamma$. A tangent vector
field $Z$ on $G$ is {\it right invariant} if it  satisfies:

\begin{itemize}
\item[--] $Z$ is $s$-vertical: $Ts(Z)=0$.
\item[--] For all $(\gamma_1,\gamma_2)$ in $G^{(2)}$, $Z(\gamma_1 \cdot
  \gamma_2)=TR_{\gamma_2}(Z(\gamma_1))$.
\end{itemize}

\noindent
Note that if $Z$ is a right invariant vector field and $h^t$ its flow
then for any $t$, the local diffeomorphism $h^t$ is a {\it local left
  translation} of $G$ that is $h^t(\gamma_1 \cdot
  \gamma_2)=h^t(\gamma_1)\cdot \gamma_2$ when it makes sense.

\medskip \noindent The Lie algebroid ${\cA}  G$ of $G$ is defined as follows:
\begin{itemize}
\item[--] The fiber bundle ${\cA}  G \rightarrow G^{(0)}$ is the restriction of
$T^sG$ to $G^{(0)}$. In other words: $\cA G=\cup_{x\in G^{(0)}} T_x
G_x$ is the union of the tangent spaces to the $s$-fiber at the corresponding unit.
\item[--] The {\it anchor} $p:{\cA}  G \rightarrow TG^{(0)}$
is the restriction of the differential $Tr$ of $r$
to ${\cA}  G$.
\item[--] If $Y:U \rightarrow {\cA}  G$ is a local section of ${\cA}  G$,
  where $U$ is an open subset of $G^{(0)}$, we define the local {\it
    right invariant vector field} $Z_Y$ {\it associated} with $Y$ by
  $$Z_Y(\gamma)=TR_{\gamma}(Y(r(\gamma))) \ \makebox{ for all } \
  \gamma \in G^U \ .$$

\noindent The Lie bracket is then defined by:
$$\begin{array}{cccc} [\ ,\ ]: & \Gamma({\cA}  G)\times \Gamma({\cA} 
  G) & \longrightarrow & \Gamma({\cA}  G) \\
 & (Y_1,Y_2) & \mapsto & [Z_{Y_1},Z_{Y_2}]_{G^{(0)}}
\end{array}$$
where $[Z_{Y_1},Z_{Y_2}]$ denotes the $s$-vertical vector field obtained
with the usual bracket and $[Z_{Y_1},Z_{Y_2}]_{G^{(0)}}$ is the
restriction of $[Z_{Y_1},Z_{Y_2}]$ to ${G^{(0)}}$.
\end{itemize}

\medskip \noindent {\bf Example} If $\Pi(\cF)$ is the homotopy
groupoid (or the holonomy groupoid) of a smooth foliation, its Lie
algebroid is the tangent space $T\cF$ to the foliation. The anchor
is  the inclusion. In particular the Lie algebroid of the pair
groupoid $M\times M$ on a smooth manifold $M$ is $TM$, the anchor
being the identity map.

\medskip \noindent Lie theory for groupoids is much trickier than
for groups. For a long time people thought that, as for Lie
algebras, every Lie algebroid integrates into a Lie groupoid \cite{Pr3}. In fact
this assertion, named {\it Lie's third theorem for Lie algebroids} is
false. This was pointed out by a counter example given by P. Molino
and R. Almeida in \cite{M0Al}. Since then, a lot of work has been
done around this problem. A few years ago 
M. Crainic and R.L. Fernandes \cite{CF} completely solved
this question by giving a  necessary and sufficient condition for the
integrability of Lie algebroids.

\subsection{Examples of groupoids involved in index theory}\label{groupoid-index-theory}
Index theory is a part of non commutative geometry where 
groupoids may play a crucial role. Index theory will be discussed later in
this course but we 
want to present here some of the groupoids which will arise.

\begin{definition}\label{def-deformation-groupoid} A smooth groupoid $G$ is called a {\it deformation
    groupoid} if: $$G=G_1\times \{0\} \cup G_2\times ]0,1]
  \rightrightarrows G^{(0)}=M\times [0,1]\ ,$$where $G_1$ and $G_2$
  are smooth groupoids with unit space $M$. That is, $G$ is obtained
  by gluing $G_2\times ]0,1]\rightrightarrows M\times ]0,1]$ which is
  the groupoid $G_2$ parametrized by $]0,1]$ with the groupoid
  $G_1\times \{0\} \rightrightarrows M\times \{0\}$.
\end{definition}

\noindent{\bf Example} Let $G$ be a smooth groupoid and let $\cA G$ be its Lie 
  algebroid. The {\it adiabatic groupoid} of $G$ \cite{Co0,MP,NWX} is a
  deformation of $G$ on its Lie algebroid: 
$$G_{ad}=\cA G\times \{0\}  \cup G\times ]0,1]   \rightrightarrows 
G^{(0)}\times [0,1],  $$ 
where One can put a natural smooth structure on $G_{ad}$.
Here, the vector 
bundle $\pi: \cA G\to G^{(0)}$ is considered as a groupoid in the 
obvious way.

\smallskip \noindent {\bf The tangent groupoid}

\noindent A special example of adiabatic groupoid is the {\it
  tangent groupoid} of A. Connes \cite{Co0}. Consider the pair groupoid
$M\times M$ on a smooth manifold $M$. We saw that its Lie
algebroid is $TM$. In this situation, the adiabatic groupoid is called
the 
{\it tangent groupoid} and is given by: $$\cG^t_M:= TM\times \{0\} \sqcup
M\times M \times ]0,1] \rightrightarrows M\times [0,1]\ .$$
The Lie algebroid is the bundle $\cA(\cG^t_M):= TM\times [0,1]
\rightarrow M\times [0,1]$ with anchor $p: (x,V,t)\in TM\times [0,1]
\mapsto (x,tV,t,0)\in TM\times T[0,1]$.\\ Choose a riemannian metric on
$M$.  The smooth structure on $\cG^t_M$ is such that the following map
: $$\begin{array}{ccc} \cU\subset TM\times [0,1] & \rightarrow
  &\cG^t_M \\ (x,V,t) & \mapsto & \left\{ \begin{array}{ll} (x,V,0)
      \mbox { if } t=0 \\ (x,exp_x(-tV),t) \mbox { elsewhere }
    \end{array} \right. \end{array}$$
is a smooth diffeomorphism on its range, where $\cU$ is an open
neighborhood of $TM\times \{0\}$.

\medskip \noindent The previous construction of the tangent groupoid
of a compact manifold generalizes to the case of conical manifold. When
$X$ is a conical manifold, its tangent groupoid is a deformation of
the pair groupoid over $\smX$ into the groupoid $T^{\fS}X$. This deformation has a nice 
description at the level of Lie algebroids. Indeed, with the notation
of \ref{ncgtangent}, the  Lie algebroid 
of $\cG^t_X$ is the (unique) Lie algebroid given 
by the fiber bundle $\cA \mathcal{G}^t_X=[ 0,1] \times \cA
(T^{\fS}X)=[ 0,1] \times T\smX \rightarrow [ 0,1] \times \smX$, 
with anchor map 
$$\begin{array}{cccc} p_{\mathcal{G}^t_X}: & \cA\mathcal{G}^t_X=[ 0,1] \times T\smX & \longrightarrow & T([ 0,1] 
\times \smX)=T[ 0,1]\times T\smX \\ & (\lambda,x,V) & \mapsto & 
(\lambda,0,x,(\tilde{\tau 
}(x)+\lambda)V)  \ . \end{array}$$ 
Such a Lie algebroid is  almost injective, thus it is integrable 
\cite{CF,Moa3}.\\ Moreover, it integrates into the {\it tangent
  groupoid} which is defined by: 
$$\mathcal{G}^t_X=\smX\times \smX\times ]0,1] \ \cup \ T^{\fS}X\times \{0\} 
{\rightrightarrows} \smX \times 
[0,1]. $$ 
Once again one can equip such a groupoid with a smooth structure
compatible with the usual one on each piece: $\smX\times \smX\times
]0,1]$ and $ T^{\fS}X\times \{0\} $ \cite{DL2}.

\smallskip \noindent {\bf The Thom groupoid}

\noindent Another important deformation groupoid for our
purpose is the {\it Thom groupoid} \cite{DLN}.\\ Let $\pi:E
\rightarrow X$ be a {\it conical vector} bundle. This means that $X$
is a conical manifolds (or a smooth manifold without vertices) and we
have a smooth vector bundle $\smpi: \smE \rightarrow \smX$ which
restriction to $X^-=L\times ]0,1[$ is equal to
$E_L\times ]0,1[$ where $E_L\rightarrow L$ is a smooth vector
bundle. If $E^+\rightarrow X^+$ denotes the bundle $\smE$ restricted
to $X^+$, then $E$ is the conical manifold: $E= cE_L\cup E^+$.

\smallskip \noindent When $X$ is a smooth manifold (with no conical
point), this boils down to the usual
notion of smooth vector bundle.

\medskip \noindent From the definition,
$\pi$  restricts to a smooth vector bundle map $\pi^{\circ}:
\smE \rightarrow \smx$. We let $\pi_{[0,1]}=\pi^{\circ} \times id
: \smE \times [0,1] \rightarrow \smx \times [0,1]$.

\smallskip \noindent We consider the tangent groupoids
$\cG^t_X \rightrightarrows \smx\times [0, 1]$ for $X$ and $\cG^t_E
\rightrightarrows \smE\times [0, 1]$ for $E$ equipped with a smooth structure
constructed using the same gluing function $\tau$ (in particular
$\tilde{\tau_X} \circ \pi = \tilde{\tau_E}$). We denote by
$\pb{\pi_{[0,1]}}(\cG^t_X)\rightrightarrows \smE \times [0,1]$ the
pull back of $\cG^t_X$ by $\pi_{[0,1]}$.

\smallskip \noindent 
We first associate to the conical vector bundle
$E$ a deformation groupoid $\cT^t_{E}$ from
$\pb{\pi_{[0,1]}}(\cG^t_X)$ to $\cG^t_E$.  More precisely, we
define:
\begin{equation*}
    \cT^t_{E}:= \cG^t_E \times \{0\} \sqcup
    \pb{\pi_{[0,1]}}(\cG^t_X)\times ]0,1]
    \rightrightarrows \smE \times [0,1] \times [0,1].
\end{equation*}
Once again, one can equip $\cT^t_E$ with a smooth structure \cite{DLN}
and the restriction of $\cT^t_E$ to
$\smE \times \{0\} \times [0,1]$ leads to a smooth groupoid:
\begin{equation*}
    \cH_{E} = T^{\fS}E\times \{0\} \sqcup
    \pb{\pi}(T^{\fS}X) \times ]0,1]
    \rightrightarrows \smE \times [0,1],
\end{equation*}
called a {\it Thom groupoid} associated to the conical vector bundle
$E$ over $X$.

\smallskip \noindent The following example explains what these constructions
become if there is no singularity.

\smallskip \noindent {\bf{Example}}
Suppose that $p:E\rightarrow M$ is a smooth vector bundle over the
smooth manifold $M$. Then we have the usual
tangent groupoids $\cG^t_E =
TE\times \{0\} \sqcup E \times E \times ]0,1] \rightrightarrows E
\times [0,1]$ and $\cG^t_M = TM \times \{0\} \sqcup M \times M
\times ]0,1] \rightrightarrows M \times [0,1]$. In this example the groupoid $\cT^t_{E}$ will
be given by
\begin{equation*}
    \cT^t_{E}= TE\times \{0 \}\times \{0 \}
    \sqcup \pb{p}(TM)\times \{0\} \times ]0,1]
    \sqcup E\times E \times ]0,1] \times [0,1]
    \rightrightarrows E\times [0,1]\times [0,1]
\end{equation*}
and is smooth. Similarly, the Thom groupoid will be given
by: $\cH_E:= TE\times \{0 \} \sqcup \pb{p}(TM) \times
]0,1]\rightrightarrows E\times [0,1]$.

\subsection{Haar systems}

A locally compact groupoid $G\rightrightarrows G^{(0)}$ can be viewed as a
family of locally compact spaces:
$$G_x=\{\gamma\in G \ \vert \ s(\gamma)=x\}$$ parametrized by $x\in
G^{(0)}$. Moreover, right translations act on these spaces. Precisely,
to any $\gamma \in G$ one associates the homeomorphism $$\begin{array}{cccc} R_\gamma: & G_y &
  \rightarrow & G_x \\ & \eta & \mapsto &
\eta \cdot \gamma  \ .\end{array}$$
This picture enables to define the right analogue of Haar measure on locally
compact groups to locally compact groupoids, namely {\it Haar systems}. The
following definition is due to J. Renault \cite{Re}.

\begin{definition} A {\it Haar system} on $G$ is a collection
  $\nu=\{\nu_x\}_{x\in G^{(0)}}$ of positive regular Borel measure on
  $G$ satisfying the following conditions:\begin{enumerate} \item
    {\it Support:} For every $x\in G^{(0)}$, the support of $\nu_x$ is
    contained in $G_x$.\\
\item {\it Invariance:} For any $\gamma \in G$, the right-translation
  operator $R_\gamma:  G_y \rightarrow  G_x$ is
  measure-preserving. That is, for all $f\in C_c(G)$: $$\int
  f(\eta)d\nu_y(\eta)=\int f(\eta \cdot \gamma) d\nu_x(\eta)\ .$$
\item {\it Continuity:} For all $f\in C_c(G)$, the map
  $$\begin{array}{ccc} G^{(0)} & \rightarrow & \CC \\ x & \mapsto &
    \int f(\gamma)d\nu_x(\gamma) \end{array} $$ is continuous.
\end{enumerate}
\end{definition}

\noindent In contrast to the case of locally compact groups, Haar systems on
groupoids may not exist. Moreover, when such a Haar system exists , it
may not be unique. In the special case of a smooth groupoid,  a Haar
system always exists \cite{Pat,Ra} and any two Haar systems $\{\nu_x\}$
and $\{\mu_x\}$ differ by a continuous and positive function $f$ on
$G^{(0)}$: $\nu_x=f(x)\mu_x$ for all $x\in G^{(0)}$.

\medskip \noindent {\bf Example:} When the source and range maps are local
homeomorphisms, a possible choice for $\nu_x$ is the counting measure
on $G_x$.


\section{$C^*$-algebras of groupoids}

\noindent This second part starts with the definition of a $C^*$-algebra
together with some results. Then we construct the
maximal and minimal $C^*$-algebras associated to a groupoid and compute
explicit examples.

\subsection{$C^*$-algebras - Basic definitions}
In this chapter we introduce the terminology and give some examples and
 properties of $C^*$-algebras. We refer the reader to \cite{Di,Pe,Ar} for
a more complete overview on this subject.

\begin{definition} A {\it $C^*$-algebra} $A$ is a complex Banach
  algebra 
with an involution $x\mapsto x^*$ such that: 
\begin{enumerate} \item $(\lambda x+ \mu y)^*=\bar{\lambda}
  x^*+\bar{\mu} y^*$ for $\lambda,\ \mu \in \CC$ and $x,\ y \in A$,\\
\item $(xy)^*=y^*x^*$ for $x,\ y \in A$, and\\
\item $\|x^*x \|=\|x\|^2$ for $x\in A$.
\end{enumerate} 
\end{definition}

\noindent Note that it follows from the definition that $*$ is
isometric.

\smallskip \noindent The element $x$ in $A$ is {\it self-adjoint} if
$x^*=x$, {\it normal} if $xx^*=x^*x$. When $1$ belongs to $A$, $x$ is
{\it unitary} if $xx^*=x^*x=1$.\\
Given two $C^{*}$-algebras $A, B$,  a homomorphism respecting the
involution is a called a $*$-homomorphism. 

\medskip \noindent {\bf Examples }{\bf 1.} Let $\cH$ be an Hilbert space. The algebra
$\cL(\cH)$ of all continuous linear transformations of $\cH$ is a
$C^*$-algebra. The involution of $\cL(\cH)$ is given by the usual
adjunction of bounded linear operators. \\ 
{\bf 2.}
Let $\cK(\cH)$ be the norm closure of finite rank operators on $\cH$. It
is the $C^*$-algebra of compact operators on $\cH$. \\
{\bf 3.} The algebra $M_n(\CC)$ is a $C^*$-algebra. It is a special
example of the previous kind, when $dim(\cH)=n$.\\
{\bf 4.} Let $X$ be a locally compact, Hausdorff, topological
space. The algebra $C_0(X)$  of continuous functions vanishing at
$\infty$, endowed with the supremum norm and the involution $f\mapsto \bar{f}$
is a commutative $C^*$-algebra. When $X$ is compact, $1$ belongs to
$C(X)=C_0(X)$.

\smallskip \noindent Conversely, the Gelfand'd theorem asserts that every commutative $C^*$-algebra $A$ is
isomorphic to $C_{0}(X)$ for some locally compact space $X$ (and it
is compact precisely when $A$ is unital). Precisely, a {\it character} $\cX$ of $A$ is a continuous
homomorphism of algebras $\cX:A\rightarrow \CC$. The set $X$ of characters
of $A$, called the {\it spectrum} of $A$, can be endowed with a
locally compact space topology.
The {\it Gelfand transform}
$\cF:A\rightarrow C_0(X)$ given by $\cF(x)(\cX)=\cX(x)$ is the desired
$*$-isomorphism.

\medskip \noindent Let $A$ be a $C^*$-algebra and $\cH$ a Hilbert space.
\begin{definition} A {\it $*$-representation} of $A$ in $\cH$ is a
  $*$-homomorphism $\pi:A\rightarrow \cL(\cH)$. The representation is
  {\it faithful} if $\pi$ is injective.
\end{definition}

\begin{theorem}(Gelfand-Naimark) If $A$ is a $C^*$-algebra, there exists
  a Hilbert space $\cH$ and a faithful representation
  $\pi:A\rightarrow \cL(\cH)$.
\end{theorem}
\noindent In other words, any $C^*$-algebra is isomorphic to a norm-closed
involutive subalgebra of $\cL(\cH)$. Moreover, when $A$ is separable, $\cH$ can
be taken to be the (unique up to isometry) separable Hilbert
space of infinite dimension.

\medskip \noindent{\bf Enveloping algebra}\\  Given a Banach
$*$-algebra $A$, consider the family $\pi_{\alpha}$ of all continuous
$*$-representations for $A$. The Hausdorff completion of $A$ for the
seminorm $\|x\|=sup_{\alpha}(\|\pi_{\alpha}(x)\|)$ is a $C^*$-algebra
called the {\it enveloping $C^*$-algebra} of $A$.

\medskip \noindent{\bf Units}\\ A $C^*$-algebra may or may not have a
unit, but it can always be embedded into a unital
$C^*$-algebra $\tilde{A}$: $$\tilde{A}:=\{x+\lambda \ \vert \ x\in
A,\ \lambda\in \CC\}$$ with the obvious product and involution. The
norm on $\tilde{A}$ is given for all $x\in \tilde{A}$ by: $\|x \|^{\sim}=Sup\{\|xy\|,\ y\in
A\ ;\ \|y\|=1\}$. On $A$ we have $\|\cdot \|=\|\cdot \|^{\sim}$. The algebra $A$
is a closed two sided ideal in $\tilde{A}$ and $\tilde{A}/A=\CC$.

\medskip \noindent {\bf Functional calculus}\\
Let $A$ be a $C^*$-algebra. If $x$ belongs to $A$, the {\it spectrum}
of $x$ in $A$ is the compact set: $$Sp(x)=\{\lambda\in \CC \ \vert \ x-\lambda
{\ is\ not\ invertible\ in\ } \tilde{A}\}$$ The {\it spectral radius} of
$X$ is the number: $$\nu(x)=sup\{ \vert \lambda \vert \ ; \
\lambda\in Sp(x)\}\ .$$
We have: $$\begin{array}{ccc} Sp(x)\subset \RR \mbox{ when } x \mbox{
    is self-adjoint } (x^*=x),\\ Sp(x)\subset \RR_+ \mbox{ when } x \mbox{
    is positive } (x=y^*y \mbox{ with } y\in A),\\ Sp(x)\subset U(1) \mbox{ when } x \mbox{
    is unitary } (x^*x=xx^*=1)\ .\end{array}$$
When $x$ is {\it normal}: $x^*x=xx^*$, these conditions on the
spectrum are equivalences.

\medskip \noindent When $x$ is normal, $\nu (x)=\|x\|$. From these,
one infers that for any
polynomial $P\in \CC[x]$, $\|P(x)\|=sup\{P(t)\ \vert \ t\in
Sp(x)\}$ (using that $Sp(P(x))=P(Sp(x))$).  We can then define $f(x)\in
A$ for every continuous function $f:Sp(x)\rightarrow \CC$. Precisely,
according to Weierstrass' theorem, there is a
sequence $(P_n)$ of polynomials which converges
uniformly to $f$ on $Sp(x)$. We simply define $f(x)=\lim P_n(x)$.

\subsection{The reduced and maximal $C^*$-algebras of a groupoid}
We  restrict our study to the case of
Hausdorff locally compact groupoids.  For the non Hausdorff case,
which is also important and not exceptional,  in particular when
dealing with foliations,  we refer the reader to \cite{Co0,Co3,KhSk} .

\smallskip \noindent From now on, $G\rightrightarrows
G^{(0)}$ is a locally compact Hausdorff groupoid equipped with a fixed Haar
system $\nu=\{\nu_x\}_{x\in G^{(0)}}$. We let $C_c(G)$ be the space of
complex valued functions with compact support on $G$. It is provided
with a structure of involutive algebra as follows. If $f$ and $g$
belong to $C_c(G)$ we define 

\noindent the {\it involution} by $$\mbox{for } \gamma\in G,\
f^*(\gamma)=\overline{f(\gamma^{-1})},$$
the {\it convolution product} by $$\mbox{for } \gamma\in G,\
f*g(\gamma)=\int_{\eta \in G_x} f(\gamma \eta^{-1})g(\eta)
d\nu_x(\eta),$$ 
\noindent where $x=s(\gamma)$. The $1$-norm on $C_c(G)$ is defined by
$$\|f\|_1=\underset{x\in G^{(0)}}{\sup}\max\left(\int_{G_x}\vert f(\gamma)
\vert d\nu_x(\gamma),\int_{G_x}\vert f(\gamma^{-1})
\vert d\nu_x(\gamma)\right).$$
The {\it groupoid full} $C^*$-algebra $C^*(G,\nu)$ is defined to be
the enveloping $C^*$-algebra of the Banach $*$-algebra
$\overline{C_c(G)}^{\|\cdot \|_1}$ obtained by
completion of $C_c(G)$ with respect to the norm $\| \cdot \|_1$.

\bigskip \noindent Given $x$ in $G^{(0)}$, $f$ in $ C_c(G)$, $\xi$ in
$L^2(G_x,\nu_x)$ and $\gamma$ in $G_x$, we set
$$\pi_x(f)(\xi)(\gamma)=\int_{\eta\in G_x}f(\gamma\eta^{-1}) \xi(\eta)
d\nu_x(\eta) \ .$$
One can show that $\pi_x$ defines a $*$-representation of $C_c(G)$ on
the Hilbert space $L^2(G_x,\nu_x)$. Moreover for any
$f\in C_c(G)$,  the inequality
$\|\pi_x(f)\|\leq \|f\|_1$ holds. The {\it reduced norm} on $C_c(G)$ is
$$\|f\|_r=\underset{x\in G^{(0)}}{sup} \{ \| \pi_x(f) \| \}$$ 
which defines a $C^*$-norm. The {\it reduced } $C^*$-algebra
$C_r(G,\nu)$ is 
defined to be the $C^*$-algebra obtained by completion of $A$ with
respect to $\| \cdot \|_r$. 

\smallskip \noindent When $G$ is smooth, the reduced and maximal
$C^*$-algebras of the 
groupoid $G$ do not depend up to isomorphism on the choice of
the Haar system $\nu$. In the general case they do not depend on $\nu$
up to Morita equivalence \cite{Re}. When there is no ambiguity on the
Haar system, we write $C^*(G)$ and $C^*_r(G)$ for the maximal
and reduced $C^*$-algebras.

\smallskip \noindent The identity map on $C_c(G)$ induces a
  surjective homomorphism from $C^*(G)$ to $C^*_r(G)$. Thus $C^*_r(G)$ is a quotient of
  $C^*(G)$. 

\noindent For a quite large class of groupoids, {\it amenable }
groupoids \cite{ARe}, the reduced and maximal $C^*$-algebras are
equal. This will be the case for all the groupoids we will meet in the
last part of this course devoted to index theory.

\medskip \noindent {\bf Examples 1.} When $X\rightrightarrows X$ is a
locally compact space, $C^*(X)=C^*_r(X)=C_0(X)$.\\
{\bf 2.} When $G\rightrightarrows e$ is a group and $\nu $ a Haar
measure on $G$, we recover the usual notion of reduced and maximal
$C^*$-algebras of a group.\\
{\bf 3.} Let $M$ be a smooth manifold and $TM \rightrightarrows M$ the
tangent bundle. Let us equip the vector bundle $TM$ with a euclidean
structure. The Fourier transform: 
 $$f\in C_{c}(TM), \ (x,w)\in T^*M, \quad  \hat{f}(x,w)=\frac{1}{(2\pi)^{n/2}}\int_{X\in
  T_xM}e^{-iw(X)}f(X)dX $$
gives rise to an isomorphism between $C^*(TM)=C^*_r(TM)$ and
      $C_0(T^{*}M)$. Here, $n$ denotes the dimension of $M$ and
      $T^{*}M$ the cotangent bundle of $M$. \\
{\bf 4.} Let $X$ be a locally compact space, with $\mu$ a measure on $X$
and consider the pair groupoid $X\times X \rightrightarrows X$. If
$f,\ g$ belongs to $C_c(X\times X)$, the convolution product is given
by: $$f*g(x,y)=\int_{z\in X}f(x,z)g(z,y)d\mu(z)$$ and a
representation of $C_c(X\times X)$ by
$$\pi:  C_c(X\times X) \rightarrow 
  \cL(L^2(X,\mu)) ; \  \pi(f)(\xi)(x)=\int_{z\in X} f(x,z)\xi(z) d\mu(z)$$
when $f\in C_c(X\times X), \xi \in L^2(X,\mu) \mbox{ and }
  x\in X  .$\\
It turns out that $C^*(X\times X)=C^*_r(X\times
X) \simeq \cK(L^2(X,\mu))$.\\
{\bf 5.} Let $M$ be a compact smooth manifold and
$\cG^t_M\rightrightarrows M\times [0,1]$ its tangent groupoid. In this
situation $C^*(\cG^t_M)=C^*_r(\cG^t_M)$ is a continuous field
$(A_t)_{t\in [0,1]}$ of $C^*$-algebras (\cite{Di}) with $A_0\simeq C_0( T^*M)$
 a commutative $C^*$-algebra and for any $t\in ]0,1]$,
  $A_t\simeq \cK(L^2(M))$  \cite{Co0}.

\medskip \noindent In the sequel we will need the two following
properties of $C^*$-algebras of groupoids.

\smallskip \noindent{\bf Properties 1.} Let $G_1$ and $G_2$ be two
locally compact
groupoids equipped with Haar systems and suppose for instance that
$G_1$ is amenable. Then 
according to \cite{ARe}, $C^*(G_1)=C^*_r(G_1)$ is {\it nuclear} - which implies
that for any $C^*$-algebra $B$ there is only one tensor product $C^*$-algebra
$C^*(G_1)\otimes B$. The groupoid $G_1\times G_2$ is locally compact
 and 
$$C^*(G_1\times G_2)\simeq C^*(G_1)\otimes C^*(G_2) \
\mbox{ and } C^*_r(G_1\times G_2)\simeq C^*(G_1)\otimes C^*_r(G_2)\ 
.$$  
{\bf 2.} Let $G\rightrightarrows G^{(0)}$ be a locally compact
groupoid with a Haar system $\nu$.\\ An open subset $U\subset G^{(0)}$
is {\it saturated} if $U$ is a union of orbits of $G$, in other words
if $U=s(r^{-1}(U))=r(s^{-1}(U))$. The set $F=G^{(0)}\setminus U$ is
then a
closed saturated subset of $G^{(0)}$. The Haar system $\nu$ can be 
restricted to the restrictions $G\vert_U:=G_U^U$ and $G\vert_F:=G_F^F$
and we have the following exact sequence of $C^*$-algebras
\cite{HS,Ra}:
$$0\rightarrow C^*(G\vert_U) \overset{i}{\rightarrow} C^*(G)
\overset{r}{\rightarrow} C^*(G\vert_F) \rightarrow 0$$ 
where $i: C_c(G\vert_U) \rightarrow C_c(G)$ is the
extension of functions by $0$ while $r:C_c(G)\rightarrow
C_c(G\vert_F)$ is the restriction of functions.


\newpage\part*{KK-THEORY}

\medskip This part on $KK$-theory starts with a historical
introduction. In order to motivate our purpose we list most
of the properties of the $KK$-functor. Sections \ref{HilbMod} to
\ref{finKK} are devoted to a detailed description of the ingredients involved in
$KK$-theory.
\smallskip \noindent
As already pointed out in the introduction  we made an intensive use of the following
references \cite{ska,Hig1,ska-DEA,WO}. Moreover a significant part of this chapter has
been written by Jorge Plazas from the lectures held in Villa de Leyva
and we would like to thank him
for his great help.

\section{Introduction to KK-theory}\label{intro-KK}

\subsection{Historical comments}
\noindent The story begins with several studies of M. Atiyah \cite{At,ATout}.

\smallskip Firstly, recall that if $X$ is a compact space, the {\it $K$-theory} of $X$
is constructed in the following way: let $\cE v$ be the set of
isomorphism classes of continuous vector bundles over $X$. Thanks to the
direct sum of bundles, the set $\cE v$ is
naturally endowed with a structure of abelian semi-group. One can then
symetrize $\cE v$ in order to get a group, this gives the $K$-theory
group of $X$:
 $$
  K^0(X)=\{[E]-[F] \ ; \ [E], [F] \in \cE v \}. 
 $$
For example the $K$-theory of a point is $\ZZ$ since a vector bundle
on a point is just a vector space and vector spaces are classified, up
to isomorphism, by their dimension.

\smallskip A first step towards $KK$-theory is the discover, made by
M. Atiyah \cite{At} and independently K. J\"anich \cite{Janich}, that
$K$-theory of a compact space $X$ can be described with Fredholm
operators.  

\smallskip \noindent When $\cH$ is an infinite dimensional separable
Hilbert space, the set $\cF(\cH)$ of {\it Fredholm operators} on $\cH$ is the open
subset of $\cL(\cH)$ made of
bounded operators $T$ on $\cH$ such that the dimension of the kernel and
cokernel of $T$ are finite. The set $\cF(\cH)$ is stable under
composition. We set
 $$
  [X,\cF(\cH)]=\{ \mbox{homotopy classes of continuous
  maps: } X\rightarrow \cF(\cH)\}. 
 $$ 
The set $[X,\cF(\cH)]$ is naturally endowed with a structure of semi-group. 
M. Atiyah and K. J\"anich, showed that $[X,\cF(\cH)]$ is
actually (a group) isomorphic to $K^0(X)$ \cite{At}. The idea of the proof is the
following. If $f:X\rightarrow \cF(\cH)$ is a continuous map, one can
choose a subspace $V$ of $\cH$ of finite codimension such that:
\begin{equation}
  \label{Atiyah-Janich-key-step}
 \forall x\in X,  V\cap \ker f_{x}=\{0\} \hbox{ and } \bigcup_{x\in
   X} \cH/f_{x}(V) \hbox{ is a vector bundle}. 
\end{equation}
Denoting by $\cH/f(V)$ the vector bundle arising in
(\ref{Atiyah-Janich-key-step}) and by $\cH/V$ the product bundle
$X\times \cH/V$, the Atiyah-Janich isomorphism is then given by:
\begin{equation}
  \label{Atiyah-Janich-iso}
  \begin{array}{ccc} \lbrack X,\cF(\cH)\rbrack & \rightarrow & K^0(X) \\ \lbrack f\rbrack & \mapsto &
 [\cH/V]-[\cH/f(V)].  \end{array}
\end{equation}
Note that choosing $V$ amounts to modify $f$ inside its homotopy
class into $\widetilde{f}$ (defined to be equal to $f$ on $V$ and
to $0$ on a supplement of $V$) such that:
\begin{equation}
  \label{proto-Generalized-Fred}
  \mbox{Ker}\tilde f:=\cup_{x\in X} \mbox{Ker}(\tilde f_x)  \hbox{
    and } \mbox{CoKer}\tilde f:=\cup_{x\in X} \cH/\tilde f_x(\cH) 
\end{equation}
are vector bundles over $X$. These constructions contain relevant
information for the sequel: the map $f$ arises as a {\sl generalized
  Fredholm} operator on the {\sl Hilbert $C(X)$-module} $C(X,\cH)$.

\medskip Later, M. Atiyah tried to describe the dual functor
$K_{0}(X)$, the $K$-homology of $X$, with the help of Fredholm operators. This
gave rise to $\rm{Ell}(X)$ whose cycles are triples $(H,\pi,F)$ where:
\begin{itemize}
\item[-] $H=H_0\oplus H_1$ is a $\ZZ_2$ graded Hilbert space,
\item[-] $\pi:C(X)\rightarrow \cL(H)$ is a representation by operators of
  degree $0$ (this means that $\pi(f)=\begin{pmatrix} \pi_0(f) & 0 \\ 0 &
    \pi_1(f) \end{pmatrix}$),
\item[-] $F$ belongs to $\cL(H)$, is of degree $1$ (thus it is of the
  form 
  $F=\begin{pmatrix}   0 & G \\ T & 0 \end{pmatrix})$ and satisfies
  $$F^2-1\in \cK(H) \mbox{ and } [\pi,F]\in \cK(H)\ .$$
In particular $G$ is an inverse of $T$ modulo compact operators.
\end{itemize}
Elliptic operators on closed manifolds produce natural examples of
such cycles.  Moreover, there exists a natural pairing between
$\rm{Ell}(X)$ and $K^0(X)$,  justifying the choice of $\rm{Ell}(X)$
as a candidate for the cycles of the $K$-homology of $X$:
\begin{equation}
  \label{pairing-Ell}
  \begin{array}{ccc} K^0(X)\times \rm{Ell}(X) & \rightarrow & \ZZ \\
  ([E],(H,\pi,F)) & \mapsto & \mbox{Index}(F_E) \end{array}
\end{equation}
where $\mbox{Index}(F_E)=\mbox{dim} (\mbox{Ker}(F_E))-\mbox{dim}(\mbox{CoKer}(F_E))$ is
the {\it index} of a Fredholm operator associated to a vector bundle
$E$ on $X$ and a cycle $(H,\pi,F)$ as follows.\\ Let $E'$ be a vector
bundle on $X$ such that $E\oplus E'\simeq \CC^N \times
X$ and let $e$ be the projection of $\CC^N\times X$ onto $E$. We can
identify $C(X,\CC^N)\underset{C(X)}{\otimes} H$ with $H^N$. Let
$\tilde{e}$ be the image of $e\otimes 1$ under this
identification. We define $F_E:=\tilde e F^N \vert_{\tilde e(H^N)}$ where
$F^N$ is the diagonal operator with $F$ in each diagonal entry. The
operator $F_E$ is the desired Fredholm operator on $\tilde e (H^N)$.

\medskip Now, we should recall that to any $C^*$-algebra $A$
(actually,  to any ring) is associated a
group $K_0(A)$.  When $A$ is unital, it can be defined as follows:

\begin{center} $K_0(A)=\{ [\cE]-[\cF] \ ; \ [\cE], [\cF] \mbox{ are isomorphism
classes of} \linebreak[4] \mbox{finitely generated projective $A$-modules}\}\ .$
\end{center} 

\smallskip \noindent Recall that a $A$-module
$\cE$ is finitely generated  and projective if there exists another $A$-module $\cG$
such that $\cE\oplus \cG \simeq A^N$ for some integer $N$.

\smallskip \noindent The Swan-Serre theorem asserts that for any
compact space $X$, the category of (complex) vector bundles over $X$ is
equivalent to the category of finitely generated projective modules
over $C(X)$,  in particular: $K^0(X)\simeq K_0(C(X))$. This fact and
the ($C^{*}$-)algebraic flavor of the constructions above leads to
the natural attempt to generalize them for noncommutative $C^{*}$-algebras.

\medskip \noindent During the 79 $\sim$ 80's G. Kasparov  
defined  with great success for any pair of $C^*$-algebras a
bivariant theory, the $KK$-theory. This theory generalizes both
$K$-theory and $K$-homology and carries a
product generalizing the pairing (\ref{pairing-Ell}). Moreover,  in many cases $KK(A,B)$ contains all the
morphisms from $K_0(A)$ to $K_0(B)$. To understand
this bifunctor,  we will study the notions of Hilbert modules, of
adjointable operators acting on them and  of generalized Fredholm
operators which generalize to arbitrary $C^{*}$-algebras the notions
encountered above for $C(X)$. Before going to this 
functional analytic part, we end this introduction by listing most of the properties of the bi-functor $KK$.

\subsection{Abstract properties of $KK(A,B)$}
Let $A$ an $B$ be two $C^*$-algebras. In order to simplify our
presentation, we assume that $A$ and $B$ are separable. Here is the
list of the most important properties of the $KK$ functor.

\smallskip \noindent {\bf $\bullet$ } $KK(A,B)$ is an abelian group.

\smallskip \noindent {\bf $\bullet$ Functorial properties} The functor $KK$ is covariant in $B$
and contravariant in $A$: if $f:B\rightarrow C$ and
$g:A\rightarrow D$ are two  homomorphisms of
$C^*$-algebras, there exist two homomorphisms of groups:
$$f_*:KK(A,B)\rightarrow KK(A,C) \mbox{ and } g^*:KK(D,B)\rightarrow
KK(A,B)\ .$$
In particular $id_*=id$ and $id^*=id$.

\smallskip \noindent {\bf $\bullet$} Each *-morphism $f:A\rightarrow B$ defines an element, denoted
by $[f]$, in $KK(A,B)$. We set $1_A:=[id_A]\in KK(A,A)$.

\smallskip \noindent {\bf $\bullet$ Homotopy invariance} $KK(A,B)$ is homotopy invariant.\\
Recall that the $C^*$-algebras $A$ and $B$ are {\it homotopic}, if there
exist two *-morphisms $f:A\rightarrow B$ and $g:B\rightarrow A$ such
that $f\circ g$ is homotopic to $id_B$ and $g\circ f$ is homotopic to
$id_A$.\\
Two homomorphisms $F,G:A\rightarrow B$ are homotopic when there exists a
$*$-morphism $H:A\rightarrow C([0,1],B)$ such that $H(a)(0)=F(a)$ and
$H(a)(1)=G(a)$ for any $a\in A$.

\smallskip \noindent {\bf $\bullet$ Stability} If $\cK$ is the algebra of compact operators on a
Hilbert space,  there are isomorphisms:
$$KK(A,B\otimes \cK)\simeq
KK(A\otimes \cK,B) \simeq KK(A,B)\ .$$

\smallskip \noindent More generally, the bifunctor $KK$ is invariant under {\it Morita
  equivalence}.

\smallskip \noindent {\bf $\bullet$ Suspension} If $E$ is a $C^*$-algebra there exists an
homomorphism $$\tau_E:KK(A,B)\rightarrow KK(A\otimes E,B\otimes E)$$
which satisfies $\tau_E \circ \tau_D =\tau_{E\otimes D}$ for any
$C^*$-algebra $D$.

\smallskip \noindent {\bf $\bullet$ Kasparov product} There is a well defined
bilinear coupling: $$\begin{array}{ccc} KK(A,D)\times KK(D,B) &
  \rightarrow & KK(A,B) \\ (x,y) & \mapsto & x\otimes y \end{array}$$
called the {\it Kasparov product}. It is associative, covariant in $B$ and
contravariant in $A$: if $f:C\rightarrow A$ and
$g:B\rightarrow E$ are two  homomorphisms of
$C^*$-algebras then $$f^*(x\otimes y)=f^*(x)\otimes y \mbox{ and }
g_*(x\otimes y)=x\otimes g_*(y).$$
If $g:D\rightarrow C$ is another *-morphism, $x\in KK(A,D)$
and $z\in KK(C,B)$ then $$h_*(x)\otimes z=x\otimes h^*(z)\ .$$
Moreover, the following equalities hold: $$f^*(x)=[f]\otimes x \ , \
g_*(z)=z\otimes [g] \mbox{ and }[f\circ h]=[h]\otimes [f]\ .$$ In
particular $$x\otimes 1_D=1_A\otimes x = x\ .$$
The Kasparov product behaves well with respect to suspensions. If $E$ is a
$C^*$-algebra: $$\tau_E(x\otimes y)=\tau_E(x)\otimes \tau_E(y)\ .$$
This enables to extend the Kasparov product: $$\begin{array}{cccc}
  \underset{D}{\otimes}: & KK(A,B\otimes D)\times KK(D\otimes C,E) &
  \rightarrow & KK(A\otimes C,B\otimes E) \\ & (x,y) & \mapsto &
  x\underset{D}{\otimes}y:=\tau_C(x)\otimes \tau_B(y) \end{array}$$

\smallskip \noindent {\bf $\bullet$} The Kasparov product $\underset{\CC}{\otimes}$ is commutative.

\smallskip \noindent {\bf $\bullet$ Higher groups} For any $n\in \NN$, let
$$KK_n(A,B):=KK(A,C_0(\RR^n)\otimes B)\ .$$ 
An alternative definition, leading to isomorphic groups,  is 
$$KK_n(A,B):=KK(A,C_{n}\otimes B), $$
where  $C_{n}$ is the Clifford algebra of $\CC^{n}$. This will be explained later.
The functor $KK$ satisfies {\bf Bott periodicity}: there is an isomorphism $$KK_2(A,B)\simeq
KK(A,B)\ .$$

\smallskip \noindent {\bf $\bullet$ Exact sequences} Consider the following exact
sequence of $C^*$-algebras:$$0\rightarrow J \overset{i}{\rightarrow}
A \overset{p}{\rightarrow} Q \rightarrow 0$$
and let $B$ be another $C^*$-algebra. Under a few more assumptions (for
example all the $C^*$-algebras are nuclear or $K$-nuclear,  or the above
exact sequence admits a completely positive norm decreasing
cross section \cite{ska-88}) we have the following two periodic exact
sequences 
$$
 \begin{CD} KK(B,J) @>i_*>> KK(B,A) @>p_*>>KK(B,Q) \\
  @A{\delta}AA     @.    @VV{\delta}V \\
  KK_1(B,Q) @<<p_*< KK_1(B,A)
  @<<i_*< KK_1(B,J) 
 \end{CD}
$$
$$
 \begin{CD} KK(Q,B) @>p^*>> KK(A,B) @>i^*>> KK(J,B) \\
  @A{\delta}AA   @.   @VV{\delta}V \\ KK_1(J,B) @<<i^*< KK_1(A,B)
  @<<p^*< KK_1(Q,B) 
 \end{CD}
$$
where the connecting homomorphisms $\delta$ are given by Kasparov
products.  

\medskip \noindent{\bf $\bullet$ Final remark} Let us go back to the end of the
introduction in order to make it more 
precise. \\
The usual $K$-theory groups appears as special cases of $KK$-groups:
$$KK(\CC,B)\simeq K_0(B), $$ 
while the {\it $K$-homology} of a
$C^*$-algebra $A$ is defined by $$K^0(A)=KK(A,\CC)\ .$$ 
Any $x\in KK(A,B)$ induces a homomorphism of groups:
$$\begin{array}{ccc}KK(\CC,A)\simeq K_0(A) & \rightarrow &
  K_0(B)\simeq KK(\CC,B) \\ \alpha & \mapsto & \alpha \otimes x
\end{array}$$
In most situations, the induced homomorphism $KK(A,B)\rightarrow
Mor(K_0(A),K_0(B))$ is surjective. Thus one can think of $KK$-elements
as homomorphisms between $K$-groups.

\medskip \noindent When $X$ is a compact space,
one has $K^0(X)\simeq K_0(C(X))\simeq KK(\CC,C(X))$ and as we will
see shortly, $K^{0}(C(X))=KK(C(X), \CC)$ is a quotient of the set $Ell(X)$ introduced by
M. Atiyah. Moreover the pairing $K^0(X)\times Ell(X) \rightarrow \ZZ$
coincides with the Kasparov product
$KK(\CC,C(X))\times KK(C(X),\CC) \rightarrow KK(\CC,\CC)\simeq\ZZ$.


\section{Hilbert modules} \label{HilbMod}

We review the main properties of Hilbert modules over
$C^*$-algebras,  necessary for a correct understanding of bivariant
$K$-theory. We closely follow the presentation given by G. Skandalis
\cite{ska-DEA}. Most proofs given below are taken from his lectures on
the subject. We are indebted to him for allowing us to use his
lectures notes.  Some of the material given below can also be found in \cite{WO}, 
where the reader will find a guide to the literature and a more detailed
presentation. 

\subsection{Basic definitions and examples}
Let $A$ be a $C^*$-algebra and $E$ be a $A$-right module.

\smallskip \noindent A sesquilinear form $(\cdot,\cdot):E\times E\to
A$ is {\it positive} if
for all $x\in E$, $(x,x)\in A_+$. Here $A_+$ denotes the set of
positive elements in $A$. It is {\it positive definite} if moreover
$(x,x)=0$ if and only if $x=0$.

\smallskip \noindent Let $(\cdot,\cdot): E\times E\to A$ be a positive sesquilinear form and set 
 $Q(x)=(x,x)$. By the polarization identity: 
 $$
 \forall x,y\in E,\quad   (x,y) = \frac{1}{4}\left( Q(x+y)-iQ(x+iy)-Q(x-y)+iQ(x-iy)\right)
 $$
we get: 
 $$ \forall x,y\in E, \quad (x,y)=(y,x)^* $$

\begin{definition}
   A {\it pre-Hilbert $A$-module} is a right $A$-module $ E$ with a
   positive definite sesquilinear map $(\cdot,\cdot): E\times E\to A$
  such that $y\mapsto (x,y)$ is $A$-linear. 
\end{definition}

\begin{proposition} Let $( E,(\cdot,\cdot))$ be a pre-Hilbert $A$-module. The
  following:
  \begin{equation}
    \label{eq:norm-hilbert}
    \forall x\in E,\quad  \| x\| = \sqrt{\|(x,x)\|}
  \end{equation}
defines a norm on $E$. 
\end{proposition}
The only non trivial fact is the triangle inequality, which results
from: 
\begin{lemma}(Cauchy-Schwarz inequality) 
 $$\forall x,y\in E,\quad (x,y)^*(x,y) \le \|x\|^2(y,y) $$
In particular: $\|(x,y)\|\le \|x\|\|y\|$. 
\end{lemma}
\noindent Set $a=(x,y)$. We have for all $t\in \RR$: $ (xa+ty,xa+ty) \ge 0$,
thus: 
\begin{equation}
  \label{eq:CSproof1}
2ta^*a \le a^*(x,x)a + t^2(y,y)  
\end{equation}
Since $(x,x)\ge 0$, we have: $a^*(x,x)a\le \|x\|^2a^*a$ (it uses the
equivalence: $z^*z\le w^*w$ if and only if $\|zx\|\le \|wx\|$ for all $x\in A$)
and choosing $t=\|x\|^2$ in (\ref{eq:CSproof1}) gives the result. 

\begin{definition}
  A {\it Hilbert $A$-module} is a pre-Hilbert $A$-module which is complete
  for the norm defined in (\ref{eq:norm-hilbert}).\\
 A {\it Hilbert $A$-submodule} of a Hilbert $A$-module is a closed
 $A$-submodule provided with the restriction of the $A$-valued scalar product. 
\end{definition}
\noindent When there is no ambiguity about the base $C^*$-algebra $A$, we simply
say  pre-Hilbert module and Hilbert module. 

\smallskip \noindent Let $(E,(\cdot,\cdot))$ be a pre-Hilbert $A$-module. From the continuity
of the sesquilinear form $(\cdot,\cdot): E\times E\to A$ and of the
right multiplication
$E\to E, x\mapsto xa$ for any
$a\in A$, we infer that the completion of $E$ for the norm
(\ref{eq:norm-hilbert}) is a Hilbert $A$-module. 
\begin{remark}
  In the definition of a pre-Hilbert $A$-module, one could remove the
  hypothesis $(\cdot,\cdot)$ is {\sl definite}. In that case,
  (\ref{eq:norm-hilbert}) defines a semi-norm and one checks that the
  Hausdorff completion of a pre-Hilbert $A$-module, in this extended
  sense, is a Hilbert $A$-module. 
\end{remark}

\noindent We continue this paragraph with classical examples.  

\noindent {\bf 1.} The algebra  $A$ is a Hilbert $A$-module with its obvious right $A$-module
structure and: 
 $$(a,b):= a^*b \ .$$ 

\noindent {\bf 2.} For any positive integer $n$, $A^n$ is a Hilbert $A$-module with its obvious right $A$-module
structure and: 
 $$ ((a_i),(b_i)):= \sum_{i=1}^na_i^*b_i \ .$$
Observe that $\sum_{i=1}^na_i^*a_i $ is a sum of positive elements in
$A$, which implies that 
 $$ 
   \|(a_i)\|=\sqrt{ \| \sum_{i=1}^na_i^*a_i \|}\ge \|a_k\| 
 $$
for all $k$. It follows that if $(a^m_1,\ldots,a^m_n)_{m}$ is a Cauchy sequence in $A^n$, the
sequences $(a^m_k)_m$ are Cauchy  in $A$, thus convergent and we
conclude that $A^n$ is complete.

\noindent {\bf 3.} Example {\bf 2.} can be extended to the direct sum of $n$ Hilbert
$A$-modules $E_1,\ldots,E_n$ with the Hilbertian product:
 $$ ((x_i),(y_i)):= \sum_{i=1}^n(x_i,y_i)_{{E_{i}}} $$

\noindent {\bf 4.} If $F$ is a closed $A$-submodule of a Hilbert $A$-module $E$ then
$F$ is a Hilbert $A$-module. 
 For instance, a closed right ideal in $A$ is a Hilbert $A$-module. 

\noindent {\bf 5.} {\bf The standard Hilbert $A$-module} is defined by
 \begin{equation}
   \label{eq:standard-A-module}
   \cH_A = \{ x=(x_k)_{k\in\NN}\in A^{\NN} \ | \ \sum_{k\in\NN} x_k^*x_k \hbox{
     converges } \}. 
 \end{equation}
The right $A$-module structure is given by $(x_k) a = (x_ka)$ and the
Hilbertian $A$-valued product is: 
\begin{equation}
  \label{eq:standard-module-product}
  \left( (x_k) , (y_k)\right) = \sum_{k=0}^{+\infty} x_k^*y_k 
\end{equation}
This sum converges for elements of $\cH_A$, indeed for all $q>p\in\NN$:
\begin{eqnarray*}
  \|\sum_{k=p}^{q} x_k^*y_k\| & = & \|\left( (x_k)_p^q ,
    (y_k)_p^q\right)_{A^{q-p}}\| \\
    & \le &  \| (x_k)_p^q\|_{A^{q-p}}\| (y_k)_p^q\|_{A^{q-p}} \quad
    (\hbox{Cauchy Schwarz inequality in} A^{q-p})\\
     & = & \sqrt{\|\sum_{k=p}^{q} x_k^*x_k\| }\sqrt{\|\sum_{k=p}^{q} y_k^*y_k\| }
\end{eqnarray*}
This implies that $ \sum_{k\ge 0} x_k^*y_k $ satisfies the Cauchy
criterion, and therefore converges, so that (\ref{eq:standard-module-product})
makes sense. Since  for all $(x_k)$, $(y_k)$ in $\cH_A$:
 $$ 
   \sum_{k\ge 0} (x_k+y_k)^*(x_k+y_k)=\sum_{k\ge 0}
   x_k^*x_k+\sum_{k\ge 0} y_k^*x_k +\sum_{k\ge 0} x_k^*y_k+\sum_{k\ge 0} y_k^*y_k
 $$
is the sum of four convergent series, we find that $(x_k)+(y_k)=(x_k+y_k)$
is in $\cH_A$. We also have, as before, that for all $a\in A$ and 
$(x_k)\in\cH_A$:
 $$
    \|\sum_{k=0}^{+\infty} (x_ka)^*(x_ka)\|\le \|a\|^2
    \|\sum_{k=0}^{+\infty} x_k^*x_k\|
 $$
Hence, $\cH_A$ is a pre-Hilbert $A$-module, and we need to check that
it is complete. Let $(u_n)_n=((u_i^n))_n$ be a Cauchy sequence in $\cH_A$. We get,
as in Example {\bf 2.}, that for all $i\in\NN$, the sequence $(u_i^n)_n$ is Cauchy in $A$,
thus converges to an element denoted $v_i$. Let us check that $(v_i)$
belongs to $\cH_A$.\\ Let $\varepsilon>0$. Choose $n_0$ such that
$$\forall p>q\ge n_0 , \ \| u_q-u_p\|_{\cH_A}\le \varepsilon/2 \ .$$ Choose
$i_0$ such that $$\forall k>j\ge i_0,\ 
\|\sum_{i=j}^{k}u^{n_0}_i{}^*u^{n_0}_i\|^{1/2}\le
\varepsilon/2 \ .$$ Then thanks to the triangle inequality in $A^{k-j}$ we
get for all $p,q\ge n_0$ and $j,k\ge i_0$:
\begin{eqnarray*}
  \|\sum_{i=j}^{k}u^{p}_i{}^*u^{p}_i\|^{1/2} & \le & 
       \| \sum_{i=j}^{k}(u^{p}_i-u^{n_0}_i)^*(u^{p}_i-u^{n_0}_i)\|^{1/2}
       +\|\sum_{i=j}^{k}u^{n_0}_i{}^*u^{n_0}_i\|^{1/2} \le  \varepsilon 
\end{eqnarray*}
Taking the limit $p\to +\infty$, we get: 
$\|\sum_{i=j}^{k}v_i^*v_i\|^{1/2} \le \varepsilon$ for all 
$j,k\ge i_0$ which implies that $(v_i)\in\cH_A$. It remains to check
that $(u_n)_n$ converges to $v=(v_i)$ in $\cH_A$. With the notations
above:
 $$
  \forall p,q\ge n_0,\ \forall I\in\NN,\quad 
  \| \sum_{i=0}^{I}(u^{p}_i-u^{q}_i)^*(u^{p}_i-u^{q}_i)\|^{1/2}\le
  \varepsilon,
 $$
taking the limit $p\to +\infty$:
 $$ 
   \forall q\ge n_0,\ \forall I\in\NN,\quad  \| \sum_{i=0}^{I}(v_i-u^{q}_i)^*(v_i-u^{q}_i)\|^{1/2}\le
  \varepsilon,
 $$
taking the limit $I\to +\infty$:
 $$ 
  \forall q\ge n_0,\quad  \| v-u_q \|\le   \varepsilon,
 $$
which ends the proof. \hfill{\small $\Box$ }

\medskip \noindent The standard Hilbert module $\cH_A$ is maybe the most important
Hilbert module. Indeed, Kasparov proved: 
\begin{theorem}
  Let $E$ be a countably generated Hilbert $A$-module. Then $\cH_{A}$ and
  $E\oplus \cH_A$ are isomorphic.
\end{theorem}
\noindent The proof can be found in \cite{WO}. 
This means that there exists a $A$-linear unitary map 
$U: E\oplus \cH_A\to \cH_{A}$. The notion of unitary uses the notion of
adjoint, which will be explained later. 

\begin{remark}
 {\bf 1.}  The algebraic sum $\underset{\NN}{\oplus} A$ is dense in
  $\cH_A$. \\
 {\bf 2.} We can replace in $\cH_A$ the summand $A$ by any sequence of Hilbert
  $A$-modules $(E_i)_{i\in\NN}$ and the Hilbertian $A$-valued product
  by:
  $$ 
      \left( (x_k) , (y_k)\right) = \sum_{k=0}^{+\infty}
      (x_k,y_k)_{E_k}
  $$
 If $E_i=E$ for all $i\in\NN$, the resulting Hilbert $A$-module is
 denoted by
 $l^2(\NN,E)$. \\
 {\bf 3.} We can generalize the construction to any family $(E_i)_{i\in I}$
 using summable families instead of convergent series. 
\end{remark}

\noindent We end this paragraph with two concrete examples. 

\smallskip \noindent {\bf a.} Let $X$ be a locally compact space and $E$ an hermitian
vector bundle. The space $C_0(X,E)$ of continuous sections of $E$
vanishing at infinity is a Hilbert $C_0(X)$-module with the module
structure given by: $$\xi.a(x)=\xi(x)a(x),\  \xi\in C_0(X,E),\ 
a\in C_0(X) $$
and the $C_{0}(X)$-valued product given by: 
 $$(\xi,\eta)(x)=(\xi(x),\eta(x))_{E_x}\ $$ 

\noindent {\bf b.} Let $G$ be a locally compact groupoid with a Haar system
$\lambda$ and $E$ a hermitian vector bundle over $G^{(0)}$. Then 
\begin{equation}
  \label{eq:hilbert-GE}
   f,g\in C_c(G,r^*E),\quad
   (f,g)(\gamma)=\int_{G_{s(\gamma)}}(f(\eta \gamma^{-1}),g(\eta))_{E_{r(\eta)}}d\lambda^{s(\gamma)}(\eta) 
\end{equation}
gives a positive definite sesquilinear $C_c(G)$-valued form which has
the correct behavior with respect to the right action of $C_c(G)$ on 
$C_c(G,r^*E)$.  This leads to two norms
$\|f\|=\|(f,f)\|_{C^*(G)}^{1/2}$ and $\|f\|_r=\|(f,f)\|_{C^*_r(G)}^{1/2}$
and two completions of $C_c(G,r^*E)$, denoted $C^*(G,r^*E)$ and 
$C^*_r(G,r^*E)$ which are Hilbert modules, respectively over $C^*(G)$
and $C^*_r(G)$. 

\subsection{Homomorphisms of Hilbert $A$-modules}

Let $E,F$ be Hilbert $A$-modules. We will need the orthogonality in
Hilbert modules: 
\begin{lemma}
  Let $S$ be a subset of $E$. The orthogonal of $S$:
 $$ S^\perp = \{ x\in E \ | \ \forall y\in S, \ (y,x)=0 \} $$
  is a Hilbert $A$-submodule of $E$.
\end{lemma}

\subsubsection{Adjoints}

Let $T: E \to F$ be a map. $T$ is {\it adjointable} if there exists a map 
$S: F\to E$ such that: 
\begin{equation}
  \label{eq:adjoint-condition}
  \forall (x,y)\in E\times F,\quad (Tx,y) = (x,Sy) 
\end{equation}
\begin{definition}
Adjointable maps are called {\it homomorphisms of Hilbert $A$-modules}.  
The set of adjointable maps from $E$ to $F$ is denoted by 
$\mor (E,F)$, and $\mor(E)=\mor(E,E)$. The space of linear continuous
maps from $E$ to $F$ is denoted by $\cL(E, F)$ and $\cL(E)=\cL(E, E)$. 
\end{definition}
\noindent The terminology will become clear after the next proposition.
\begin{proposition}
  Let $T\in\mor(E,F)$. 
  \begin{itemize}
  \item[(a)] The operator satisfying (\ref{eq:adjoint-condition}) is
    unique. It is denoted by $T^*$ and called the adjoint of $T$. One has
    $T^*\in\mor(F,E)$ and $(T^*)^*=T$.
  \item[(b)] $T$ is linear, $A$-linear and continuous. 
  \item[(c)] $\|T\|=\|T^*\|$, $\|T^*T\|=\|T\|^2$, $\mor(E,F)$ is a
    closed subspace of $\cL(E,F)$. In particular $\mor(E)$ is a
    $C^*$-algebra.
  \item[(d)] If $S\in\mor(E,F)$ and $T\in\mor(F,G)$ then
    $TS\in\mor(E,G)$ and $(TS)^*=S^*T^*$. 
  \end{itemize}
\end{proposition}
\begin{proof}
\noindent  (a)  Let $R,S$ be two maps satisfying
(\ref{eq:adjoint-condition}) for $T$. Then:
$$ \forall x\in E, y\in F,\quad (x, Ry-Sy)= 0 $$
and taking $x=Ry-Sy$ yields $Ry-Sy=0$. The remaining part of the
assertion is obvious. 

\noindent (b) $\forall x,y\in E, z\in F, \lambda\in \CC,$
$$ \quad 
 (T(x+\lambda y),z)=(x+\lambda y,T^*z)=(x,T^*z)+\overline{\lambda}
 (y,T^*z)= (Tx,z)(\lambda Ty,z) $$
thus $T(x+\lambda y)=Tx+\lambda Ty$ and $T$ is linear. Moreover:
$$\forall x\in E, y\in F, a\in A,\quad
(T(xa),y)=(xa,T^*y)=a^*(x,T^*y)=((Tx)a,y), $$
which gives the $A$-linearity. Consider the set 
 $$ S = \{ (-T^*y,y)\in E\times F \ | y\in F\} \ .$$
Then 
\begin{eqnarray*}
  (x_0,y_0)\in S^\perp &\Leftrightarrow & \forall y\in F,\
  (x_0,-T^*y)+(y_0,y) = 0 \\
  &\Leftrightarrow & \forall y\in F,\ (y_0-Tx_0,y)=0
\end{eqnarray*}
Thus $ G(T) = \{ (x,y)\in E\times F\ | \ y=Tx\}=S^\perp$ is closed and
the closed graph theorem implies that $T$ is continuous. 

\noindent (c) We have:
 $$ \| T\|^2 = \sup_{\|x\|\le 1} \| Tx\|^2 = \sup_{\|x\|\le
   1}(x,T^*Tx)\le \| T^*T\| \le  \| T^*\| \|T\| \ .
 $$
Thus $\| T\|\le \| T^*\|$ and switching $T$ and $T^*$ gives the
equality.\\ One has also proved: 
 $$  \| T\|^2 \le \| T^*T\| \le \| T^*\| \|T\| = \| T\|^2 $$
thus $\| T^*T\| = \| T\|^2 $ and the norm of $\mor(E)$ satisfies the
$C^*$-algebraic equation.

\smallskip \noindent  Let $(T_n)_n$ be a  sequence in $\mor(E,F)$, which converges to
 $T\in\cL(E,F)$. Since $\|T\|=\|T^*\|$
 and since $T\to T^*$ is (anti-)linear, the
 sequence $(T^*_n)_n$ is a Cauchy sequence, and therefore converges to an
 operator $S\in\cL(F,E)$. It then immediately follows that $S$ is the adjoint
 of $T$. This proves that $\mor(E,F)$ is closed, in particular
 $\mor(E)$ is a $C^*$-algebra. 

\noindent (d) Easy. 
\end{proof}
\begin{remark}
  There exist continuous linear and $A$-linear maps $T: E\to F$ which
  do not have an adjoint. For instance, take $A=C([0,1])$,
  $J=C_0(]0,1])$ and $T: J\hookrightarrow A$ the inclusion. Assuming that $T$ is
  adjointable, a one line computation proves that $T^*1 =
  1$. But $1$ does not belong to $J$. Thus $J\hookrightarrow A$ has no
  adjoint. \\
  One can also take $E=C([0,1])\oplus C_0(]0,1])$ and
  $T: E\to E,x+y\mapsto y+0$ to produce an example of $T\in \cL(E)$
  and $T\not\in \mor(E)$.  
\end{remark}
\noindent One can characterize self-adjoint and positive elements in the
$C^*$-algebra $\mor(E)$ as follows. 
\begin{proposition}\label{char-self-adjt-pos}
 Let $T\in\mor(E)$. \\
  (a)  $\ds T=T^* \Leftrightarrow \forall x\in E,\ (x,Tx)=(x,Tx)^*$\\
  (b) $\ds T\ge 0 \Leftrightarrow \forall x\in E,\ (x,Tx)\ge 0$
\end{proposition}
\begin{proof}
 (a) The implication $(\Rightarrow)$ is obvious. Conversely, set 
$Q_T(x)=(x,Tx)$. Using the polarization identity:
 $$ (x,Ty)= \frac{1}{4}\left(
   Q_T(x+y)-iQ_T(x+iy)-Q_T(x-y)+iQ_T(x-iy)\right)$$
one easily gets $(x,Ty)=(Tx,y)$ for all $x,y\in E$, thus $T$ is
self-adjoint. 

\noindent (b) If $T$ is positive, there exists $S\in\mor(E)$ such that
 $T=S^*S$. Then $(x,Tx)=(Sx,Sx)$ is positive for all $x$. Conversely,
 if $(x,Tx)\ge 0$ for all $x$ then $T$ is self-adjoint using (a) and
 there exist positive elements $T_+,T_-$ such that:
 $$ T = T_+-T_- ,\ T_+T_-=T_-T_+=0 $$
It follows that:
\begin{eqnarray*}
  \forall x\in E,\ (x,T_+x) & \ge & (x,T_-x) \\
  \forall z\in E,\ (T_-z,T_+T_-z) & \ge & (T_-z,T_-T_-z) \\
  \forall z\in E,\ (z,(T_-)^3z) & \le & 0 
\end{eqnarray*}
Since $T_-$ is positive, $T_-^3$ is also positive and the last line above implies
$T^3_-=0$. It follows that $T_-=0$ and then $T=T_+\ge 0$. 
\end{proof}

\subsubsection{Orthocompletion}

Recall that for any subset $S$ of $E$, $S^\perp$ is a Hilbert
submodule of $E$. 
It is also worth noticing that any orthogonal submodules: $F\perp G$ of $E$ are direct summands.

\smallskip \noindent The following properties are left to check as an exercise:
\begin{proposition} Let $F,G$ be $A$-submodules of $E$.
  \begin{itemize}
  \item $E^\perp=\{0\}$ and $\{0\}^\perp=E$.
  \item $\ds F \subset G \Rightarrow G^\perp \subset F^\perp$.
  \item $\ds F\subset F^{\perp\perp}$.
  \item If $F\perp G$ and $F\oplus G=E$ then $F^\perp=G$ and
    $G^\perp=F$. In particular $F$ and $G$ are Hilbert submodules. 
  \end{itemize}
\end{proposition}

\begin{definition}
  A Hilbert $A$-submodule $F$ of $E$ is said to be {\it orthocomplemented} if
  $F\oplus F^\perp = E$. 
\end{definition}

\begin{remark}
  A Hilbert submodule is not necessarily orthocomplemented, even if it
  can be topologically complemented. For instance consider $A=C([0,1])$
  and $J=C_0(]0,1])$ as a Hilbert $A$-submodule of $A$. One easily
  check that $J^\perp=\{0\}$, thus $J$ is not orthocomplemented. On
  the other hand: $A=J\oplus \CC$.
\end{remark}

\begin{lemma}
  Let $T\in \mor(E)$. Then 
  \begin{itemize}
  \item $\ker T^*= (\im T)^\perp$
  \item $\overline{\im T}\subset (\ker T^*)^\perp$
  \end{itemize}
\end{lemma}
\noindent The proof is obvious. Note the difference in the second point with the
case of bounded operators on Hilbert spaces (where equality always occurs). Thus, in general,
 $\ker T^*\oplus\overline{\im T}$ is not the whole of $E$. Such a
 situation can occur  when $\overline{\im T}$ is not orthocomplemented.
 
\smallskip \noindent Let us point out that we can have $T^*$ injective without having $\im T$
dense in $E$ (for instance: $T: C[0,1]\to C[0,1], f\mapsto tf$). Nevertheless,
we have:
\begin{theorem}\label{0isolated}
  Let $T\in \mor(E,F)$. The following assertions are equivalent:
  \begin{enumerate}
  \item $\im T$ is closed,
  \item $\im T^*$ is closed,
  \item $0$ is isolated in $\mathop{\mathrm{spec}}(T^*T)$ (or
    $0\not\in \mathop{\mathrm{spec}}(T^*T)$),
  \item $0$ is isolated in $\mathop{\mathrm{spec}}(TT^*)$ (or
    $0\not\in \mathop{\mathrm{spec}}(TT^{*})$), 
  \end{enumerate}
and in that case $\im T$, $\im T^*$ are orthocomplemented.  
\end{theorem}
\noindent Thus, under the assumption of the theorem $\ker T^*\oplus\im T=F$, 
$\ker T\oplus\im T^*=E$. Before proving the theorem, we gather some
technical preliminaries into a lemma: 
\begin{lemma}
  Let $T\in \mor(E,F)$. Then 
  \begin{enumerate}
  \item $T^*T\ge 0$. We set $|T|=\sqrt{T^*T}$.
  \item $\overline{\im T^*}=\overline{\im |T|}=\overline{\im T^*T}$
  \item Assume that $T(E_1)\subset F_1$ for some Hilbert
    submodules $E_1,F_1$. Then $T|_{E_1}\in \mor(E_1,F_1)$.
  \item If $T$ is onto then $TT^*$ is invertible (in $\mor(F)$) and 
   $E =\ker T\oplus\im T^*$. 
  \end{enumerate}
\end{lemma}
\begin{proof}[Proof of the lemma]
  (1) is obvious.\\
  (2) On has $T^*T(E)\subset T^*(F)$. Conversely: 
 $$ T^*= \lim T^*(1/n+TT^*)^{-1}TT^* \ .$$
This is a convergence in norm since: 
 $$ \| T^*(1/n+TT^*)^{-1}TT^*-T^*\|=\|\frac{1}{n}T^*(\frac{1}{n}+TT^*)^{-1}\|=O(1/\sqrt{n}).$$
It follows that $T^*(F)\subset \overline{T^*T(E)}$
and thus $\overline{\im T^*}=\overline{\im T^*T}$. Replacing $T$ by
$|T|$ yields the other equality. \\
  (3) Easy. \\
  (4) By the open mapping theorem, there exists a positive real number
  $k>0$ such that each $y\in F$ has a preimage $x_y$ by $T$ with $\|y\|\ge
  k\|x_y\|$. Using Cauchy-Schwarz for $T^*y$ and $x_y$, we get:
 $$(*)\quad  \|T^*y\|\ge k\|y\| \quad \forall y\in F \ .$$
Recall that in a $C^*$-algebra, the inequality $a^*a\le b^*b$ is
equivalent to:  $\|ax\|\le \|bx\|$ for all $x\in A$. It can be
adapted to Hilbert modules to show that $(*)$ implies $TT^*\ge k^2$ in
$\mor(F)$, so that $TT^*$ is invertible. Then 
$p=T^*(TT^*)^{-1}T$ is an idempotent and $E=\ker p\oplus \im
p$. Moreover $(TT^*)^{-1}T$ is onto from which it follows that $\im
p= \im T^*$. On the other hand, $T^*(TT^*)^{-1}$ is injective, so that $\ker p = \ker T$. 
\end{proof}

\begin{proof}[Proof of the theorem] 
  Let us start with the implication $(1)\Rightarrow (3)$. By point (3) of the lemma
 $S:=(T: E\to TE)\in \mor(E,TE)$ and by  point (4) of the lemma
 $SS^*$ is invertible. Since the
 spectra of $SS^*$ and $S^*S$ coincide outside $0$ and since
 $S^*S=T^*T$, we get (3). \\
 The implication $(4)\Rightarrow (1)$. Consider the functions $f,g:\RR\to \RR$
 defined by $f(0)=g(0)=0$, $f(t)=1,g(t)=1/t$ for $t\not =0$. 
Thus $f$ and $g$ are continuous on the spectrum of
 $TT^*$. Using the equalities $f(t)t=t$ and $tg(t)=f(t)$, we get
 $f(TT^*)TT^*=TT^*$ and $TT^*g(TT^*)=f(TT^*)$ from which we deduce
 $\im f(TT^*)=\im TT^*$. But $f(TT^*)$ is a projector (self-adjoint
 idempotent), hence $\im TT^*$ is closed and orthocomplemented. 
 Using point (2) of the lemma and the inclusion $\im TT^*\subset \im T$,
 yields (1) (and also the orthocomplementability of $\im T$). \\
At this point we have the following equivalences $(1)\Leftrightarrow (3)\Leftrightarrow (4)$. Replacing
$T$ by $T^*$ we get $(2)\Leftrightarrow (3)\Leftrightarrow (4)$.
\end{proof}
\noindent Another result which deserves to be stated is: 
\begin{proposition} Let $H$ be a Hilbert submodule of $E$ and 
$T:E\to F$ a $A$-linear map. 
  \begin{itemize}
  \item $H$ is orthocomplemented if and only if
    $i: H\hookrightarrow  E\in\mor(H,E)$.
  \item $T\in \mor(E,F)$ if and only if the graph of $T$: 
     $$ \{ (x,y)\in E\times F\ | y=Tx \} $$
    is orthocomplemented. 
  \end{itemize}
\end{proposition}

\subsubsection{Partial isometries}
The following easy result is left as an exercise: 
\begin{proposition}(and definition). 
  Let $u\in \mor(E,F)$. The following assertions are equivalent:
  \begin{enumerate}
  \item $u^*u$ is an idempotent,
  \item $uu^*$ is an idempotent,
  \item $u^*=u^*uu^*$,
  \item $u=uu^*u$. 
  \end{enumerate}
$u$ is then called a partial isometry, with initial support 
$I=\im u^*$ and final support $J=\im u$. 
\end{proposition}
\begin{remark}
  If $u$ is a partial isometry, then $\ker u=\ker u^*u$, 
 $\ker u^*=\ker uu^*$, $\im u=\im uu^*$ and $\im u^*=\im u^*u$. 
 In particular $u$ has closed range and $E=\ker u\oplus \im u^*$,
 $F= \ker u^*\oplus \im u$ where the direct sums are orthogonal. 
\end{remark}

\subsubsection{Polar decompositions}

\smallskip \noindent All homomorphisms do not admit a polar decomposition. For instance,
consider: $T\in\mor(C[-1,1])$ defined by $Tf=t.f$ (here $C[-1,1]$ is
regarded as a Hilbert $C[-1,1]$-module). $T$ is self-adjoint and 
$|T|: f\mapsto |t|.f$. The equation $T=u|T|$, $u\in\mor(C[-1,1])$
leads to the constraint $u(1)(t)= \mbox{sign}(t)$, so 
$u(1)\not \in C[-1,1]$ and $u$ does not exist.

\smallskip \noindent  The next result
clarifies the requirements for a polar decomposition to exist:
\begin{theorem}
  Let $T\in\mor(E,F)$ such that $\overline{\im T}$ and
  $\overline{\im T^*}$ are orthocomplemented. Then there exists a
  unique $u\in \mor(E,F)$, vanishing on $\ker T$, such that 
 $$ T = u |T| $$
Moreover, $u$ is a partial isometry with initial support
$\overline{\im T^*}$ and final support $\overline{\im T}$.
\end{theorem}
\begin{proof}
  We first assume that $T$ and $T^*$ have dense range. Setting 
 $u_n=T(1/n+T^*T)^{-1/2}$ we get a bounded sequence ($\|u_n\|\le 1$)
 such that for all $y\in F$, 
$u_n(T^*y)=T(1/n+T^*T)^{-1/2}T^*y\to \sqrt{TT^*}(y)$.
Thus, by density of $\im T^*$, $u_n(x)$ converges for all $x\in
E$. Let $v(x)$ denotes the limit. Replacing above $T$ by $T^*$, we also
have that $u_n^*(y)$ converges for all $y\in F$, which yields
$v\in\mor(E,F)$. A careful computation shows that $u_n|T|-T$ goes to $0$
in norm. Thus $v|T|=T$. The homomorphism $v$ is unique by density of 
$\im |T|$ and unitary since $u_n^*u_n(x)\to x$ for all 
$x\in \im T^*T$, which proves $v^*v=1$ and similarly for $vv^*$. 

\noindent Now consider the general case and set $E_1=\overline{\im T^*}$, 
$F_1=\overline{\im T}$. One applies the first step to the restriction 
$T_1\in\mor(E_1,F_1)$ of $T$, and we call $v_1$ the unitary
constructed. We set $u(x)=v_1(x)$ if $x\in E_1$ and $u(x)=0$ if 
$x\in E_1^\perp=\ker T$. This definition forces the uniqueness, and it is
clear that $u$ is a partial isometry with the claimed initial/final
supports. 
\end{proof}
\begin{remark}
  $u$ is the strong limit of $T(1/n+T^*T)^{-1/2}$.
\end{remark}

\subsubsection{Compact homomorphisms}

Let $x\in E,y\in F$ and define $\theta_{y,x}\in\mor(E,F)$ by 
 $$ \theta_{y,x}(z)=y.(x,z) \ .$$
The adjoint is given by $\theta_{y,x}^*=\theta_{x,y}$. Then 
\begin{definition} We define 
  $\cK(E,F)$ to be the closure of the linear span of 
$\{ \theta_{y,x}; \ x\in E,y\in F\}$ in $\mor(E,F)$.  
\end{definition}
One easily checks that 
\begin{itemize}
\item $\|\theta_{y,x}\|\le\|x\|\|y\|$ and $\|\theta_{x,x}\|=\|x\|^2$,
\item $T\theta_{y,x}=\theta_{Ty,x}$ and $\theta_{y,x}S=\theta_{y,S^*x}$,
\item $\cK(E):=\cK(E,E)$ is a closed two-sided ideal of $\mor(E)$ (and
  hence a $C^*$-algebra).
\end{itemize}
We also prove: 
\begin{proposition}
  $$ \cM(\cK(E))\simeq \mor(E) $$
where $\cM(A)$ denotes the multiplier algebra of a $C^*$-algebra $A$. 
\end{proposition}
\begin{proof}
 One can show that for all $x\in E$ there is a unique $y\in E$ such that
$x=y.<y,y>$ (a technical exercise: show that the limit $y=\lim x.f_n(\sqrt{(x,x)})$ with 
$f_n(t)=t^{1/3}.(1/n+t)^{-1}$ exists and satisfies the desired assertion).

\smallskip \noindent Consequently, $E$ is a non degenerate $\cK(E)$-module
(ie, $\cK(E).E=E$), indeed $x=y.<y,y>=\theta_{y,y}(y)$.
Using an approximate unit $(u_\lambda)_\Lambda$ for $\cK(E)$, we can
extend the $\cK(E)$-module structure of $E$ into a
$\cM(\cK(E))$-module structure: 
 $$ \forall T\in\cM(\cK(E)), x\in E,\quad T.x = \lim_\Lambda T(u_\lambda).x$$
The existence of the limit is a consequence of $x=\theta_{y,y}(y)$
and $T(u_\lambda).\theta_{y,y}=T(u_\lambda\theta_{y,y})\to T(\theta_{y,y})$.
The limit is $T(\theta_{y,y}).y$. By the uniqueness of $y$,
this module structure, extending that of $\cK(E)$ is unique. 

\smallskip \noindent Hence each $m\in\cM(\cK(E))$ gives rise to a map
$M: E\to E$.  For any $x,z$  in $E$, 
 $$
 (z,M.x)=(z,(m\theta_{y,y}).y)=((m\theta_{y,y})^*(z),y)$$
thus $M$ has an adjoint: $M\in\mor(E)$ and $M^*$ corresponds to
$m^*$. The map $\rho: m \to M$ provides a $*$-homomorphism from
$\cM(\cK(E))$ to $\mor(E)$ which is the identity on $\cK(E)$. On the
other hand let $\pi:\mor(E)\to\cM(\cK(E)$ be the unique $*$-homomorphism,
equal to identity on $\cK(E)$, associated to the inclusion $\cK(E)\subset \mor(E)$ as a closed
ideal. We have $\pi\circ\rho=Id$, and by unicity of the
$\cM(\cK(E))$-module structure of $E$, $\rho\circ\pi=Id$.
\end{proof}
\noindent Let us  give some generic examples:
\begin{enumerate}
\item Consider $A$ as a Hilbert  $A$-module. We know that for any $a\in A$, there
  exists $c\in A$ such that $a=cc^*c$. It follows that the
  map $\gamma_a: A\to A, b\mapsto ab$ is equal to $\theta_{c,c^*c}$
  and thus is compact. We get a $*$-homomorphism 
$\gamma: A\to \cK(A), a\mapsto \gamma_a$ which has dense image (the
linear span of the $\theta$'s is dense in $\cK(A)$) and clearly
injective, because $yb=0$ for all $b\in A$ implies $y=0$. Thus
$\gamma$ is an isomorphism; 
 $$ \cK(A)\simeq A \ .$$
In particular, $\mor(A)\simeq \cM(A)$, and if $1\in A$, then 
$ A\simeq \mor(A)=\cK(A)$.

\item For any $n$, one has in a similar way $\cK(A^n)\simeq M_n(A)$
  and $\mor(A^n) \simeq M_n(\cM(A))$. If moreover $1\in A$,
 $$ 
   (i)\ \mor(A^n)= \cK(A^n)\simeq M_n(A) \ .
 $$
For any  Hilbert $A$-module $E$, we also have $\cK(E^n)\simeq M_n(\cK(E))$.
\end{enumerate}
Relations $(i)$ can be extended to arbitrary finitely generated
 Hilbert $A$-modules: 

\begin{proposition}\label{K-M}
  Let $A$ be a unital $C^*$-algebra and $E$ a $A$-Hilbert module. Then
  the following are equivalent:
  \begin{enumerate}
  \item $E$ is finitely generated. 
  \item $\cK(E)=\mor(E)$. 
  \item $Id_E$ is compact. 
  \end{enumerate}
In that case, $E$ is also projective (ie, it is a direct summand of
$A^n$ for some $n$).
\end{proposition}
For the proof we refer to \cite{WO}. 

\subsection{Generalized Fredholm operators}
\noindent Atkinson's theorem claims that for any bounded linear
operator on a Hilbert
space $H$, the assertion: \\
\centerline{$\ker F$ and $\ker F^*$ are finite dimensional,} \\
is equivalent to: \\
\centerline{there exists a linear bounded operator $G$ such that
  $FG-\id, GF-\id$ are compact .}

\noindent This situation is a little more subtle on Hilbert $A$-modules, since
first of all the kernel of homomorphisms are $A$-modules which are non necessarily
free and secondly, replacing the condition ``finite
dimensional'' by ``finitely generated'', is not enough to recover the
previous equivalence. This is why one uses the second assertion as a definition
of Fredholm operator in the context of Hilbert modules,  and we will
see how to adapt  Atkison's classical result to this new setup. 

\begin{definition} The homomorphism 
  $T\in\mor(E,F)$ is a {\it generalized Fredholm operator} if there exists
  $G\in \mor(F,E)$ such that:
 $$ 
   GF - \id \in\cK(E) \quad \text{ and }\quad FG-\id \in \cK(F) \ .
 $$
\end{definition}

\noindent The following theorem is important to understand the next
chapter on $KK$-theory.
\begin{theorem}\label{comp-pert-with-closed-range}
   Let $A$ be a unital $C^*$-algebra,  $\cE$ a countably
   generated Hilbert $A$-module and $F$ a generalized Fredholm
     operator on $\cE$.
 \begin{enumerate}
   \item If  $\im F$ is closed, then $\ker F$ and $\ker F^*$ are
     finitely generated Hilbert modules. 
   \item  There exists a compact perturbation $G$ of $F$ such
     that $\im G$ is closed. 
 \end{enumerate}
\end{theorem}
\begin{proof}
  (1) Since $\im F$ is closed, so is $\im F^*$ and both are
  orthocomplemented by, respectively, $\ker F^*$ and $\ker F$. Let 
  $P\in\mor(\cE)$ be the orthogonal projection on $\ker F$. Since
  $F$ is a generalized Fredholm operator, there exists $G\in\mor(\cE)$
  such that $Q=1-GF$ is compact. In particular, $Q$ is
  equal to $\id$ on $\ker F$ and: 
  $$
    QP: \cE= \ker F\oplus \im F^* \to \cE, \ x\oplus y \mapsto
    x\oplus 0.
  $$
  Since $QP$ is compact, its restriction: $QP|_{\ker F}:\ker F\to\ker F$  
 is also compact, but $QP|_{\ker F}=\id_{\ker F}$ hence Proposition
 \ref{K-M} implies that $\ker F$ is finitely
  generated. The same argument works for $\ker F^*$. \\
 (2) Let us denote by $\pi$ the projection homomorphism: 
 $$\pi:\mor(\cE)\to C(\cE):=\mor(\cE)/\cK(\cE)\ .$$ Since $\pi(F)$ is
 invertible in $C(\cE)$ it has a polar decomposition:
 $\pi(F)=\omega.|\pi(F)|$. Any unitary of
 $C(\cE)$ can be lifted to a partial isometry of $\mor(\cE)$ \cite{WO}. Let $U$ be such
 a lift of the unitary $\omega$. Using
 $|\pi(F)|=\pi(|F|)$, it follows that:
   $$ F = U|F| \mod \cK(\cE) \ . $$  
 Since $\pi(|F|)$ is also invertible, and positive, we can form 
$\log(\pi(|F|))$ and choose a self-adjoint $H\in\mor(\cE)$ with 
$\pi(H)=\log(\pi(|F|))$. Then:
 $$ \pi(Ue^H)= \omega\pi(|F|) = \pi(F)$$
that is, $Ue^H$ is a compact perturbation of $F$ (and thus is a 
generalized Fredholm operator). $U$ is a partial isometry, hence has a closed image,
and $e^H$ is invertible in $\mor(\cE)$, hence $Ue^H$ has closed image
and the theorem is proved. 
\end{proof}

\subsection{Tensor products}

\subsubsection{Inner tensor products}
Let $E$ be a Hilbert $A$-module, $F$ a Hilbert $B$-module and $\pi:
A\to\mor(F)$ a $*$-homomorphism. We define a sesquilinear
form on $E\otimes_A F$ by setting: 
 $$ 
  \forall x,x'\in E, y,y'\in F,\quad (x\otimes y, x'\otimes
  y')_{E\otimes F}:= (y, (x, x')_{E}\cdot y')_F
 $$
where we have set $a\cdot y= \pi(a)(y)$ to lighten the formula. 
This sesquilinear form is a $B$-valued scalar
product: only the positivity axiom needs some explanation. Set: 
 $$ b = (\sum_i x_i\otimes y_i,\sum_i x_i\otimes y_i) =
 \sum_{i,j}(y_i,(x_i,x_j).y_j)
 $$
where $\pi$ has been omitted. Let us set 
$P=((x_i,x_j))_{i,j}\in M_n(A)$. The matrix $P$ provides a
(self-adjoint) compact homomorphism of $A^n$, which is positive since: 
 $$ \forall a\in A^n, (a,Pa)_{A^n}=\sum_{i,j}a_i^*(x_i,x_j)a_j=(\sum_i
 x_ia_i,\sum_j x_ja_j)\ge 0 \ .
 $$
This means that $P=Q^*Q$ for some $Q\in M_n(A)$. On the other hand,
one can consider $P$ as a homomorphism on $F^n$ and setting 
 $y=(y_1,\ldots,y_n)\in F^n$ we have: 
 $$ b = (y,Py)=(Qy,Qy)\ge 0 \ .$$
Thus $E\otimes_A F$ is a pre-Hilbert module in the generalized sense
(i.e. we do not require the inner product to be definite) and the
Hausdorff completion of $E\otimes_A F$ is a Hilbert $B$-module denoted
in the same way.

\begin{proposition} Let $T\in\mor(E)$ and  $S\in \mor(F)$.
  \begin{itemize}
  \item
 $T\otimes 1: x\otimes y\mapsto Tx\otimes y$ defines a homomorphism of
 $E\otimes_A F$. 
\item  If $S$ commutes with $\pi$ then
  $1\otimes S: x\otimes y\mapsto x\otimes Sy$ is a homomorphism which
  commutes with any $T\otimes 1$.
  \end{itemize}
\end{proposition}
\begin{remark}
 {\bf 1.} Even if $T$ is compact, $T\otimes 1$ is not compact in
    general. The same is true for $1\otimes S$ when defined. \\
 {\bf 2.} In general $1\otimes S$ is not even defined.
\end{remark}

\subsubsection{Outer tensor products}
Now forget the homomorphism $\pi$ and consider the tensor product over
$\CC$ of $E$ and $F$. We set: 
  $$ 
  \forall x,x'\in E, y,y'\in F,\quad (x\otimes y, x'\otimes
  y')_{E\otimes F}:= (x,x')_E\otimes (y, y')_F\in A\otimes B \ .
 $$ 
This defines a pre-Hilbert $A\otimes B$-module in the generalized
sense (the proof of positivity uses similar arguments), where
$A\otimes B$ denotes the spatial tensor product (as it will be the
case in the following, when not otherwise specified). The
Hausdorff completion will be denoted $E\otimes_\CC F$. 

\begin{examples}
  Let $H$ be a separable Hilbert space. Then: 
 $$ H\otimes_\CC A \simeq H_A$$
\end{examples}

\subsubsection{Connections}
We turn back to internal tensor products. We keep
notations of the corresponding subsection. A. Connes and G.
Skandalis \cite{CoS} introduced the notion of connection to bypass
in general the non existence of $1\otimes S$. 
\begin{definition}
Consider two $C^{*}$-algebras $A$ and $B$.  Let  $E$ be a Hilbert $A$-module and 
$F$ be a Hilbert $B$-module.
Assume there is a $*$-morphism 
\begin{eqnarray*}
\pi: A \rightarrow \mathcal{L} (F)
\end{eqnarray*}
and take the inner tensor product $E \otimes_{A} F$.
Given $x \in  E$ we define a homomorphism
\begin{eqnarray*}
T_{x}:   E & \rightarrow &  E \otimes_{A} F \\
 y &\mapsto & x\otimes y
\end{eqnarray*}
whose adjoint is given by 
\begin{eqnarray*}
T_{x}^{*}:   E\otimes_{A} F & \rightarrow &  F \\
 z\otimes y &\mapsto & \pi((x, z))  y\ . 
\end{eqnarray*}
If $S \in  \mathcal{L}(F)$, an {\it $S$-connection} on $E\otimes_{A} F$ is given by an element
\begin{eqnarray*}
G &\in &  \mathcal{L}(E\otimes_{A} F) 
\end{eqnarray*}
such that for all $x  \in  E$: 
\begin{eqnarray*}
T_{x} S - GT_{x } &\in &  \mathcal{K}( F,  E\otimes_{A} F ) \\
S T_{x} ^{*} - T_{x }^{*} G&\in &  \mathcal{K}( E\otimes_{A} F,  F )
\ .
\end{eqnarray*}

\end{definition}

\begin{proposition}
  \begin{itemize}
  \item[(1)] If $[\pi,S]\subset\cK(F)$ then there are $S$-connections. 
  \item[(2)] If $G_i$, $i=1,2$ are $S_i$-connections, then $G_1+G_2$ is a 
    $S_1+S_2$-connection and $G_1G_2$ is a $S_1S_2$-connection. 
  \item[(3)] For any $S$-connection $G$, 
  $[G,\cK(E)\otimes 1]\subset\cK(E\otimes_A F)$.
  \item[(4)] The space of $0$-connections is exactly:
 $$ \{ G\in\mor(F,E\otimes_A F)\ | \ (\cK(E)\otimes 1)G \text{ and }
 G(\cK(E)\otimes 1) \text{ are subsets of } \cK(E\otimes_A F)\} 
 $$
  \end{itemize}
\end{proposition}
\noindent All these assertions are important for the construction of the
Kasparov product. For the proof, see \cite{CoS}

\section{KK-Theory}\label{finKK}
\subsection{Kasparov modules and homotopies}
Given two $C^{*}$-algebras $A$ and $B$  a {\it Kasparov
  $A$-$B$-module} (abbreviated ``Kasparov module'') is given by a triple 
\begin{eqnarray*}
x &=& (\mathcal{E},\pi,F) 
\end{eqnarray*}
where  $\mathcal{E}=\mathcal{E}^{0} \bigoplus \mathcal{E}^{1}$ is a
$(\mathbb{Z}/2\mathbb{Z})$-graded countably generated Hilbert
$B$-module,  
$\pi: A \rightarrow \mathcal{L}(\mathcal{E})$ is a $*$-morphism of
degree $0$ with respect to the grading, and $F \in
\mathcal{L}(\mathcal{E})$ is of degree $1$. These data are required to satisfy the following properties:
\begin{eqnarray*}\label{K-axioms}
\pi (a) (F^{2}-1) &\in&  \mathcal{K}(\mathcal{E})  \qquad \text{for all } a\in A  \\
\left[ \pi (a),F \right] &\in&  \mathcal{K}(\mathcal{E})  \qquad \text{for all } a\in A.  
\end{eqnarray*}
We denote the set of Kasparov $A$-$B$-modules by $E(A,B)$.

Let us immediately define the equivalence
relation leading to $KK$-groups. We denote $B(\left[ 0,1 \right] ):=C(\left[ 0,1 \right],B)$. 

\begin{definition}
A {\it homotopy} between  two Kasparov $A$-$B$-modules
$x=(\mathcal{E},\pi,F)$ and $x'=(\mathcal{E}',\pi',F')$ is a Kasparov
$A$-$B(\left[0,1 \right])$-module $\tilde{x}$ such that: 
\begin{eqnarray}\label{homotopy-axioms}
(ev_{t=0})_{*}(\tilde{x})&=& x,\\
(ev_{t=1})_{*}(\tilde{x})&=& x'.\nonumber
\end{eqnarray}
Here $ev_{t=\cdot}$ is the evaluation map at $t=\cdot$. Homotopy
between Kasparov $A$-$B$-modules is an equivalence 
relation. If there exists a homotopy between $x$ and $x'$ we write
$x\sim_{h}x'.$  
\\ The set of homotopy classes of Kasparov $A$-$B$-modules is denoted $KK(A,B)$.
\end{definition}

There is a natural {\it sum} on $E(A,B)$:
if $x=(\mathcal{E},\pi,F)$ and $x'=(\mathcal{E}',\pi',F')$ belong to $E(A,B)$,
their sum $x+x' \in E(A,B)$ is defined by
\begin{eqnarray*}
x+x' &=& (\mathcal{E}\oplus \mathcal{E}' , \pi \oplus \pi' , F\oplus F').
\end{eqnarray*}

A Kasparov $A-B$-module $x=(\mathcal{E},\pi,F)$ is called {\it degenerate}
if for all $a\in A$, $\pi (a) (F^{2}-1)=0$ and
$\left[ \pi (a),F \right] =0$. It follows:

\begin{proposition}
  Degenerate elements of $E(A,B)$ are homotopic to
  $(0,0,0)$.\\  The sum of Kasparov $A-B$-modules provides  $KK(A,B)$
  with a
  structure of abelian group. 
\end{proposition}
\begin{proof}
  Let $x=(\mathcal{E},\pi,F)\in E(A,B)$ be a degenerate element. Set
 $\tilde{x}= (\tilde{\cE},\tilde{\pi},\tilde{F}) \in E(A,B(\left[ 0,1 \right] ))$ with  
\begin{eqnarray*}
 \tilde{\cE} &=& C_0(\left[ 0,1 \right[,\mathcal{E})\\
\tilde{\pi}(a)\xi(t)&=& \pi(a)\xi(t),\\
\tilde{F}\xi(t)&=& F \xi(t).
\end{eqnarray*}
Then $\tilde{x}$ is a homotopy between $x$ and $(0,0,0)$.

One can easily show that the sum of Kasparov modules makes sense at the level of their
homotopy classes. Thus $KK(A,B)$ admits a commutative
semi-group structure with $(0,0,0)$ as a neutral element. Finally,
the opposite in $KK(A,B)$ of $x=(\mathcal{E},\pi,F)\in E(A,B)$ may be
represented by:
\begin{eqnarray*}
(\mathcal{E}^{op},\pi,-F).
\end{eqnarray*}
where $\cE^{op}$ is $\cE$ with the opposite graduation:
$(\cE^{op})^i=\cE^{1-i}$. 
Indeed,  the module 
$(\mathcal{E},\pi,F)\oplus (\mathcal{E}^{op},\pi,-F)$ is homotopically equivalent to the degenerate module 
\begin{eqnarray*}
(\mathcal{E}\oplus \mathcal{E}^{op},\pi\oplus \pi, \left(
\begin{array}{cc}
  0 & \id \\
  \id & 0 \\
\end{array}
\right)
)
\end{eqnarray*}
This can be realized with the homotopy 
\begin{eqnarray*}
G_{t} &=& \cos (\frac{\pi t}{2})\left(
\begin{array}{cc}
  F & 0 \\
  0 & -F \\
\end{array}
\right) +  \sin (\frac{\pi t}{2})\left(
\begin{array}{cc}
  0 & \id\\
  \id & 0 \\
\end{array}
\right)
\end{eqnarray*}
\end{proof}

\subsection{Operations on Kasparov modules}
\noindent Let us explain the functoriality of $KK$-groups with respect to
its variables. The following two operations on Kasparov modules make
sense on $KK$-groups:

\begin{itemize}
\item{\bf Pushforward along $*$-morphisms: covariance in the second variable.}

Let $x=(\mathcal{E},\pi,F)\in E(A,B)$ and let $g:B\rightarrow C$ be a  $*$-morphism. We define an element $g_{*} (x)\in E(A,C)$ by
\begin{eqnarray*}
g_{*} (x) &=& (\mathcal{E}\otimes_{g}C ,\pi \otimes 1, F\otimes \id),
\end{eqnarray*}
where $\cE\otimes_{g}C$ is the inner tensor product of the 
Hilbert $B$-module  $\cE$ with the Hilbert $C$-module $C$ endowed with
the left action of $B$ given by $g$.

\item{\bf Pullback along $*$-morphisms: contravariance in the first variable.}

Let $x=(\mathcal{E},\pi,F)\in E(A,B)$ and let $f:C\rightarrow A$ be a  $*$-morphism. We define an element 
$f^{*} (x)\in E(C,B)$ by
\begin{eqnarray*}
f^{*} (x) &=& (\mathcal{E} , \pi \circ f , F).
\end{eqnarray*}
\end{itemize}

\noindent Provided with these operations, $KK$-theory is a bifunctor from the
category (of pairs) of $C^{*}$-algebras to the
category of abelian groups.

We recall another useful operation in $KK$-theory: 

\begin{itemize}
\item {\bf Suspension}:

Let $x=(\mathcal{E},\pi,F)\in E(A,B)$ and let $D$ be a  $C^{*}$-algebra. We define an element 
$\tau_{D} (x)\in E(A\otimes D,B\otimes D)$ by
\begin{eqnarray*}
\tau_{D} (x) &=& (\mathcal{E}\otimes_{\mathbb{C}}D , \pi \otimes 1, F\otimes id).
\end{eqnarray*}
Here we take the external tensor product
$\mathcal{E}\otimes_{\mathbb{C}}D$, which is a $B\otimes D$-Hilbert
module. 
\end{itemize}

\subsection{Examples of Kasparov modules and of homotopies between them}

\subsubsection{Kasparov modules coming from homomorphisms between $C^*$-algebras} 

Let $A,B$ be two $C^*$-algebras and $f: A\to B$ a
$*$-homomorphism. Since $\cK(B)\simeq B$, the following:
 $$ [f]:=(B,f,0)$$ 
defines a Kasparov $A-B$-module. If $A$ and $B$ are $\ZZ_{2}$-graded,
$f$ has to be a homomorphism of degree $0$ (ie, respecting the grading).

\subsubsection{Atiyah's Ell}\label{atiyah-Ell}

 Let $X$ be a compact Hausdorff topological space. Take $A=C(X)$ be 
 the algebra of continuous functions on $X$ and let
 $B=\mathbb{C}$. Then  
\begin{eqnarray*}
E(A,B)=Ell(X)
\end{eqnarray*}
the ring of generalized elliptic operators on $X$ as defined by
M. Atiyah. Below we give two concrete examples of such Kasparov
modules:
\begin{itemize}
\item  Assume $X$ is a compact smooth manifold, let $A=C(X)$ as above
  and let $B=\mathbb{C}$ . Let $E$ and $E'$ be  two smooth vector
  bundles over $X$ and denote by $\pi$ the action of $A=C(X)$ by
  multiplication on $L^{2}(X,E)\oplus L^{2}(X,E')$. Given a zero order
  elliptic pseudo-differential operator  
\begin{eqnarray*}
P: C^{\infty}(E) \rightarrow C^{\infty}(E')
\end{eqnarray*}
with parametrix $Q: C^{\infty}(E') \rightarrow C^{\infty}(E)$  the triple
\begin{eqnarray*}
x_{P} &=& \left( L^{2}(X,E)\oplus L^{2}(X,E'), \; \pi, \;  \left(
\begin{array}{cc}
  0 & Q \\
  P & 0 \\
\end{array}
\right)
\right)
\end{eqnarray*}
defines an element in  $E(A,B)=E(C(X), \mathbb{C})$.   
\item  Let $X$ be a compact spin$^c$ manifold  of dimension $2n$, let
  $A=C(X)$ be as above and let $B=\mathbb{C}$. Denote by 
 $S=S^{+}\oplus S^{-}$ the complex spin bundle over  
$X$ and let 
\begin{eqnarray*}
D\mkern -11.5mu /: L^{2}(X,S) \rightarrow L^{2}(X,S)
\end{eqnarray*}
be the corresponding Dirac operator. Let $\pi$ be the action of $A=C(X)$
by multiplication on $L^{2}(X,S)$. Then, the triple
\begin{eqnarray*}
x_{D\mkern -11.5mu /} &=& \left(  L^{2}(X,S), \; \pi, \; \frac{D\mkern
    -11.5mu /}{\sqrt{1+ D\mkern -11.5mu /\, ^2}}    \right)
\end{eqnarray*}
defines an element in  $E(A,B)=E(C(X), \mathbb{C})$.   
\end{itemize}

\subsubsection{Compact perturbations}
 Let $x=(\mathcal{E},\pi,F)\in E(A,B)$.
  Let  $P \in \mor(\mathcal{E})$ which satisfy:
  \begin{equation}
    \label{quasi-compact-def}
 \forall a\in A, \ \pi(a).P\in \cK(\cE)\text{ and } P. \pi(a)\in \cK(\cE)    
  \end{equation}
Then:
\begin{eqnarray*}
x&\sim_{h}& (\mathcal{E},\pi,F+P).
\end{eqnarray*}
The homotopy is the obvious one: $(\mathcal{E}\otimes C([0, 1]), \pi\otimes\id, F+tP)$. 
In particular, when $B$ is unital,  we can always choose a representative  
 $(\mathcal{E},\pi,G)$ with $\im G$ closed (cf. Theorem
 \ref{comp-pert-with-closed-range}).

 \subsubsection{(Quasi) Self-adjoint representatives}
 There exists a representative $(\mathcal{E},\pi,G)$ of $x=(\mathcal{E},\pi,F)\in E(A,B)$
satisfying:
\begin{equation}\label{extra-axiom-defKK} 
 \pi(a)(G-G^{*})\in\cK(\cE) \ .
\end{equation}
Just take $(\mathcal{E}\otimes C([0, 1]), \pi\otimes\id,F_{t})$ as a
homotopy where 
 $$ F_{t}=(tF^{*}F+1)^{1/2}F(tF^{*}F+1)^{-1/2}$$
Then $G=F_{1}$ satisfies (\ref{extra-axiom-defKK}). Now,
$H=(G+G^{*})/2$ is self-adjoint and $P=(G-G^{*})/2$ satisfies
(\ref{quasi-compact-def}) thus $(\cE, \pi, H)$ is another
representative of  $x$. 

Note that (\ref{extra-axiom-defKK}) is often useful in practice and is added as an axiom in many definitions of
$KK$-theory,  like the original one of Kasparov.  It was observed in
\cite{ska} that it could be omitted.  

\subsubsection{Stabilization and Unitarily equivalent modules}
Any Kasparov module
$(E,\pi,F)\in E(A, B)$ is homotopic to a Kasparov module
$(\widehat{\cH}_{B},\rho, G)$ where
$\widehat{\cH}_{B}=\cH_{B}\oplus\cH_{B}$ is the standard graded
Hilbert $B$-module. Indeed,  add to $(E,\pi, F)$ the degenerate
module $(\widehat{\cH}_{B}, 0, 0)$ and consider a grading preserving isometry
$u:E\oplus\widehat{\cH}_{B}\to\widehat{\cH}_{B}$ provided by Kasparov
stabilization theorem. Then,  set $\widetilde{E}=E\oplus\widehat{\cH}_{B}$,  $\widetilde{F}=F\oplus 0$,
$\widetilde{\pi}=\pi\oplus 0$,  $\rho=u\widetilde{\pi} u^{*}$, $G=u\widetilde{F}u^{*}$ and consider the homotopy:
\begin{equation}\label{stab-kasp-mod}
\left( \widetilde{E}\oplus \widehat{\cH}_{B}, \widetilde{\pi}\oplus\rho,
  \begin{pmatrix}
    \cos (\frac{t\pi}{2}) & -u^{*}\sin (\frac{t\pi}{2}) \\ 
    u\sin (\frac{t\pi}{2}) & \cos (\frac{t\pi}{2})
  \end{pmatrix}
\begin{pmatrix}
  \widetilde{F} & 0 \\
  0 & J
\end{pmatrix}
\begin{pmatrix}
  \cos (\frac{t\pi}{2}) & u^{*}\sin (\frac{t\pi}{2}) \\ 
    -u\sin (\frac{t\pi}{2}) & \cos (\frac{t\pi}{2}) 
\end{pmatrix}
\right)
\end{equation}
between $(E,\pi, F)\oplus (\widehat{\cH}_{B}, 0,
0)=(\widetilde{E},\widetilde{\pi},\widetilde{F})$ and
$(\widehat{\cH}_{B},\rho, G)$.  Above, $J$ denotes the operator 
$
\begin{pmatrix}
  0 & 1 \\ 1 & 0 
\end{pmatrix}
$
defined on $\widehat{\cH}_{B}$. 

On says that two Kasparov modules $(E_{i}, \pi_{i}, F_{i})\in
E(A,B)$,  $i=1, 2$ are unitarly equivalent when there exists
 a grading preserving isometry $v:E_{1}\to E_{2}$ such that: 
$$ vF_{1}v^{*}=F_{2} \hbox{ and } \forall a\in A, \
v\pi_{1}(a)v^{*}-\pi_{2}(a)\in \cK(\cE_{2})$$ 
Unitarily equivalent Kasparov modules are homotopic. Indeed,  one can
replace $(E_{i}, \pi_{i}, F_{i})$, $i=1, 2$,  by homotopically,
equivalent modules $(\widehat{\cH}_{B}, \rho_{i}, G_{i})$,
$i=1,2$. It follows from the construction above that the new modules $(\widehat{\cH}_{B}, \rho_{i}, G_{i})$
remain unitarly equivalent and one adapts 
immediately (\ref{stab-kasp-mod}) into a homotopy between then. 

\subsubsection{Relationship with ordinary $K$-theory}\label{link2K-theory} 

Let $B$ be a unital $C^*$-algebra. A finitely generated ($\ZZ/2\ZZ$-graded)
projective $B$-module $\cE$ is a submodule of some $B^N\oplus B^{N}$ and can then be endowed
with a structure of Hilbert $B$-module. On the other hand, $\id_{\cE}$ is a compact
morphism (prop. \ref{K-M}), thus: 
  $$ (\cE, \iota, 0)\in E(\CC, B) $$
where $\iota$ is just multiplication by complex numbers.  This
provides a group homomorphism $K_{0}(B)\to KK(\CC, B)$. 

Conversely, let $(\cE, 1, F)\in E(\CC, B)$ be any Kasparov
module where we have chosen $F$ with closed range (see
above): $\ker F$ is then a finitely generated $\ZZ/2\ZZ$-graded projective
$B$-module. Consider $\widetilde{\cE}=\{\xi\in C([0, 1], \cE)\ |\
\xi(1)\in\ker F\}$ and $\widetilde{F}(\xi):t\mapsto F(\xi(t))$. The
triple $(\widetilde{\cE}, 1, \widetilde{F})$ provides a homotopy between
$(\cE, 1, F)$ and $(\ker F, 1, 0)$. This also gives an inverse of the
previous group homomorphism. 

\subsubsection{A non trivial generator of $KK(\CC, \CC)$}
In the special case $B=\CC$,  we get $KK(\CC, \CC)\simeq
K_{0}(\CC)\simeq\ZZ$ and under this isomorphism,  the following
triple: 
\begin{equation}
  \label{generator-KK-C-C}
 \left( L^{2}(\RR)^{2}, 1, \frac{1}{\sqrt{1+H}}
  \begin{pmatrix}
    0 & -\partial_{x}+x\\ \partial_{x}+x & 0
  \end{pmatrix}
\right) \quad \text{where}\quad H=-\partial_{x}^{2}+x^{2}
\end{equation}
corresponds to $+1$. The reader can check as an exercise that
$\partial_{x}+x$ and  $H$ are essentially self-adjoint as unbounded
operators on $L^{2}(\RR)$,  that $H$ has a compact resolvant and that 
$\partial_{x}+x$ has
a Fredholm index equal to $+1$.  It follows that the Kasparov module
in (\ref{generator-KK-C-C}) is well defined and satisfies the required claim.

\subsection{Ungraded Kasparov modules and $KK_{1}$}

Triple $(\cE,\pi,F)$ satisfying  properties (\ref{K-axioms})
can arise with no natural grading for $\cE$, and consequently with no
diagonal/antidiagonal decompositions for $\pi,F$. 
We refer to those as ungraded Kasparov
$A$-$B$-modules and  the corresponding set is denoted by
$E^{1}(A,B)$. The direct sum is defined in the same way, as well as the
homotopy, which is this times an element of $E^1(A,B[0,1])$.
The homotopy defines an equivalence relation on $E^{1}(A, B)$ and the
quotient inherits a structure of abelian
group as before.   

\smallskip \noindent Let $C_1$ be the complex Clifford Algebra of the
vector space $\CC$ provided with the obvious quadratic form \cite{LM1989}. It
is the $C^*$-algebra $\CC\oplus \varepsilon\CC$ generated by $\varepsilon$ satisfying
$\varepsilon^*=\varepsilon$ and $\varepsilon^2=1$.  Assigning to
$\varepsilon$ the degree $1$ yields a $\ZZ/2\ZZ$-grading on $C_{1}$. We have:
\begin{proposition}
The following map:
\begin{equation}
  \label{eq:ungraded2graded}
\begin{matrix}
  E^1(A,B) & \longrightarrow & E(A,B\otimes C_1) \\
   (\cE,\pi,F) & \longmapsto & (\cE\otimes C_1,\pi\otimes \id ,F\otimes
   \varepsilon)
\end{matrix}
\end{equation}
induces an isomorphism between the quotient of $E^1(A,B)$ under homotopy and $KK_{1}(A,B)= KK(A,B\otimes C_1)$. 
\end{proposition}
\begin{proof}
The grading of $C_{1}$ gives the one of $\cE\otimes C_{1}$ and the
map (\ref{eq:ungraded2graded}) easily gives  a homomorphism $c$ from
$KK_{1}(A,B)$ to $KK(A,B\otimes C_1)$. 

Now let $y=(\cE, \pi, F)\in E(A,B\otimes C_1)$. The multiplication by
$\varepsilon$ on the right of $\cE$ makes sense,  even if $B$ is not
unital,  and one has $\cE_{1}=\cE_{0}\varepsilon$.  It follows that  $\cE = \cE_{0}\oplus\cE_{1}\simeq
\cE_{0}\oplus\cE_{0}$ and any $T\in\mor(\cE)$,  thanks to the $B\otimes
C_{1}$-linearity, has the following expression: 
 $$ T=
 \begin{pmatrix}
   Q & P \\ P & Q
 \end{pmatrix} \quad P, Q\in\mor_{B}(\cE_{0})
$$
Thus $F=
\begin{pmatrix}
  0 & P \\ P & 0
\end{pmatrix}
$,  $\pi=
\begin{pmatrix}
  \pi_{0} & 0 \\ 0 & \pi_{0}
 \end{pmatrix}
$ and $c^{-1}[y]=[\cE_{0}, \pi_{0}, P]$. 
\end{proof}

\begin{remark}
  The opposite of $(\cE,\pi,F)$ in $KK_{1}(A,B)$ is represented by
  $(\cE,\pi,-F)$. One may wonder why we have to decide if a Kasparov
  module is graded or not. Actually, If we forget the $\ZZ/2\ZZ$ grading of a graded Kasparov
  $A-B$-module $x=(\cE,\pi,F)$ and  consider it as an
  ungraded module, then we get the trivial class in
  $KK_{1}(A,B)$. Let us prove this claim.  

  The grading of $x$ implies that $\cE$ has
  a decomposition $\cE=\cE_{0}\oplus\cE_{1}$ for which $F$ has degree
  $1$,  that is: $F=\begin{pmatrix}0& Q\\P& 0\end{pmatrix}$.  Now:
  \begin{equation}
    \label{graded-elt-trivial-as-ungraded}
    G_{t} = \cos(t\pi/2)F+ \sin(t\pi/2)\begin{pmatrix} 1& 0\\ 0& -1 \end{pmatrix}
  \end{equation}
  provides an homotopy in $KK_{1}$ between $x$ and
  $(\cE,\pi,\begin{pmatrix} 1& 0\\ 0& -1 \end{pmatrix})$.  Since the
  latter is degenerate,  the claim is proved. 
\end{remark}

\begin{examples}
  Take again the example of the Dirac operator $D\mkern -11.5mu /$ introduced in
  (\ref{atiyah-Ell}) on a spin$^c$ manifold $X$ whose dimension is
  odd. There is no natural $\ZZ/2\ZZ$ grading for the spinor
  bundle. The previous triple $x_{D\mkern -11.5mu /}$ provides this
  time an interesting class in $E^{1}(C(X),\CC)$.
\end{examples}

\subsection{The Kasparov product}
In this section we construct the product
\begin{eqnarray*}
KK(A,B) \otimes KK(B,C) \rightarrow KK(A,C) \ .
\end{eqnarray*}
It satisfies the properties given in Section \ref{intro-KK}. Actually: 
\begin{theorem}
\label{theoconnect}
 Let $x=(\mathcal{E},\pi,F)$ in $E(A,B)$ and
 $x=(\mathcal{E}',\pi',F')$ in $ E(B,C)$ be two Kasparov modules. Set 
 \begin{eqnarray*}
\mathcal{E}'' &=& \mathcal{E} \otimes_{B} \mathcal{E}'  
\end{eqnarray*}
and 
\begin{eqnarray*}
 \pi'' &=& \pi   \otimes 1 
\end{eqnarray*}
Then there exists a unique, up to homotopy, $F'$-connection on $\mathcal{E}''$ denoted by $F''$ such that 
\begin{itemize}
\item $(\mathcal{E}'',\pi'',F'') \in E(A,C)$
\item $\pi''(a) \left[  F'' , F\otimes 1\right] \pi''(a)$ is nonnegative modulo $\mathcal{K}(\mathcal{E}'')$ for all $a\in A$.
\end{itemize}
$(\mathcal{E}'',\pi'',F'')$ is the Kasparov product of $x$ and
$x'$. It enjoys all the properties described in Section \ref{intro-KK}. 
\end{theorem}

\begin{proof}[Idea of the proof]
  
We only explain the construction of the operator 
$F''$. For a complete proof, see for instance \cite{Ka1,CoS}. A very naive idea for $F''$ could be $F\otimes 1 + 1 \otimes F'$ but the
trouble is that the operator $1\otimes F'$ is in general not well defined. We can overcome
this first difficulty by replacing the not well defined
$1\otimes F'$ by any $F'$-connection $G$ on $\mathcal{E}''$, and try 
$F\otimes 1 + G$. We stumble on a second problem,  namely  that the
properties of Kasparov module are not satisfied in general with this candidate for
$F''$: for instance 
$(F^2-1)\otimes 1\in \cK(\cE)\otimes 1\not\subset\cK(\cE'')$ as soon as $\cE''$ is not finitely generated. 

\smallskip \noindent The case of tensor products of elliptic self-adjoint differential
operators on a closed manifold $M$, gives us a hint towards the right way. If
$D_1$ and $D_2$ are two such operators and $H_1$,$H_2$ the natural
$L^2$ spaces on which they act, then the bounded operator on
$H_1\otimes H_2$: 
\begin{equation}
  \label{eq:bad-def}
\frac{D_1}{\sqrt{1+D_1^2}}\otimes 1 + 1\otimes
 \frac{D_2}{\sqrt{1+D_2^2}}
\end{equation}
inherits the same problem as $F\otimes 1 + G$ but: 
$$ D'':= \frac{1}{\sqrt{2+D_1^2\otimes 1+1\otimes D_2^2}}(D_1\otimes 1
+ 1\otimes D_2)
 $$
has better properties: $D''^2-1$ and $[C(M),D'']$ belong to $\cK(H_1\otimes
H_2)$. Note that 
 $$ D'' = \sqrt{M}.\frac{D_1}{\sqrt{1+D_1^2}}\otimes 1 + \sqrt{N}.1\otimes
 \frac{D_2}{\sqrt{1+D_2^2}}$$
with $$M= \frac{1+D_1^2\otimes 1}{2+D_1^2\otimes 1+1\otimes D_2^2}
\mbox{ and } N=\frac{1+1\otimes D_2^2}{2+D_1^2\otimes 1+1\otimes
  D_2^2} \ .$$
The operators $M,N$ are bounded on $H_1\otimes H_2$, positive, and
satisfy $M+N=1$. We thus see that in that case, the naive idea
(\ref{eq:bad-def}) can be corrected by combining the involved
operators with some adequate ``partition of unity''.

Turning back to our problem, this calculation leads us to look for
an adequate operator $F''$ in the following form: 
 $$ F''=\sqrt{M}.F\otimes 1 + \sqrt{N}G \ .$$
We need to have that $F''$
is a $F'$-connection, and satisfies $a.(F''^2-1)\in\cK(E'')$ and 
$[a,F'']\in\cK(E'')$ for all $a\in A$ (by $a$ we mean $\pi''(a)$). Using the previous form for
$F''$, a small computation shows that these assertions become true if
all the following conditions hold:
\begin{itemize}
\item[$(\imath)$] $M$ is a $0$-connection (equivalently, $N$ is a $1$-connection),
\item[$(\imath\imath)$] $[M,F\otimes1],\ N.[F\otimes 1,G],\ [G,M],\ N(G^2-1)$ belong to
  $\cK(E'')$,
\item[$(\imath\imath\imath)$] $[a,M],\ N.[G,a]$ belong to  $\cK(E'')$.
\end{itemize}

At this point there is a miracle:

\begin{theorem}[Kasparov's technical theorem]
Let $J$ be a $C^{*}$-algebra and denote by $ \mathcal{M}(J)$ its multipliers algebra.  Assume there are two subalgebras 
$A_{1}, \, A_{2} $ of $\mathcal{M}(J)$ and a linear subspace $\bigtriangleup  \subset  \mathcal{M}(J)$ such that 
\begin{eqnarray*}
A_{1} \, A_{2} &\subset  &  J,  \\
\left[ \bigtriangleup, A_{1} \right] &\subset &  J. 
\end{eqnarray*}
Then there exist two nonnegative elements $M, N \in   \mathcal{M}(J)$ with $M+N =1$ such that 
\begin{eqnarray*}
M \, A_{1} &\subset  &  J,  \\
N \, A_{2} &\subset  &  J,  \\
\left[ M, \bigtriangleup \right] &\subset &  J.  
\end{eqnarray*}
\end{theorem}
\noindent For a proof, see \cite{Hig2}. 

\noindent Now, to get $(\imath),(\imath\imath), (\imath\imath\imath)$, we apply this theorem with: 
\begin{eqnarray*}
 A_{1} &= &  C^{*} \langle   \mathcal{K}(\mathcal{E}) \otimes 1, \, \mathcal{K}(\mathcal{E}'')  \rangle,   \\
A_{2} &= &  C^{*} \langle   G^{2} - 1, \, \left[  G, F \otimes 1     \right]  ,\,  \left[  G, \pi ''  \right] \rangle,   \\
\bigtriangleup &=& Vect  \langle  \pi'' (A), \, G, \, F\otimes 1 \rangle. 
\end{eqnarray*}

This gives us the correct $F''$. 
\end{proof}

\subsection{Equivalence and duality in $KK$-theory}
With the Kasparov product come the following notions: 
\begin{definition}
  Let $A,B$ be two $C^*$-algebras. 
  \begin{itemize}
  \item One says that $A$ and $B$ are {\it $KK$-equivalent} if there exist
    $\alpha\in KK(A,B)$ and $\beta\in KK(B,A)$ such that: 
    $$ \alpha\otimes\beta= 1_{A}\in KK(A,A) \text{ and }
    \beta\otimes\alpha=1_{B}\in KK(B,B). $$
    In that case, the pair $(\alpha,\beta)$ is called a
    {\it $KK$-equivalence} and it gives rise to isomorphisms 
 $$ KK(A\otimes C,D)\simeq KK(B\otimes C,D) \text{ and } KK(C,A\otimes
 D)\simeq KK(C,B\otimes D) $$
  given by Kasparov products for all $C^*$-algebras $C,D$.
\item One says that $A$ and $B$ are {\it $KK$-dual (or Poincar\'e dual)} if
  there exist $\delta\in KK(A\otimes B,\CC)$ and $\lambda\in KK(\CC,A\otimes B)$ such that: 
    $$ \lambda\underset{B}{\otimes}\delta=1\in KK(A,A) \text{ and
    }\lambda\underset{A}{\otimes}\delta=1\in KK(B,B) \ .$$
  In that case, the pair $(\lambda,\delta)$ is called a
    $KK$-duality and it gives rise to isomorphisms 
 $$ KK(A\otimes C,D)\simeq KK(C,B\otimes D) \text{ and } KK(C,A\otimes
 D)\simeq KK(B\otimes C,B\otimes D) $$
  given by Kasparov products for all $C^*$-algebras $C,D$.
  \end{itemize}
\end{definition}

We continue this paragraph with classical computations illustrating these notions. 

\subsubsection{Bott periodicity}\label{bott-periodicity}
 
Let $\beta\in KK(\mathbb{C},C_{0}(\mathbb{R}^{2}))$
be represented by the Kasparov module:
\begin{eqnarray*}
\left(\cE, \pi, C \right)=\left( C_{0}(\mathbb{R}^{2})\oplus
  C_{0}(\mathbb{R}^{2}), \; 1, \;  \frac{1}{\sqrt{1+c^{2}}}\left(
\begin{array}{cc}
  0 &  c_{-} \\
  c_{+} & 0 \\
\end{array}
\right)
\right).
\end{eqnarray*}
where $c_{+}, c_{-}$ are the operators given by pointwise multiplication by
$x-\imath y$ and $x+\imath y$ respectively and 
$ c=
\begin{pmatrix}
  0 & c_{-}\\ c_{+} & 0
\end{pmatrix}
$. 

Let $\alpha \in KK(C_{0}(\mathbb{R}^{2}),\mathbb{C})$ be represented by the Kasparov module:
\begin{eqnarray*}
\left(\cH, \pi, F \right)=\left( L^{2}(\mathbb{R}^{2})\oplus
  L^{2}(\mathbb{R}^{2}), \; \pi, \;  \frac{1}{\sqrt{1+ D ^{2}}}
  \left(
\begin{array}{cc}
  0 & D_{-} \\
  D_{+} & 0 \\
\end{array}
\right)
\right)
\end{eqnarray*}
where $\pi: C_{0}(\mathbb{R}^{2}) \rightarrow \mathcal{L}(L^{2}(\mathbb{R}^{2})\oplus L^{2}(\mathbb{R}^{2}))$ is the action given by multiplication of functions and the operators $D_{+}$  and $D_{-}$ are given by 
\begin{eqnarray*}
  D_{+} &=&  \partial_{x} + \imath \partial{y} \\
    D_{-}  &=&  -\partial_{x} + \imath \partial{y}.
\end{eqnarray*} and  $ D= \left(
\begin{array}{cc}
  0 & D_{-} \\
  D_{+} & 0 \\
\end{array}
\right)$.
\begin{theorem}
  $\alpha$ and $\beta$ provide a $KK$-equivalence between $C_0(\RR^2)$
  and $\CC$
\end{theorem}
This is the Bott periodicity Theorem in the bivariant $K$-theory framework. 
\begin{proof}
  Let us begin with the computation of
$\beta\otimes\alpha\in KK(\CC,\CC)$. 
We have an identification: 
\begin{equation}
  \label{Bott-EprimetensE}
  \cE\underset{C_{0}(\RR^{2})}{\otimes}\cH \simeq \cH\oplus\cH
\end{equation}
where on the right,  the first copy of $\cH$ stands for
$\cE_{0}\underset{C_{0}(\RR^{2})}{\otimes}\cH_{0}\oplus
\cE_{1}\underset{C_{0}(\RR^{2})}{\otimes}\cH_{1}$ and the second for 
 $\cE_{0}\underset{C_{0}(\RR^{2})}{\otimes}\cH_{1}\oplus\cE_{1}\underset{C_{0}(\RR^{2})}{\otimes}\cH_{0}$. 
One checks directly that under this identification the following
operator
\begin{equation}
  \label{Bott-Fconnection}
  G = \frac{1}{\sqrt{1+D^{2}}}
  \begin{pmatrix}
    0 & 0 & D_{-} & 0 \\
    0 & 0 & 0 &   -D_{+} \\
    D_{+} & 0 & 0 & 0 \\
    0 & -D_{-} & 0 & 0 
  \end{pmatrix}
\end{equation}
is an $F$-connection.  On the other hand,  under the identification
(\ref{Bott-Fconnection}),  the operator $C\otimes 1$ gives:
\begin{equation}
  \label{Bott-Ctens1}
    \frac{1}{\sqrt{1+c^{2}}}
  \begin{pmatrix}
    0 & 0 &0 & c_{-} \\
    0 & 0 & c_{+} &   0 \\
    0 & c_{-} & 0 & 0 \\
    c_{+} & 0 & 0 & 0 
  \end{pmatrix}
\end{equation}
It immediately follows that $\beta\otimes\alpha$ is represented by: 
\begin{equation}
  \label{betaprodalpha}
  \delta=\left(\cH\oplus\cH, 1,
    \frac{1}{\sqrt{1+c^{2}+D^{2}}}\mathbf{D}
\right), 
\end{equation}
where $\mathbf{D}=\begin{pmatrix}
      0 & \mathbf{D}_{-} \\ \mathbf{D}_{+}&0
    \end{pmatrix}$; $\mathbf{D}_{+}=\begin{pmatrix}
      D_{+} & c_{-}\\ c_{+}& -D_{-}\end{pmatrix}$ and $\mathbf{D}_{-}=\mathbf{D}_{+}^{*}$. 
Observe that, denoting by
$\rho$ the rotation in $\RR^{2}$ of angle $\pi/4$,  we have:
\begin{eqnarray*}
\begin{pmatrix}
  \rho^{-1} & 0 \\ 0 & \rho
\end{pmatrix}
\begin{pmatrix}
  0 & \mathbf{D}_{-}\\ \mathbf{D}_{+} & 0 
\end{pmatrix}
\begin{pmatrix}
  \rho & 0 \\ 0 & \rho^{-1}
\end{pmatrix} &=&
\begin{pmatrix}
   0 & \rho^{-1}\mathbf{D}_{-}\rho^{-1}\\ \rho\mathbf{D}_{+}\rho & 0 
\end{pmatrix} \\ &=&
\begin{pmatrix}
  0& 0& \imath(\partial_{y}-y) & -\partial_{x}+x \\
  0 & 0& \partial_{x}+x & -\imath(\partial_{y}+y) \\
  \imath(\partial_{y}+y) & -\partial_{x}+x & 0 & 0 \\
  \partial_{x}+x & \imath(-\partial_{y}+y) & 0 & 0 
\end{pmatrix} \\ 
 &=&
 \begin{pmatrix}
   0 & x-\partial_{x} \\ x+\partial_{x} & 0
 \end{pmatrix}\otimes 1+1\otimes 
\begin{pmatrix}
   0 & \imath(\partial_{y}-y) \\ \imath(\partial_{y}+y)& 0 
 \end{pmatrix}. 
\end{eqnarray*}
Of course 
 $$ \delta \sim_{h} \left(\cH\oplus\cH, 1,
   \frac{1}{\sqrt{1+c^{2}+D^{2}}}
\begin{pmatrix}
   0 & \rho^{-1}\mathbf{D}_{-}\rho^{-1}\\ \rho\mathbf{D}_{+}\rho & 0 
\end{pmatrix}\right)
$$
and the above computation shows that $\delta$ coincides with the
Kasparov product $u\otimes u$ with $u\in KK(\CC, \CC)$ given by: 
 $$ u=\left(L^{2}(\RR)^{2}, 1, \frac{1}{\sqrt{1+x^{2}+\partial_{x}^{2}}}\begin{pmatrix}
   0 & x-\partial_{x} \\ x+\partial_{x} & 0
 \end{pmatrix}\right).  $$
A simple exercise shows that $\partial_{x}+x:L^{2}(\RR)\to
L^{2}(\RR)$ is essentially self-adjoint with one dimensional kernel
and zero dimensional cokernel,  thus $1=u=u\otimes u\in KK(\CC, \CC)$. 

Let us turn to the computation of $\alpha\otimes\beta\in
KK(C_{0}(\RR^{2}),C_{0}(\RR^{2}))$: it is a Kasparov product over
$\CC$,  thus it commutes:
\begin{equation}
  \label{Bott-prodoverC}
  \alpha\otimes\beta=\tau_{C_{0}(\RR^{2})}(\beta)\otimes\tau_{C_{0}(\RR^{2})}(\alpha)
\end{equation}
but we must observe that the two copies of $C_{0}(\RR^{2})$ in
$\tau_{C_{0}(\RR^{2})}(\beta) $ and $\tau_{C_{0}(\RR^{2})}(\alpha) $
play a different r\^ole: on should think of the first copy as
functions of the variable $u\in\RR^{2}$  and of the variable
$v\in\RR^{2}$ for the second.  It follows that one can not directly factorize $\tau_{C_{0}(\RR^{2})}$
 on  the right hand side of (\ref{Bott-prodoverC}) in order to use the
value of $\beta\otimes\alpha$.  This is where a classical argument,
known as the rotation trick of Atiyah,  is necessary:
\begin{lemma}
  Let $\phi: C_{0}(\RR^{2})\otimes C_{0}(\RR^{2})\to
  C_{0}(\RR^{2})\otimes C_{0}(\RR^{2}) $ be the flip automorphism:
$\phi(f)(u,v)=f(v,u)$.
Then: 
  $$[\phi]=1\in
  KK(C_{0}(\RR^{2})\otimes C_{0}(\RR^{2}),C_{0}(\RR^{2})\otimes
  C_{0}(\RR^{2}) ) $$ 
\end{lemma}
\begin{proof}[Proof of the lemma]
  Let us denote by $I_{2}$ the identity matrix of $M_{2}(\RR)$.
  Use a continuous path of isometries of $\RR^{4}$ connecting 
$\begin{pmatrix}
    0 & I_{2} \\ I_{2} & 0 
  \end{pmatrix}
$ to 
$\begin{pmatrix}
    I_{2} & 0 \\
      0 & I_{2}  \end{pmatrix}
$. 
This gives a homotopy $(C_{0}(\RR^{2})\otimes C_{0}(\RR^{2}), \phi,
0)\sim_{h}(C_{0}(\RR^{2})\otimes C_{0}(\RR^{2}), \id, 0) $. 
\end{proof}
Now 
\begin{eqnarray}
  \label{Bott-end}
  \alpha\otimes\beta&=&\tau_{C_{0}(\RR^{2})}(\beta)\otimes\tau_{C_{0}(\RR^{2})}(\alpha)
                    =
                    \tau_{C_{0}(\RR^{2})}(\beta)\otimes[\phi]\otimes\tau_{C_{0}(\RR^{2})}(\alpha)\\
                    &=&
                    \tau_{C_{0}(\RR^{2})}(\beta\otimes\alpha)
                    = \tau_{C_{0}(\RR^{2})}(1)=1\in KK(C_{0}(\RR^{2})
                    , C_{0}(\RR^{2})). \nonumber
\end{eqnarray}  
\end{proof}
\subsubsection{Self duality of $C_{0}(\RR)$} With the same notations
as before,  we get:
\begin{corollary}
 The algebra $C_{0}(\RR)$ is Poincar\'e dual to itself. 
\end{corollary}
Other examples of Poincar\'e dual algebras will be given later. 
\begin{proof}
The automorphism $\psi$ of $C_{0}(\RR)^{\otimes^{3}}$ given by
$\psi(f)(x, y, z)=f(z, x, y)$ is homotopic to the identity thus: 
\begin{eqnarray}
  \label{Self-duality-C0R}
   \beta \underset{C_{0}(\RR)}{\otimes}\alpha
     &=&\tau_{C_{0}(\RR)}(\beta)\otimes\tau_{C_{0}(\RR)}(\alpha)
                    =
                    \tau_{C_{0}(\RR)}(\beta)\otimes[\psi]\otimes\tau_{C_{0}(\RR)}(\alpha)\\
                    &=&
                    \tau_{C_{0}(\RR)}(\beta\otimes\alpha)
                    = \tau_{C_{0}(\RR)}(1)=1\in KK(C_{0}(\RR) ,
                    C_{0}(\RR)).  \nonumber
\end{eqnarray}  
\end{proof}

\begin{exercise} With $C_{1}=\CC\oplus\varepsilon\CC$ the Clifford
  algebra of $\CC$,  consider:
$$ \beta_{c}=\left(C_{0}(\RR)\otimes C_{1}, 1,
  \frac{x}{\sqrt{x^{2}+1}}\otimes\varepsilon \right)\in KK(\CC,
C_{0}(\RR)\otimes C_{1}), $$
$$ \alpha_{c}=\left(L^{2}(\RR, \Lambda^{*}\RR), \pi,
  \frac{1}{\sqrt{1+\Delta}}(d+\delta)\right)\in KK(C_{0}(\RR)\otimes
C_{1}, \CC), $$
where $(d+\delta)(a+bdx)=-b'+a'dx$,  $\Delta=(d+\delta)^{2}$ and $\pi(f\otimes\varepsilon)$
sends $a+bdx$ to $f(b+adx)$. \\
 Show that $\beta_{c}, \alpha_{c}$ provide a $KK$-equivalence
  between $\CC$ and $C_{0}(\RR)\otimes C_{1}$ (Hints: compute
  directly $\beta_{c}\otimes\alpha_{c}$,  then use the commutativity
  of the Kasparov product over $\CC$ and  check that the flip
  of $(C_{0}(\RR)\otimes C_{1})^{\otimes^2}$ is $1$ to conclude about the computation of
  $\alpha_{c}\otimes\beta_{c}$ ).  
\end{exercise}

\subsubsection{A simple Morita equivalence}
Let $\imath_{n}=(M_{1, n}(\CC), 1, 0)\in E(\CC, M_{n}(\CC))$ where
the $M_{n}(\CC)$-module structure is given by multiplication by
matrices on the right. Note that $[\imath_{n}]$ is also the class of
the homomorphism $\CC\to M_{n}(\CC)$ given by the left up
corner inclusion. Let also 
$\jmath_{n}=(M_{n, 1}(\CC), m, 0)\in E(M_{n}(\CC), \CC)$ where $m$ is
 multiplication by matrices on the left. It follows in a
 straightforward way that:
 $$ \imath_{n}\otimes\jmath_{n}\sim_{h}(\CC, 1, 0)\text{ and }
 \jmath_{n}\otimes\imath_{n}\sim_{h}(M_{n}(\CC), 1, 0)$$
thus $\CC$ and $M_{n}(\CC)$ are $KK$-equivalent and this is an
example of a Morita equivalence.  The map in $K$-theory associated
with $\jmath$: 
 $ \cdot\otimes\jmath_{n}: K_{0}(M_{n}(\CC))\to\ZZ$ is just the trace
 homomorphism.  
Similarly,  let us consider the Kasparov elements  $\imath\in E(\CC,
\cK(\cH))$ associated to the homomorphism $\imath: \CC\to \cK(\cH)$ given by the choice of
a rank one projection and $\jmath=(\cH, m, 0)\in E(\cK(\cH), \CC)$
where $m$ is just the action of compact operators on $\cH$: they
provide a $KK$-equivalence between $\cK$ and $\CC$.

\subsubsection{$C_{0}(\RR)$ and $C_{1}$. }
We leave the proof of the following result as an exercise:
\begin{proposition}
  The algebras $C_{0}(\RR)$ and $C_{1}$ are $KK$-equivalent.  
\end{proposition}
\begin{proof}[Hint for the proof:]
Consider 
$$\widetilde{\alpha}=\left(L^{2}(\RR, \Lambda^{*}\RR), m,
  \frac{1}{\sqrt{1+\Delta}}(d+\delta)\right)\in
KK(C_{0}(\RR), C_{1})$$
where $d, \delta, \Delta$ are defined in the previous exercise,  
$m(f)(\xi)=f\xi$,  and the $C_{1}$-right module structure of 
$L^{2}(\RR, \Lambda^{*}\RR)$ is given by
$(a+bdx)\cdot\varepsilon=-ib+iadx$. Consider also: 
$$\widetilde{\beta}=\left(C_{0}(\RR)^{2},\varphi, \frac{x}{\sqrt{1+x^{2}}}
  \begin{pmatrix}
    0&1\\1&0
  \end{pmatrix}
\right)\in
KK(C_{1}, C_{0}(\RR)) $$
where $\varphi(\varepsilon)(f, g)=(-ig, if)$. 
Prove that they provide the desired $KK$-equivalence. 
\end{proof}
\begin{exercise}
\begin{enumerate}
\item Check that $\tau_{C_{1}}:KK(A, B)\to KK(A\otimes C_{1},
  B\otimes C_{1})$ is an isomorphism. 
\item Check that under $\tau_{C_{1}}$ and the
  Morita equivalence $M_{2}(\CC)\sim\CC$,  the elements $\alpha_{c},
  \beta_{c}$ of the previous exercise coincide with
  $\widetilde{\alpha},\widetilde{\beta}$
and recover the $KK$-equivalence between $C_{1}$ and $C_{0}(\RR)$.
\end{enumerate}
 \end{exercise}
 \begin{remark}
   At this point, one sees that 
   $KK_{1}(A,B)=KK(A,B(\RR))$ , ($B(\RR):=C_0(\RR)\otimes B$) can also be presented in the following different ways:
$$  E_1(A, B)/\!\sim_h \simeq KK(A, B\otimes C_{1})\simeq KK(A\otimes C_{1},
B)\simeq KK(A(\RR), B) $$
 \end{remark}
 
\subsection{Computing the Kasparov product without its definition}

Computing the product of two Kasparov modules is in general quite
hard, but we are very often in one of the following situations.  

\subsubsection{Use of the functorial properties}

Thanks to the functorial properties listed in Section \ref{intro-KK}, many products
can be deduced from known, already computed,  ones. For instance,  in
the proof of Bott periodicity (the $KK$-equivalence between
$\CC$ and $C_{0}(\RR^{2})$) one had to compute two products:
the first one was directly computed,  the second one was deduced from
the first using the properties of the Kasparov product and a simple
geometric fact.  There are numerous examples of this kind. 

\subsubsection{Maps between $K$-theory groups}

Let $A, B$ be two unital (if not,  add a unit) $C^{*}$-algebras,
$x\in KK(A, B)$ be given by a Kasparov module $(\cE, \pi, F)$ where
$F$ has a closed range and
assume that we are interested in the map $\phi_{x}:K_{0}(A)\to
K_{0}(B)$ associated with $x$ in the following way: 
$$ y\in K_{0}(A)\simeq KK(\CC, A);\quad \phi_{x}(y)=y\otimes x$$
This product takes a particularly simple form when $y$ is represented
by $(\cP, 1, 0)$ with $\cP$ a finitely generated projective
$A$-module (see \ref{link2K-theory}): 
 $$ y\otimes x = \left(\cP\underset{A}{\otimes}\cE, 1\otimes \pi,
 \id\otimes F\right)=\left(\ker(\id\otimes F), 1, 0\right). $$

\subsubsection{Kasparov elements constructed from homomorphisms}

Sometimes,  Kasparov classes $y\in KK(B, C)$ can be explicitly
represented as Kasparov products of classes of homomorphisms with inverses
of such classes. Assume for instance that
$y=[e_{0}]^{-1}\otimes[e_{1}]$ where $e_{0}:\cC\to B$,  $e_{1}:\cC\to
C$ are  homomorphisms   of $C^{*}$-algebras and $e_{0}$ produces an
invertible element in $KK$-theory (for instance: $\ker e_{0}$ is
$K$-contractible and: $B$ is nuclear or $\cC, B$ $K$-nuclear, see \cite{ska-88,cun-ska}). 
Then computing a Kasparov product $x\otimes y$ where $x\in KK(A, B)$
amounts to lifting $x$ to $KK(A, \cC)$,  that is to finding $x'\in KK(A, \cC)$
such that $(e_{0})_{*}(x')=x$ and restrict this lift to $KK(A, C)$,
that is evaluate $x"=(e_{1})_{*}(x')$.  It follows from the
properties of the product that $x"=x\otimes y$.  
\begin{examples}
  Consider the tangent groupoid $\cG_{\RR}$ of $\RR$ and let
  $\delta=[e_{0}]^{-1}\otimes[e_{1}]\otimes\mu$ be the associated
  deformation element: $e_{0}:C^{*}(\cG_{\RR})\to C^{*}(T\RR)\simeq
  C_{0}(\RR^{2})$ is evaluation at $t=0$,  $e_{1}:C^{*}(\cG_{\RR})\to
  C^{*}(\RR\times\RR)\simeq\cK(L^{2}(\RR))\simeq\cK$ is evaluation at $t=1$
  and $\mu=(L^{2}(\RR),m, 0)\in KK(\cK, \CC)$ gives the Morita
  equivalence $\cK\sim\CC$.  

  Let $\beta\in KK(\CC, C_{0}(\RR^{2}))$ be the element used in
  paragraph \ref{bott-periodicity}. Then $\beta\otimes\delta\in KK(\CC,\CC)$ is easy
  to compute. The lift $\beta'\in KK(\CC, C^{*}(\cG_{\RR}))$ is
  produced using the pseudodifferential calculus for groupoids (see
  below) and can be presented as a family $ \beta'=(\beta_{t})$ with: 
 $$  \beta_{0}=\beta; \ t>0,
 \beta_{t}=\left(C^{*}(\RR\times\RR,\frac{dx}{t}), 
   1,\frac{1}{\sqrt{1+x^{2}+t^{2}\partial_{x}^{2}}}\begin{pmatrix} 
   0 & x-t\partial_{x} \\ x+t\partial_{x} & 0
 \end{pmatrix}  \right) $$ 
After restricting at $t=1$ and applying the Morita equivalence;
only the index of the Fredholm operator appearing in
$\beta_{1}$ remains,  that is $+1$,  and this proves $\beta\otimes\delta=1$.  

Observe that by uniqueness of the inverse, we conclude that
$\delta=\alpha$ in $KK(C_{0}(\RR^{2}), \CC)$.  
\end{examples}
\begin{examples}(Boundary homomorphisms in long exact sequences)
Let 
 $$ 0\to I\underset{i}{\to}A\underset{p}{\to}B\to 0 $$
be a short exact sequence of $C^{*}$-algebras.  We assume that either
it admits a completely positive, norm
decreasing linear section or $I, A, B$ are $K$-nuclear
(\cite{ska-88}).  
Let $C_{p}=\{(a, \varphi)\in A\oplus C_{0}([0, 1[, B) \ |
p(a)=\varphi(0)\}$ be the cone of the homomorphism $p:A\to B$ and
denote by $d$ the homomorphism: $C_{0}(]0, 1[, B)\hookrightarrow
C_{p}$ given by $d(\varphi)=(0, \varphi)$ and by $e$ the
homomorphism: $I\to C_{p}$ given by $e(a)=(a, 0)$. Thanks to the
hypotheses, $[e]$ is invertible in $KK$-theory. 
One can set $\delta=[d]\otimes[e]^{-1}\in KK(C_{0}(\RR)\otimes B,
I)$ and using the Bott periodicity
$C_{0}(\RR^{2})\underset{\hbox{\tiny{KK}}}{\sim}\CC$ in order to identify:
$$  KK_2(C,D)=KK(C_{0}(\RR^{2})\otimes C, D)\simeq KK(C, D), $$
the connecting maps in the long exact sequences:
$$\cdots\to KK_{1}(I, D)\to KK(B, D)\overset{i^{*}}{\to}KK(A,
D)\overset{p^{*}}{\to}KK(I, D)\to KK_{1}(B, D)\to \cdots, $$
$$\cdots\to KK_{1}(C, B)\to KK(C, I)\overset{i_{*}}{\to}KK(C,
A)\overset{p_{*}}{\to}KK(C, B)\to KK_{1}(C, I)\to \cdots$$
are given by the appropriate Kasparov products with $\delta$.

\end{examples}

\newpage \part*{INDEX THEOREMS}

\section{Introduction to pseudodifferential operators on groupoids}\label{pseudo}
The historical motivation for developing  pseudodifferential
calculus on groupoids comes from
A. Connes, who implicitly introduced  this notion for foliations. Later
on, this calculus was axiomatized and studied on general groupoids by
several authors \cite{MP,NWX,Vas}. 

The following example illustrates how  pseudodifferential calculus
on groupoids arises in our approach of index theory. If $P$ is a partial differential operator
on $\RR^n$:
 $$ 
   P(x,D)= \sum_{|\alpha|\le d} c_\alpha(x)D_x^\alpha 
 $$
we may associate to it the following asymptotic operator:
 $$ P(x,tD)=\sum_{|\alpha|\le d} c_\alpha(x)(tD_x)^\alpha $$
by introducing a parameter $t\in]0,1]$ in front of each
$\partial_{x_j}$.  Here we use the usual convention:
$D_x^\alpha=(-i\partial_{x_{1}})^{\alpha_{1}}\ldots(-i\partial_{x_{n}})^{\alpha_{n}}
$.  We would like
to give a (interesting) meaning to the limit $t\to 0$. Of course we
would not be very happy with $tD\to 0$. 

To investigate this question, let us look at $P(x,tD)$ as a left multiplier on
$C^\infty(\RR^n\times\RR^n\times]0, 1])$ rather as a linear operator on 
$C^\infty(\RR^n)$:
\begin{eqnarray*}
 P(x,tD_x)u(x,y,t) &=& \int e^{(x-z).\xi} P(x,t\xi)u(z,y,t)dzd\xi \\
 &=& \int e^{\frac{x-z}{t}.\xi} P(x,\xi)u(z,y,t)\frac{dzd\xi}{t^n} \\
 &=& \int e^{(X-Z).\xi} P(x,\xi)u(x-t(X-Z),x-tX,t)dZd\xi. 
\end{eqnarray*}
In the last line we introduced the notation $X=\frac{x-y}{t}$ and performed
the change of variables $Z=\frac{z-y}{t}$.  

At this point, assume that $u$ has the following behaviour near $t=0$: 
 $$ 
   u(x,y,t) = \widetilde{u}(y,\frac{x-y}{t},t) \hbox{ where }\widetilde{u}\in
   C^\infty(\RR^{2n}\times[0,1]). 
 $$
It follows that: 
\begin{eqnarray*}
  P(x,tD_x)u(x,x-tX, t) &=& \int e^{(X-Z).\xi}
  P(x,\xi)\widetilde{u}(x-tX,Z,t)dZd\xi\\
  &\overset{t\to 0}{\longrightarrow}& \int e^{(X-Z).\xi}
  P(x,\xi)\widetilde{u}(x,Z,0)dZd\xi\\
  &=& P(x,D_X)\widetilde{u}(x,X,0). 
\end{eqnarray*}
Observations

\begin{itemize}

 \item  $P(x,D_X)$ is a partial differential operator in the
variable $X$ with constant coefficients, depending smoothly on a
parameter $x$ and with symbol coinciding with
the one of $P(x,D_x)$ in the sense that: 
 $\sigma(P(x,D_X)(x,X,\xi)=P(x,\xi)$. In particular, $P(x,D_X)$ is
 invariant under the translation $X\mapsto X+X_0$. Of course, $P(x,D_X)$
 is nothing else, up to a Fourier transform in $X$, than the symbol
 $P(x,\xi)$ of $P(x,D_x)$. In other words, denoting by
 $S_X(T\RR^n)$ the space of smooth functions $f(x,X)$ rapidly
 decreasing in $X$ and by $\cF_X$ the Fourier transform with respect
 to the variable $X$, we have a commutative diagram:
$$
 \xymatrix{
   S_X(T\RR^n)\ar[r]^{P(x,D_X)}\ar[d]_{\cF_X} &  S_X(T\RR^n)\ar[d]_{\cF_X}\\
   S_\xi(T^*\RR^n)\ar[r]^{P(x,\xi)} & S_\xi(T^*\RR^n)}
$$
where  $P(x,D_X)$ acts as a left multiplier on the
convolution algebra $S_X(T\RR^n)$ and $P(x,\xi)$ acts as a left multiplier on the
function algebra $S_\xi(T^*\RR^n)$ (equipped with the pointwise
multiplication of functions).

 \item $u$ and $\widetilde{u}$ are related by the bijection: 
  $$ 
  \begin{matrix}
     \phi: & \RR^{2n}\times[0,1] & \longrightarrow & \cG_{\RR^{n}} \\
     &  (x,X,t) & \longmapsto & (x-tX,x,t) \text{ if } t>0 \\
     &  (x,X,0) & \longmapsto &  (x,X,0)
  \end{matrix}
  $$
($\phi^{-1}(x,y,t)=(y,(x-y)/t,t),\ \phi^{-1}(x,X,0)=(x,X,0)$). 
In fact, the smooth structure of the tangent groupoid $\cG_{\RR^n}$
of the manifold $\RR^{n}$ (see Paragraph \ref{groupoid-index-theory})
is defined by requiring that  $\phi$ is a diffeomorphism. Thus 
$\widetilde{u}\in C^\infty(\RR^{2n}\times[0,1])$ means 
$u\in C^\infty(\cG_{\RR^n})$. 
\end{itemize}

Thus $P(x,D_X)$ is another way to look
at, and even, another way to {\bf define}, the symbol of $P(x,D_x)$. 
What is important for us is that it arises as a ``limit'' of a
family $P_t$ constructed with $P$, and the pseudodifferential calculus on the tangent groupoid
of $\RR^n$ will enable us to give a rigorous meaning to this
limit and perform interesting computations.

The material below is taken from \cite{MP,NWX,Vas}. 
Let $G$ be a Lie groupoid, with unit space $G^{(0)}=V$ and with a
smooth (right) Haar system $d\lambda$. We assume
that $V$ is a compact manifold and that the $s$-fibers $G_x$, 
$x\in V$, have no boundary. We denote by $U_\gamma$ the map induced on
functions by right multiplication by $\gamma$, that is: 
 $$ U_\gamma: C^\infty(G_{s(\gamma)})\longrightarrow
 C^\infty(G_{r(\gamma)}); \ U_\gamma f(\gamma')=f(\gamma'\gamma).$$
\begin{definition}
   A $G$-operator is a continuous linear map 
 $P:C^\infty_c(G)\longrightarrow C^\infty(G)$ such that: 
 \begin{itemize}
 \item[(i)] $P$ is given by a family $(P_x)_{x\in V}$ of linear
   operators $P_x: C^\infty_c(G_x)\to C^\infty(G_x)$ and: 
   $$ \forall f\in C^\infty_c(G),\quad P(f)(\gamma)=
   P_{s(\gamma)}f_{s(\gamma)}(\gamma) $$
where $f_x$ stands for the restriction $f|_{G_x}$. 
\item[(ii)] The following invariance property holds:
   $$ U_{\gamma}P_{s(\gamma)} =P_{r(\gamma)}U_{\gamma}.  $$
 \end{itemize}
 \end{definition}

Let $P$ be a $G$-operator and denote by $k_x\in C^{-\infty}(G_x\times
G_x)$ the Schwartz kernel of $P_x$, for each $x\in V$, as obtained
from the Schwartz kernel theorem applied to the manifold $G_x$
provided with the measure $d\lambda_x$. 

Thus, using the property [i]: 
$$ 
 \forall \gamma\in G, f\in C^\infty(G), \quad Pf(\gamma)=\int_{G_x}
 k_x(\gamma,\gamma')f(\gamma')d\lambda_x(\gamma'), \qquad
 (x=s(\gamma)). 
$$
Next:
$$ 
 U_\gamma Pf(\gamma')=Pf(\gamma'\gamma)=\int_{G_x}
 k_x(\gamma'\gamma,\gamma'')f(\gamma'')d\lambda_x(\gamma''), \qquad
 (x=s(\gamma)), 
$$
and 
\begin{eqnarray*}
P(U_\gamma f)(\gamma') &=& \int_{G_y}
 k_y(\gamma',\gamma'')f(\gamma''\gamma)d\lambda_y(\gamma''), \qquad
 (y=r(\gamma))\\
 &\overset{\eta=\gamma''\gamma}{=}& \int_{G_x}
 k_y(\gamma',\eta\gamma^{-1})f(\eta)d\lambda_x(\eta), \qquad
 (x=s(\gamma))\\
\end{eqnarray*}
where the last line uses the invariance property of Haar
systems. Axiom [ii] is equivalent to the following equalities of
distributions  
on $G_x\times G_x$, for all $x\in V$:
 $$ 
  \forall \gamma\in G,\quad  
     k_x(\gamma'\gamma,\gamma'')=k_y(\gamma',\gamma''\gamma^{-1})
     \quad (x=s(\gamma),\ y=r(\gamma) ). 
 $$
Setting  $k_P(\gamma):=k_{s(\gamma)}(\gamma,s(\gamma))$, we get 
$ k_x(\gamma,\gamma')=k_P(\gamma\gamma'^{-1})$, and the linear
operator $P: C^\infty_c(G)\to C^\infty(G)$ is given by:
$$ 
  P(f)(\gamma)=
  \int_{G_x}k_P(\gamma\gamma'^{-1})d\lambda_x(\gamma')\quad
  (x=s(\gamma)). 
$$
We may consider $k_P$ as a single distribution on $G$ acting on
smooth functions on $G$ by convolution. With a slight abuse of
terminology, we will refer to $k_P$ as the Schwartz (or convolution)
kernel of $P$. 

We say that $P$ is {\sl smoothing} if $k_P$ lies in $C^\infty(G)$ and is
{\sl compactly supported} or {\sl uniformly supported} if $k_P$ is
compactly supported (which implies that each $P_x$ is properly
supported). 

Let us develop some examples of $G$-operators. 
\begin{examples}
  \begin{enumerate}
  \item if $G=G^{(0)}=V$ is just a set, then $G_x=\{x\}$ for all $x\in V$. The property [i] is empty and the property
     [ii] implies that  a $G$-operator is
    given by pointwise multiplication by a smooth function $P\in C^\infty(V)$:
    $Pf(x)=P(x).f(x)$. 
  \item $G= V\times V$ the pair groupoid, and the Haar system
    $d\lambda$ is given in the obvious way by a single measure $dy$ on
    $V$: 
     $$ d\lambda_x(y) = dy \text{ under the identification }
     G_x=V\times\{x\}\simeq V
     $$
    It follows that for any $G$-operator $P$: 
 $$ Pg(z,x)= \int_{V\times\{x\}}k_P(z,y)g(y,x)d\lambda_{x}(y,x) =
 \int_{V}k_P(z,y)g(y,x)dy
 $$ 
 which immediately proves that $P_x=P_y$ are equal as linear operators
 on $C^\infty(V)$ under the obvious identifications 
$V\simeq V\times\{x\}\simeq V\times \{y\}$. 
\item Let $p: X\to Z$ a submersion, and 
$G=X\underset{Z}{\times}X = \{(x,y)\in X\times X\ | p(x)=p(y)\}$ the
associated subgroupoid of the pair groupoid $X\times X$. The
manifold $G_x$ can be identified with the fiber $p^{-1}(p(x))$.
Axiom [ii] implies that for any $G$-operator $P$, we have 
$P_x=P_y$ as linear operators on $ p^{-1}(p(x))$ as soon as
$y\in p^{-1}(p(x))$. Thus, $P$ is actually given by a family 
$\tilde{P}_z$, $z\in Z$ of operators on $p^{-1}(z)$, with 
the relation $P_x=\tilde{P}_{p(x)}$.
\item Let $G=E$ be the total space of a (euclidean, hermitian) vector
  bundle $p: E\to V$, with $r=s=p$. The Haar system $d_xw$, $x\in V$, is given by the
  metric structure on the fibers of $E$. We have here: 
 $$ 
   Pf(v) = \int_{E_x} k_P(v-w)f(w)d_xw \qquad (x=p(v))
 $$
Thus, for all $x\in V$, $P_x$ is a convolution operator on
the linear space $E_x$.
\item Let $G=\cG_V=TV\times\{0\}\sqcup V\times V\times ]0,1]$ be the
  tangent groupoid of $V$. It can be viewed as a family of groupoids
  $G_t$ parametrized
  by $[0,1]$, where $G_0=TV$ and $G_t=V\times V$
  for $t>0$. A $\cG_V$-operator is given by a family $P_t$ of
  $G_t$-operators, and $(P_t)_{t>0}$ is a family  of operators on
  $C^\infty_c(V)$ parametrized by $t$
  while $P_0$ is a family  of translation invariant operators on
  $T_xV$ parametrized by $x\in V$. 
  The $\cG_V$-operators are thus a blend of Examples 2 and 4. 
  \end{enumerate}
\end{examples}

We now turn to the definition of pseudodifferential operators on a Lie
groupoid $G$.
\begin{definition}
  A $G$-operator $P$ is a $G$-pseudodifferential operator
  of order $m$ if: 
  \begin{enumerate}
  \item The Schwartz kernel $k_P$ is smooth outside $G^{(0)}$. 
  \item For every distinguished chart $\psi: U\subset G\to \Omega\times s(U)\subset\RR^{n-p}\times\RR^{p}$ of $G$:
 $$  \xymatrix{U\ar[rr]^{\psi}\ar[dr]_{s} & &\Omega\times s(U)\ar[dl]^{p_{2}} \\
        &s(U)&  }$$ 
 the operator $(\psi^{-1})^{*}P\psi^{*}:
 C^{\infty}_{c}(\Omega\times s(U))\to C^{\infty}_{c}(\Omega\times
 s(U))$ is a smooth family parametrized by $s(U)$ of  
 pseudodifferential operators of order $m$ on $\Omega$.     
  \end{enumerate}
\end{definition}
We will use very few properties of this calculus and only provide some examples and a list of properties. The reader can
find a complete presentation in \cite{Vas,Vas-these,NWX,MP,Mo2003}.

\begin{examples}
  In the previous five examples,  a
  $G$-pseudodifferential operator is: 
  \begin{enumerate}
  \item  an operator given by
    pointwise multiplication by a smooth function on $V$; 
  \item  a single pseudodifferential operator on $V$;
  \item a smooth family parametrized by $Z$ of pseudodifferential
    operators in the fibers: this coincides with the notion of \cite{AS4};
  \item a family parametrized by $x\in V$ of convolution operators in
    $E_{x}$ such that the underlying distribution $k_{P}$ identifies with
    the Fourier transform of a symbol
    on $E$ (that is,  a smooth function on $E$ satisfying the
    standard decay conditions with respect to its variable in the
    fibers);
  \item the data provided by an asymptotic pseudodifferential
    operator on $V$ together with its complete symbol,  the choice
    of it depending on the gluing in $\cG_{V}$: this is quite close
    to the notions studied in \cite{ge,BF,ENN}. 
  \end{enumerate}
\end{examples}
It turns out that the space $\Psi_c^*(G)$ of
compactly supported $G$-pseudodifferential operators is an involutive algebra.  

The principal symbol of a $G$-pseudodifferential operator $P$ of
order $m$ is defined as a function $\sigma_{m}(P)$ on $A^{*}(G)\setminus G^{(0)}$ by: 
 $$ \sigma_{m}(P)(x, \xi) = \sigma_{pr}(P_{x})(x, \xi)$$
where $\sigma_{pr}(P_{x})$ is the principal symbol of the pseudodifferential
operator $P_{x}$ on the manifold $G_{x}$.  
Conversely,  given a symbol $f$ of order $m$ on $A^{*}(G)$ together with the
following data:
\begin{enumerate}
\item A smooth embedding  $\theta: \cU\to AG$,  where $\cU$ is a open
  set in $G$ containing $G^{(0)}$, such that
  $\theta(G^{(0)})=G^{(0)}$, $(d\theta)|_{G^{0}}=\hbox{Id}$ and $\theta(\gamma)\in A_{s(\gamma)}G$
  for all $\gamma\in \cU$;
\item A smooth compactly supported map $\phi:G\to \RR_{+}$ such that
  $\phi^{-1}(1)=G^{(0)}$;
\end{enumerate}
we get a $G$-pseudodifferential operator $P_{f, \theta, \phi}$ by the
formula:
$$u\in C^{\infty}_{c}(G),\ P_{f, \theta, \phi}u(\gamma)=
 \int_{\gamma'\in G_{s(\gamma)},\atop \xi\in A^{*}_{r(\gamma)}(G)}
 e^{-i\theta(\gamma'\gamma^{-1}). \xi}f(r(\gamma),
 \xi)\phi(\gamma'\gamma^{-1})u(\gamma')d\lambda_{s(\gamma)}(\gamma')$$ 
The principal symbol of $P_{f, \theta, \phi}$ is just the leading
part of $f$.  

The principal symbol map respects pointwise product while the product
law for total symbols is much more involved. An operator is {\it
  elliptic} when its principal symbol never vanishes and in that
case,  as in the classical situation, 
it has a parametrix inverting it modulo
$\Psi^{-\infty}_{c}(G)=C^{\infty}_{c}(G)$. 

Operators of negative order in $\Psi_c^*(G)$ are actually in
$C^{*}(G)$,  while zero order operators are in the multiplier algebra
$\cM(C^{*}(G))$. 
  
All these definitions and properties immediately extend to the case of operators acting
between sections of bundles on $G^{(0)}$ pulled back to $G$ with the
range map $r$. The space of compactly supported pseudodifferential
operators on  $G$ acting on sections of $r^*E$ and taking values in sections of $r^*F$
will be noted $\Psi_c^*(G,E, F)$. If $F=E$ we get an algebra
denoted by $\Psi_c^*(G,E)$. 

\begin{examples}\label{symbols-and-tangentgroupoid}
  \begin{enumerate}
  \item The family given by $P_{t}=P(x, tD_{x})$ for $t>0$ and
    $P_{0}=P(x, D_{X})$ described in the introduction of this section
    is a $G$-pseudodifferential operator with $G$ the tangent
    groupoid of $\RR^{n}$. 
  \item More generally,  let $V$ be a closed manifold endowed with a
    riemannian metric. We note $\exp$ the exponential map associated
    with the metric.  Let $f$ be a
    symbol on $V$.  We get a  $\cG_{V}$-pseudodifferential operator
    $P$ by setting:
 $$(t>0)\quad P_{t}u(x, y, t)=\int_{z\in V, \xi\in T^{*}_{x}V}
 e^{\frac{\exp^{-1}_{x}(z)}{t}.\xi}f(x,\xi)u(z,y)\frac{dzd\xi}{t^n}$$
 $$P_{0}u(x, X, 0)=\int_{Z\in T_{x}V, \xi\in
   T^{*}_{x}V}e^{(X-Z).\xi}f(x,\xi)u(x,Z)dZd\xi $$
Moreover,  $P_{1}$ is a
pseudodifferential operator on the manifold $V$ which admits $f$ as a
complete symbol.  
\end{enumerate}  
\end{examples}

\section{Index theorem for smooth manifolds}\label{debutInd}

The purpose of this last lecture is to present a proof of the
Atiyah-Singer index theorem using deformation groupoids 
and show how it generalizes to conical pseudomanifolds. The results
presented here come from recent works of the authors together with a
joint work with V. Nistor \cite{DL2,DLN,Nous2}, we refer to
\cite{DL2,DLN} for the proofs.

\medskip \noindent{\bf The $KK$-element associated to a deformation groupoid}

\noindent Before going to the description of the index maps, let us
describe a useful and classical construction \cite{Co0,HS}. \\ 
Let $G$ be a smooth deformation groupoid (definition \ref{def-deformation-groupoid}): 
$$G= G_1 \times \{0\} \cup G_2\times ]0,1] 
\rightrightarrows G^{(0)}=M\times [0,1].$$ 
One can consider the saturated 
open subset $M\times ]0,1]$ of $G^{(0)}$. Using the isomorphisms 
$C^*(G\vert_{M\times ]0,1]}) \simeq C^*(G_2)\ot C_0(]0,1])$ and 
$C^*(G\vert_{M\times\{0\}})\simeq C^*(G_1)$,  we obtain the following 
exact sequence of $C^*$-algebras: 
$$\begin{CD} 0 @>>> C^*(G_2)\ot C_0(]0,1]) @>{i_{M\times]0,1]}}>> C^*(G) 
@>{ev_0}>>  C^*(G_1)  @>>> 0 \end{CD} 
$$ 
where $i_{M\times]0,1]}$ is the inclusion map and $ev_0$ is the {\it evaluation map} at 
$0$, that is $ev_0$ is the map coming from the 
restriction of functions to $G\vert_{M\times\{0\}}$.

\smallskip \noindent We assume now that $C^*(G_1)$ is nuclear. Since 
the $C^*$-algebra $C^*(G_2)\ot C_0(]0,1])$ is contractible, 
 the long exact sequence in $KK$-theory shows that the group homomorphism 
$(ev_0)_*=\cdot {\ot}[ev_0]:KK(A,C^*(G)) \rightarrow KK(A, C^*(G_1))$ 
is an isomorphism for each $C^*$-algebra $A$. 
 
\noindent In particular with $A=C^*(G)$
we get that $[ev_0]$ is invertible in 
$KK$-theory: there is an element $[ev_0]^{-1}$ in 
$KK(C^*(G_1), C^*(G))$ such that $[ev_0] {\ot} [ev_0]^{-1}=1_{C^*(G)}$ 
and $[ev_0]^{-1} {\ot} [ev_0]=1_{C^*(G_1)}$. 
 
\smallskip \noindent Let $ev_1:C^*(G) \rightarrow C^*(G_2)$ be the 
evaluation map at $1$ and $[ev_1]$ the corresponding element of 
$KK(C^*(G),C^*(G_2))$. 
 
\smallskip \noindent The {\it $KK$-element associated to the 
  deformation groupoid} $G$ is defined by: 
$$\delta=[ev_0]^{-1} {\ot} [ev_1]\in  KK(C^*(G_1),C^*(G_2)) \ . $$

\noindent We will meet several examples of this construction in the sequel.

\medskip \noindent {\bf The analytical index}

\noindent Let $M$ be a closed manifold and consider its tangent groupoid:
$$\cG^t_M:=TM\times \{0\} \cup M\times M \times ]0,1]
\rightrightarrows M\times [0,1]$$
It is a deformation groupoid and the construction above provides us a
$KK$-element:
$$\partial_M=(e^M_1)_*\circ (e^M_0)_*^{-1}\in KK(C_0(T^*M),\cK)\simeq KK(C_0(T^*M),\CC), $$
where $e^{M}_{i}: C^*(\cG^t_M)\to C^*(\cG^t_M)|_{t=i} $ are evaluation
homomorphisms. 

\noindent The analytical index is then \cite{Co0} $$\begin{array}{cccc} Inda_M:= (e^M_1)_*\circ (e^M_0)_*^{-1}:&
KK(\CC,C_0(T^*M)) &
  \rightarrow  & KK(\CC,\cK(L^2(M)) \\ & \simeq
  K_0(C_0(T^*M)) & & \simeq \ZZ \end{array} $$ or in terms of Kasparov
  product $$Inda_M= \cdot \otimes \partial _M \ .$$

\noindent Using the notion of pseudodifferential calculus for
$\cG_{M}^{t}$, it is easy to justify that this map is the usual
analytical index map.  Indeed,  let $f(x,\xi)$ be an
elliptic zero order symbol and consider the $\cG_{M}^{t}$-pseudodifferential
 operator, $P_{f}=(P_{t})_{0\le t\le 1}$,  defined as in Example \ref{symbols-and-tangentgroupoid}. 
Then $f$ provides a $K$-theory class $[f]\in K_{0}(C^*(TM))\simeq
K_0(C_0(T^*M))$ while $P$ provides a $K$-theory class $[P]\in
K_0(C^*(\cG^t_M))$ and:
  $$  (e^M_0)_*([P])=[f]\in K_0(C^*(TM)) $$

\noindent Thus: 
 $$ [f]\otimes [e^M_0]^{-1} \otimes [e^M_1] = [P_1]\in K_0(\cK) $$
and $[P_1]$ coincides with $\mbox{Ind}(P_1)$ under 
$K_0(\cK)\simeq \ZZ$. 

\noindent Since $P_1$ has principal symbol equal to the leading part of
$f$, and since every class in $K_0(C_0(T^*M))$ can be obtained from
a zero order elliptic symbol,  the claim is justified.

\medskip \noindent To be complete,  let us explain that the 
analytical index map is the Poincar\'e dual of the homomorphism in
$K$-homology associated with the obvious map: $M\to\{\cdot\}$.  Indeed, thanks to the obvious homomorphism
$\Psi:C^*(TM)\otimes C(M)\rightarrow C^*(TM)$ given by multiplication, $\partial_M$  can be lifted into an element  
$D_M=\Psi_*(\partial_M)\in KK(C^*(TM)\otimes
C(M),\CC)=K^0(C^*(TM)\otimes C(M))$, called the {\it Dirac
  element}. This Dirac element yields the well known Poincar\'e
duality between 
$C_{0}(T^{*}M)$ and $C(M)$ (\cite{CoS,Ka2, DL2}),  and in particular
it gives an isomorphism: 
 $$
\cdot\underset{C^*(TM)}{\otimes}D_M: \  K_0(C^*(TM))\overset{\simeq
}{\longrightarrow}K^0(C(M)) 
 $$
whose inverse is induced by the principal
symbol map.  

\noindent One can then easily show the following proposition: 
\begin{proposition}\label{PD-index-smooth} 
Let $q:M\to \cdot$ be the projection onto a point.  The following
diagram commutes:
$$ 
\begin{CD}
   K^0(T^*M) @>\mathrm{PD}>> K_0(M) \\
    @V\mbox{Ind}_aVV  @VVq_*V  \\
     \ZZ @>=>> \ZZ
\end{CD}
$$
\end{proposition}

\medskip \noindent{\bf The topological index}

\noindent Take an embedding $M \rightarrow \RR^n$, and let $p:N \rightarrow M$
be the normal bundle of this embedding. The vector bundle $TN \rightarrow TM$ admits a complex
structure, thus we have a Thom isomorphism:
$$ T:  K_0(C^*(TM))\overset{\simeq}{\longrightarrow}
K_0(C^*(TN))$$
given by a $KK$-equivalence:
$$T\in KK(C^*(TM),C^*(TN))\ .$$ $T$ is called the {\it Thom\ element} \cite{Ka1}.

\smallskip \noindent The bundle $N$ identifies with an open neighborhood of $M$
into $\RR^n$, so we have the excision map:  
 $$j: C^*(TN) \rightarrow C^*(T\RR^n).$$
Consider also: $\ds B: K_0(C^*(T\RR^n))\to \ZZ$
given by the isomorphism $C^*(T\RR^n)\simeq C_0(\RR^{2n})$ together with Bott periodicity. 

\smallskip \noindent The {\it topological index map} $\mbox{Ind}_{t}$
is the composition:  
$$\hspace{-1.2cm}K(C^*(TM))
  \overset{ T }{{\longrightarrow}} K(C^*(TN)) 
\overset{j_*}{\longrightarrow} K(C^*(T\RR^n))
\overset{B}{\underset{\simeq}{\longrightarrow}} \ZZ. $$ 

\medskip \noindent This classical construction can be reformulated with
groupoids. 

\smallskip \noindent 
First, let us give a description of $T$, or rather of its inverse, in terms
of groupoids. Recall the construction of the Thom groupoid.
We begin by pulling back $TM$ over $N$ in the groupoid sense:
$$\text{Let}:\qquad \pb{p}(TM) =
N\underset{M}{\times}TM\underset{M}{\times}N\rightrightarrows N.$$
$$\text{Let}:\qquad  \cT_N = TN\times\{0\} \sqcup\pb{p}(TM) \times ]0,1]
\rightrightarrows N\times [0,1].$$ 
This {\sl Thom groupoid} and the Morita equivalence between $\pb{p}(TM)$
and $TM$ provides the $KK$-element:
 $$\tau_{N}\in KK(C^*(TN),C^*(TM)) \ .$$ 
This element is defined exactly as $\partial_M$ is. Precisely, the evaluation
map at $0$, $\tilde{e}_0: C^*(\cT_N)\rightarrow C^*(TN)$ defines an
invertible $KK$-element. We let $\tilde{e}_1: C^*(\cT_N)\rightarrow
C^*(\pb{p}(TM))$ be the evaluation map at $1$. The Morita equivalence
between the groupoids $TM$ and $\pb{p}(TM)$ leads to a Morita
equivalence between the corresponding $C^*$-algebra and thus to a
$KK$-equivalence $\cM \in KK(C^*(\pb{p}(TM)),C^*(TM))$. Then $$\tau_N:=
[\tilde{e}_0]^{-1} \otimes [\tilde{e}_1] \otimes \cM \ .$$

\noindent We have the following:
 \begin{proposition}\cite{DLN} If $T$ is the $KK$-equivalence giving the Thom
   isomorphism then: $$\tau_{N}=T^{-1}.$$
 \end{proposition}

\noindent This proposition also applies to interpret the isomorphism
$B: K_0(C^*(T\RR^n))\to 
\ZZ \ .$  

\smallskip \noindent Indeed, consider the embedding $\cdot\hookrightarrow\RR^n$. The
normal bundle is just $\RR^n\to \cdot$ and we get as before: 
   $$\tau_{\RR^n}\in KK(C^*(T\RR^n),\CC)$$
Using the previous proposition we get: $\ds B = \cdot\otimes \tau_{\RR^n}$.

\smallskip \noindent Remark also that $\cT_{\RR^n}=\cG_{\RR^n}$ so
that $\tau_{\RR^n}=[e_{0}^{\RR^{n}}]^{-1}\otimes[e_{1}^{\RR^{n}}]$. 

\medskip \noindent Finally the topological index: 
  $$ \mbox{Ind}_{t} =  \tau_{\RR^n} \circ j_*\circ\tau_{N}^{-1} $$
is entirely described using (deformation) groupoids.

\bigskip \noindent{\bf The equality of the indices} 

A last groupoid is necessary in order to prove the equality of index
maps. Namely, this groupoid is obtained by recasting the construction
of the Thom groupoid at the level of tangent groupoids: 
\begin{equation}
  \label{super-thom-groupoid}
  \widetilde{\cT}_{N}=\cG_{N}\times\{0\}\sqcup\pb{(p\otimes\hbox{Id}_{[0, 1]})}(\cG_{M})\times]0, 1]
\end{equation}
As before,  this yields a class:
 $$
  \widetilde{\tau}_{N}\in KK(C^{*}(\cG_{N}), C^{*}(\cG_{M})).  
 $$
\medskip \noindent All maps in the following diagram: 
\begin{equation}
  \label{AS-diagram}
\hspace{-0.5cm}\begin{CD}
  \ZZ  @= \ZZ @= \ZZ \\
   @Ae^{M}_{1}AA @Ae^{N}_{1}AA @Ae^{\RR^{n}}_{1}AA \\
    K_0(C^*(\cG_M))@<\otimes\widetilde{\tau}_{N}<<K_0(C^*(\cG_N))@>\widetilde{j}_{*}>>K_0(C^*(\cG_{\RR^n}))\\  
   @Ve^{M}_{0}V\simeq V @Ve^{N}_{0}V\simeq V @Ve^{\RR^{n}}_{0}V\simeq V \\
    K_0(C^*(TM))@<\otimes\tau_{N}<\simeq <K_0(C^*(TN))@>j_{*}>>K_0(C^*(T\RR^n))
\end{CD}
\end{equation}  
are given by Kasparov products with:  
\begin{enumerate}
 \item classes of homomorphisms coming from restrictions or inclusions 
   between groupoids,  
 \item  inverses of such classes, 
 \item explicit Morita equivalences. 
\end{enumerate} 
This easily yields the commutativity of diagram (\ref{AS-diagram}). 
Having in mind the previous description of index maps using
groupoids, this commutativity property just implies: 
  $$\mbox{Ind}_{a}=\mbox{Ind}_{t}$$

\section{The case of  pseudomanifolds with isolated singularities}
As we explained earlier,  the proof of the $K$-theoretical
form of the Atiyah-Singer presented in these lectures extends very easily to the case of
pseudomanifolds with isolated singularities. This is achieved
provided one uses the correct notion of {\sl tangent space} of the 
pseudomanifold and for a pseudomanifold $X$ with
one conical point (the case of several isolated singularities is
similar),  this is the noncommutative tangent space defined in
section \ref{ncgtangent}: 
$$ T^{\fS}X=X^-\times X^- \cup T\overline{X^+} \rightrightarrows \smX$$
In the sequel,  it will replace the ordinary tangent space of a smooth
manifold. Moreover, it gives rise to another deformation groupoid which
will replace the ordinary tangent groupoid of a smooth manifold:
$$\cG^t_X=T^{\fS}X\times \{0\} \cup \smX \times \smX \times
]0,1] \rightrightarrows \smX \times [0,1]$$
We call $\cG^t_X$ the {\sl tangent groupoid} of $X$. It can be
provided with a smooth structure such that $T^{\fS}X$ is a smooth
subgroupoid.  Moreover both are amenable so their  reduced and
maximal $C^*$-algebras coincide and are nuclear. 

With these choices of $T^{\fS}X$ as a tangent space for $X$ and of
$\cG^t_X$ as a tangent groupoid,  one can follow step by step all the
constructions made in the previous section. 

\subsection{The analytical index}
Using the partition $\smX\times [0,1]=\smX\times \{0\} \cup \smX
\times ]0,1]$ into saturated open and closed subsets of the units
space of the tangent groupoid, we define the $KK$-element associated
to the tangent groupoid of $X$:

$$\partial_X:= [e_0]^{-1}\otimes [e_1]\in KK(C^*(T^{\fS}X),\cK)\simeq
KK(C^*(T^{\fS}X),\CC)\ ,$$

\noindent where $e_0:C^*(\cG^t_X)\rightarrow C^*(\cG^t_X\vert_{\smX\times
  \{0\}})\simeq C^*(T^{\fS}X)$ is the evaluation at $0$ and
$e_1:C^*(\cG^t_X)\rightarrow C^*(\cG^t_X\vert_{\smX\times 
  \{1\}})\simeq \cK(L^2(X))$ is  the evaluation at $1$. 

 Now we can define the analytical index exactly as we did for closed smooth
 manifolds. Precisely,  the {\it analytical index} for $X$ is set to
 be the map: $$\mbox{Ind}_{a}^{X}=\cdot \otimes \partial_X:
 KK(\CC,C^*(T^{\fS}X)) \rightarrow KK(\CC,\cK(L^2(\smX)))\simeq \ZZ\
 .$$
The interpretation of this map as the Fredholm index of an
appropriate class of elliptic operators is possible and carried out in \cite{L-cone2}.

\subsection{The Poincar\'e duality}
Pursuing the analogy with smooth manifolds, we explain in this
paragraph that the
analytical index map for $X$ is Poincar\'e dual to the index map in
$K$-homology associated to the obvious map : $X\to \{. \}$. 

\smallskip \noindent The algebras $C(X)$ and $C^{\bullet}(X):=\{f\in C(X) \ \vert
\ f \mbox{ is constant on } cL\}$ are isomorphic. If $g$ belongs to
$C^{\bullet}(X)$ and $f$ to $C_c(T^{\fS}X)$, let $g\cdot f$ be the
element of $C_c(T^{\fS}X)$ defined by $g\cdot f
(\gamma)=g(r(\gamma))f(\gamma)$. This induces a *-morphism $$\Psi:
C(X)\otimes C^*(T^{\fS}X) \rightarrow C^*(T^{\fS}X)\ . $$
The {\it Dirac element} is defined to be  $$D_X:= [\Psi]\otimes \partial_X \in
KK(C(X)\otimes C^*(T^{\fS}X),\CC) .$$
We recall
\begin{theorem}\cite{DL2}  There exists a (dual-Dirac) element
  $\lambda_X \in KK(\CC,C(X)\otimes C^*(T^{\fS}X))$ such that
  $$\lambda_X \underset{C(X)}{\otimes} D_X = 1_{C^*(T^{\fS}X)} \in
  KK(C^*(T^{\fS}X),C^*(T^{\fS}X))\ , $$ $$ \lambda_X
  \underset{C^*(T^{\fS}X)}{\otimes} D_X = 1_{C( X)} \in KK(C(X),C(X))
  \ . $$ 
This means that $C(X)$ and $ C^*(T^{\fS}X)$ are Poincar\'e dual.    
\end{theorem}
\begin{remark}
 The explicit construction of $\lambda_{X}$,  which is heavy gooing and
   technical,  can be avoided. In fact,  the definitions of
   $T^{\fS}X$,  $\cG^t_X$  and thus of  $D_{X}$,  can be
   extended in a very natural way to the case of an arbitrary
   pseudomanifold and the proof of Poincar\'e duality can be done
   using a recursive argument on the depth of the stratification,
   starting with the case depth$=0$,  that is with the case of smooth
   closed manifolds.  This is the subject of \cite{Nous2}. 
\end{remark}

 The theorem implies that:
$$ \begin{array}{ccc} KK(\CC,C^*(T^{\fS}X))\simeq K_0(C^*(T^{\fS}X)) &
  \rightarrow & K(C(X),\CC)\simeq K^0(C(X)) \\ x & \mapsto & x
  \underset{C^*(T^{\fS}X)}{\otimes} D_X \end{array} $$
 is an isomorphism. In \cite{L-cone2}, it is explained how to in
 interpret its inverse as a principal symbol map,  and one also get
 the analogue of Proposition \ref{PD-index-smooth}:
\begin{proposition} 
Let $q:X\to \cdot$ be the projection onto a point.  The following
diagram commutes:
$$ 
\begin{CD}
   K_{0}(C^*(T^{\fS}X)) @>\mathrm{PD}>> K_0(X) \\
    @V\mbox{Ind}_a^{X}VV  @VVq_*V  \\
     \ZZ @>=>> \ZZ
\end{CD}
$$
\end{proposition}

\subsection{The topological index}

\ 

\noindent {\bf Thom isomorphism} Take an {\it embedding }  
$X \hookrightarrow c\RR^n=\RR^n\times[0,+\infty[/\RR^n\times
\{0\}$. This means that we have a map which restricts to an embedding
$\smx \rightarrow \RR^n\times]0,+\infty[$ in the usual sense and which
sends $c$ to the image of $\RR^n\times\{0\}$ in $c\RR^n$. Moreover we
require the embedding on $X^-=L\times ]0,1[$ to be of the form $j\times
\id$ where $j$ is an embedding of $L$ in $\RR^n$.

\smallskip \noindent Such an embedding provides a {\it conical normal
  bundle}. Precisely, let $p: N^{\circ} \rightarrow \smx$ be the normal bundle associated 
with $ \smx\hookrightarrow \RR^n\times]0,+\infty[$. We can
identify $N^{\circ}|_{X^-}\simeq N^{\circ}\vert_L\times ]0,1[$, 
and set :
 $$N=\bar cN^{\circ}|_L \cup N^{\circ}|_{X^+} \ .$$
Thus $N$ is the pseudomanifold with an isolated singularity obtained
by gluing the closed cone $\bar cN^{\circ}|_L:= N^{\circ}|_L\times
[0,1]/N^{\circ}|_L\times \{0\}$ with $ N^{\circ}|_{X^+}$ along their
common boundary $N^{\circ}|_L\times \{1\}= N^{\circ}|_{\partial
  X^+}$. Moreover $p: N\rightarrow X$ is a conical vector bundle.

\smallskip \noindent The {\it Thom groupoid} is then: 
$$
\cT_{N}=T^{\fS}N \times\{0\}\sqcup\pb{p}(T^{\fS}X)\times]0,1] \ .
$$ 
It is a deformation groupoid. The corresponding $KK$-element gives the {\it inverse Thom element}: 
$$\tau_N\in KK(C^*(T^{\fS}N),C^*(T^{\fS}X))\ .$$

\begin{proposition}\cite{DLN} The following map is an isomorphism.
$$ 
  K(C^*(T^{\fS}N)) \overset{\cdot \otimes \tau_N}{\longrightarrow}
  K(C^*(T^{\fS}X))
$$ 
\end{proposition}

Roughly speaking, the inverse of $\cdot \otimes \tau_N$ is the {\it Thom
isomorphism} for the  
``vector bundle'' $T^\fS N$ ``over'' $T^\fS X$. One can show that it really restricts to
usual Thom homomorphism on regular parts.

\medskip \noindent {\bf Excision } The groupoid $T^{\fS}N$ is identified with an open subgroupoid of
$T^{\fS}c\RR^n$ and we have an excision map:

 $$j:C^*(T^{\fS}N)\rightarrow C^*(T^{\fS} \RR^n)\ .$$
 
\medskip \noindent {\bf Bott element} Consider $c\hookrightarrow c\RR^n$.\\ The (conical) normal bundle is $c\RR^n$ itself. 
Remark that $\cG^t_{c\RR^n}=\cT_{c\RR^n}$. 
Then
  $$ \tau_{c\RR^n}\in KK(C^*(T^{\fS}c\RR^n),\CC) $$
gives an isomorphism:
 $$B=(\cdot\otimes\tau_{c\RR^n}): K_0(C^*(T^{\fS}c\RR^n))\to \ZZ$$

\begin{definition} The {\it topological index} is the
morphism $$\mbox{Ind}_{t}^{X}=B\circ j_*\circ \tau_{N}^{-1}:\
K_0(C^*(T^{\fS}X))\to \ZZ$$
\end{definition}

The following index theorem can be proved along the same lines as in
the smooth case: 
\begin{theorem} The following equality holds:
 $$\mbox{Ind}_{a}^{X}=\mbox{Ind}_{t}^{X}$$
\end{theorem}


\end{document}